\documentclass[11pt,a4paper]{article}
\usepackage{jheppub}
\usepackage{physics}
\usepackage{comment}
\usepackage{bm}
\usepackage{mathrsfs}  
\usepackage{amsmath}  
\usepackage{comment}  
\usepackage{pict2e}  
\usepackage{slashed}
\usepackage{color}
\usepackage[dvipsnames]{xcolor}
\usepackage{cancel}
\usepackage{soul}
\usepackage{framed}
\usepackage{parskip}
\usepackage{tikz-cd}
\usepackage{enumitem}
\usepackage{soul}

\definecolor{purple}{rgb}{0.6,0,0.9}
\definecolor{bl}{rgb}{0.3,0.1,0.9}
\definecolor{gr}{rgb}{0.133,0.545,0.133}
\definecolor{darkred}{rgb}{0.9,0,.2}
\newcommand{\bl}[1]{\textcolor{bl}{#1}}
\newcommand\ol{\overline}

\newcommand\Aut{\mathop{\mathrm{Aut}}}
\newcommand\disc{\mathop{\mathrm{disc}}}
\newcommand\Hom{\mathop{\mathrm{Hom}}}
\newcommand\rk{\mathop{\mathrm{rk}}}
\newcommand\SU{\mathop{\mathrm{SU}}}

\newcommand{\C}{\mathbb C}
\newcommand{\F}{\mathbb F}
\newcommand{\N}{\mathbb N}
\renewcommand{\P}{\mathbb P}
\newcommand{\R}{\mathbb R}
\newcommand{\Z}{\mathbb Z}
\newcommand{\Q}{\mathbb Q}

\newcommand{\CCC}{\mathcal C}
\newcommand{\eps}{\varepsilon}
\newcommand{\w}{\text{w}}

\newcommand\m{\mathcal}

\newcommand{\hf}{\frac{1}{2}}
\newcommand\wt{\widetilde}
\newcommand\wh{\widehat}

\newcommand{\ch}[1]{\stackrel{\star}{\raisebox{0pt}{$#1$}}}
\newcounter{dom}\newcommand{\domino}{\;\;\;\makebox[0pt]{$\refstepcounter{dom}\thedom$}}

\newcommand{\be}{\begin{equation}}
\newcommand{\ee}{\end{equation}}
\newcommand{\bea}{\begin{eqnarray}}
\newcommand{\eea}{\end{eqnarray}}

\makeatletter
\newcommand{\xleftrightarrow}[2][]{\ext@arrow 3359\leftrightarrowfill@{#1}{#2}}
\newcommand{\xdashrightarrow}[2][]{\ext@arrow 0359\rightarrowfill@@{#1}{#2}}
\newcommand{\xdashleftarrow}[2][]{\ext@arrow 3095\leftarrowfill@@{#1}{#2}}
\newcommand{\xdashleftrightarrow}[2][]{\ext@arrow 3359\leftrightarrowfill@@{#1}{#2}}
\def\rightarrowfill@@{\arrowfill@@\relax\relbar\rightarrow}
\def\leftarrowfill@@{\arrowfill@@\leftarrow\relbar\relax}
\def\leftrightarrowfill@@{\arrowfill@@\leftarrow\relbar\rightarrow}
\def\arrowfill@@#1#2#3#4{%
  $\m@th\thickmuskip0mu\medmuskip\thickmuskip\thinmuskip\thickmuskip
   \relax#4#1
   \xleaders\hbox{$#4#2$}\hfill
   #3$%
}
\makeatother

\newcommand{\xdownarrow}[1]{%
  {\left\downarrow\vbox to #1{}\right.\kern-\nulldelimiterspace}
}

\newtheorem{Claim}{Claim}[subsection]
\newcommand{\bc}{\begin{framed}\vspace*{-0.3em}\begin{Claim}\hspace*{\fill}\\}
\newcommand{\ec}{\end{Claim}\end{framed}}

\newtheorem{Corollary}[Claim]{Corollary}
\newcommand{\bcor}{\begin{framed}\vspace*{-0.3em}\begin{Corollary}\hspace*{\fill}\\}
\newcommand{\ecor}{\end{Corollary}\end{framed}}

\newtheorem{Construction}[Claim]{Construction}
\newcommand{\bconst}{\begin{framed}\vspace*{-0.3em}\begin{Construction}\hspace*{\fill}\\}
\newcommand{\econst}{\end{Construction}\end{framed}}

\newtheorem{Lemma}[Claim]{Lemma}
\newcommand{\blem}{\begin{framed}\vspace*{-0.3em}\begin{Lemma}\hspace*{\fill}\\}
\newcommand{\elem}{\end{Lemma}\end{framed}}

\newtheorem{Theorem}[Claim]{Theorem}
\newcommand{\bthm}{\begin{framed}\vspace*{-0.3em}\begin{Theorem}\hspace*{\fill}\\}
\newcommand{\ethm}{\end{Theorem}\end{framed}}

\newtheorem{Remark}[Claim]{Remark}
\newcommand{\brmk}{\begin{Remark}\hspace*{\fill}\\}
\newcommand{\ermk}{\end{Remark}}

\newtheorem{Example}[Claim]{Example}
\newcommand{\bexpl}{\begin{Example}\hspace*{\fill}\\}
\newcommand{\expl}{\end{Example}}

\newtheorem{Proposition}[Claim]{Proposition}
\newcommand{\bprop}{\begin{framed}\vspace*{-0.3em}\begin{Proposition}\hspace*{\fill}\\}
\newcommand{\eprop}{\end{Proposition}\end{framed}}

\newcommand{\bpf}{\textsl{Proof}:\hspace*{\fill}\\}
\newcommand{\epf}{\hspace*{\fill}$\Box$\\}

\author[a]{Kasia Budzik}
\author[b]{, Anne Taormina}
\author{, Mara Ungureanu}
\author[c]{, Katrin Wendland}
\author[d]{, Ida G. Zadeh}

\affiliation[a]{Center for the Fundamental Laws of Nature, Harvard University, Cambridge, MA 02138, USA} 
\affiliation[b]{Department of Mathematics, King's College London, Strand,\,London WC2R 2LS, United Kingdom}
\affiliation[c]{School of Mathematics, Trinity College Dublin, The University of Dublin, Dublin 2, D02 PN40,
Ireland}
\affiliation[d]{Mathematical Sciences and STAG Research Centre, University of Southampton, Highfield, Southampton SO17 1BJ, U.K.}

\emailAdd{kbudzik@fas.harvard.edu}
\emailAdd{anne.taormina@kcl.ac.uk}
\emailAdd{mara.ungureanu@gmail.com}
\emailAdd{wendland@maths.tcd.ie}
\emailAdd{ida.ghazvini-zadeh@soton.ac.uk}

\title{Tracking the symmetries of $\mathbb Z_3$-orbifold K3s within the Mathieu groups}
\abstract{ 
For $\Z_3$-orbifold limits of K3, we provide a counterpart to the extensive
studies by Nikulin and others of the geometry and symmetries of 
classical Kummer surfaces. In particular, we determine the group of holomorphic symplectic automorphisms of $\Z_3$-orbifold limits of K3. We moreover track this group within two of the Mathieu groups, 
which involves a variation of Kondo's lattice techniques that Taormina 
and Wendland introduced earlier in their study of the symmetries of Kummer surfaces 
and the genesis of their symmetry surfing programme. 
Specifically,  we realise the finite group of symplectic automorphisms of this 
class of K3 surfaces as a subgroup of the sporadic groups Mathieu $12$ and Mathieu $24$ in terms 
of permutations of $12$, resp. $24$ elements.  As a proof of concept, we construct an embedding
that yields the largest Mathieu group when the symmetry group of $\Z_3$-orbifold K3s is
combined with all symmetries of Kummer surfaces.
}


\begin{document}
\maketitle
\section{Motivation, Overview and Summary of Results}
%
K3 surfaces continue to fascinate and inspire mathematicians and theoretical 
physicists long after their name was coined, probably by Andr\'e Weil in 1958 \cite{weil2009final}. 
Their geometry and symmetries provide valuable insights in a variety of contexts, not the 
least that of superstring theory, which is the backdrop of our motivation here. 

By definition, a K3 surface is a compact, complex manifold of dimension $2$ which is simply 
connected and has trivial canonical bundle. It therefore carries a holomorphic volume 
form $\eta$.
We are interested in finite groups of symplectic automorphisms of K3 surfaces, 
i.e. automorphisms that fix $\eta$ as well as some K\"ahler class $\omega$, 
and in how they are realised as subgroups of the Mathieu group $M_{24}$.

Indeed, 
the observation in \cite{eguchi2011notes} that the elliptic genus of K3 reveals
a connection with $M_{24}$ -- a  fact mathematically proven subsequently on the level of representation theory and (Mock) modular forms in 
\cite{gannon2016much} -- remains unexplained by geometry and string theory. 
This is surprising (and haunting!), since
the observation was first made in the context of superstring theory 
compactified on a K3 surface, where one deals with K3 theories, i.e. superconformal 
field theories with $\mathcal{N} = (2, 2)$ worldsheet supersymmetry at central charges $(c, \ol c) = (6, 6)$, 
with spacetime supersymmetry, integer $\mathfrak{u}(1)$-charges 
and without holomorphic BPS states at weight 1/2 \cite{se88,asmo94,nawe00}.  
The $\mathfrak{u}(1)$-charges
are the eigenvalues of the zero modes in the two $\mathfrak{u}(1)_1$ subalgebras of the left 
and right $\mathcal{N} = 2$ superconformal 
algebras, which by the operators of two-fold spectral flow are extended to affine subalgebras 
of type $\mathfrak{su}(2)_1$. 
Thereby, any K3 theory has $\mathcal{N}=(4,4)$ worldsheet superconformal symmetry, and its conformal field 
theoretic elliptic genus  can be expressed
in terms of $\mathcal{N}=4$ superconformal characters \cite{Eguchi:1988vra} and equals  the complex 
elliptic genus of K3 surfaces. 

The `Mathieu Moonshine' observation is based on the fact that in the 
decomposition of the superconformal field theoretic elliptic genus into characters of irreducible
$\mathcal N=4$ representations, all coefficients of massive characters are positive 
\cite{eguchi2011notes,gannon2016much}. That a deeper connection exists, at least on the level of
representation theory, becomes apparent upon considering the corresponding twining 
elliptic genera which are obtained from the elliptic genus of K3 by replacing the dimensions of
$M_{24}$ representations by their characters with respect to an $M_{24}$ group element. 
They
were first computed in \cite{cheng2010k3,gaberdiel2010twining,gaberdiel2010mathieu,eguchi2011note} and their 
behaviour indicates the existence of an underlying representation of the entire group $M_{24}$,
which was proven in \cite{gannon2016much}.
However,  no K3 theory can have $M_{24}$ as its symmetry
group, and even more disconcertingly, K3 theories can have symmetry groups that are not 
subgroups of $M_{24}$ \cite{ghv12}.
Indeed,  Gaberdiel, Hohenegger and Volpato classified the 
symmetry groups of K3 non-linear sigma models 
preserving the $\mathcal N=(4,4)$ symmetry as well as the spectral flow operators, 
showing that all of these symmetry groups are naturally
realized as subgroups of the Conway group $\mbox{Co}_1$ \cite{ghv12}.

This is in contrast to the fact that any finite group of symplectic 
automorphisms of a K3 surface is isomorphic to a 
subgroup of the Mathieu group $M_{23}$, as was proven 
in \cite{mukai1988finite} almost $25$ years before \cite{eguchi2011notes}. 
The \textsc{Kummer surfaces} form a  family of special K3 surfaces which are obtained
by $\Z_2$-orbifolding an arbitrary
complex two-torus, 
where the $\Z_2$-action is the standard one.
In \cite{tawe11} and \cite{tawe13}, two of the authors used lattice constructions 
inspired by the seminal works \cite{ni80b} and \cite{ko98a} to realise and naturally 
combine   finite groups of 
symplectic automorphisms of Kummer surfaces as subgroups of $M_{24}$. This kick-started our 
symmetry-surfing programme, which is supported by results in 
\cite{tawe12, gaberdiel2017mathieu, keller2020lifting, taormina2018not}, albeit not fully 
completed yet.
Apart from the proposal to combine symmetry groups of distinct models, this programme is
built on the idea that the `Mathieu Moonshine' phenomenon arises through symplectic automorphisms
of K3 surfaces only, as opposed to the general symmetries of K3 theories. This is in accord with the
findings of Gaberdiel, Hohenegger and Volpato \cite{ghv12} mentioned above: while every symplectic automorphism of 
a given K3 surface $X$ belongs to $M_{24}$ and 
induces `geometric' symmetries of the K3 sigma models with geometric interpretation 
on $X$, such sigma models usually have plenty of additional, `non-geometric' symmetries which need not belong
to $M_{24}$. The idea is also in line with
the arguments of \cite{ka05} which suggest that compatibility of symmetries of K3 sigma models with an infinite
volume limit is required for such symmetries to contribute to `Mathieu Moonshine', thus enforcing
the restriction to `geometric' symmetries, see e.g.~\cite[\S4]{we14} and \cite[\S4.5]{we17} for further details.

For these reasons, in our work we focus on `geometric' symmetries and their action as symplectic
automorphisms of K3 surfaces.
In \cite{tawe13}, it was argued that the group
$(\Z_2)^4\rtimes A_8$, which is a maximal subgroup of $M_{24}$, arises as the combined holomorphic 
symplectic automorphism group of all Kummer surfaces whose K\"ahler class is induced by the 
underlying complex torus. The prescription for such combination of symmetries was given 
in \cite{tawe11} by means of overarching maps between pairs of Kummer surfaces, a technique which 
allows to symmetry-surf the moduli space of K3 theories. 
While representing symmetries on the lattice of integral cohomology as suggested by
\cite{ni80b,ko98a}, 
the crucial variation to these techniques which was introduced in \cite{tawe11}
uses the Kummer lattice rather than the orthogonal complement of the invariant sublattice. 

Encouraged by the above results concerning the specific class of Kummer K3 surfaces, that is,
the $\Z_2$-orbifold limits of K3, and by using analogous techniques, we set out to 
realise the group of symplectic automorphisms of $\Z_3$-orbifold limits of K3 (also known as $\Z_3$-orbifold K3s)  
as a subgroup of a Mathieu group. 
This leads us to exhibit the group of symplectic automorphisms of 
$\Z_3$-orbifold K3s, which is $(\Z_3)^2\rtimes \Z_4$ or $(\Z_3)^2\rtimes \Z_2$, depending on the 
symmetries of the underlying torus, 
as a subgroup of the Mathieu group $M_{12}$, but also as a subgroup of $M_{24}$. 
By combining the generators of the groups $(\Z_2)^4\rtimes A_8$ and 
$(\Z_3)^2\rtimes \Z_4$, one generates $M_{24}$. 
Although the latter step currently  relies on  ad hoc choices, it contributes another piece to  the Mathieu Moonshine puzzle, 
which we plan to solve in the future by providing  a geometric or conformal field theoretic justification for how these 
symmetry groups can be combined via symmetry-surfing. 
This would have to include an explanation of the action of these symmetries on quarter BPS states in the spectrum
of K3 sigma models.
We view our results, which pertain to geometry and 
group theory, as valuable independently of Mathieu Moonshine: 
based on Nikulin's seminal work \cite{ni75},
there exists a vast body of literature dedicated to Kummer surfaces,  
but we are unaware of a similarly detailed 
treatment for $\Z_3$-orbifold K3s. 
We aim to close that gap in this work. 
\bigskip

In the following,  by $X=\widetilde{T/\Z_3}$ we denote a $\Z_3$-orbifold K3, with $T$ 
a complex two-torus allowing a faithful action of $\Z_3$ by holomorphic symplectic 
automorphisms, such that $X$ is the minimal resolution of all singularities in 
the $\Z_3$-quotient. To introduce the notion of symmetries of $X$, we equip $X$ with 
a K\"ahler class $\wt\omega$ that is induced by any $\Z_3$-invariant K\"ahler class on $T$; a
\textsc{symmetry} of $X$ then is a biholomorphic automorphism of $X$ which preserves $\wt\omega$.
We now summarize the results obtained in this paper:
\begin{enumerate}
\item 
The construction of $X$ by 
means of a $\Z_3$-orbifolding procedure 
naturally generalises the classical Kummer construction: it is known and reviewed in subsection \ref{subsec:Z3orbi}. 
An alternative construction of $\Z_3$-orbifold K3s
(presumably known to experts), 
which consists in first blowing 
up sufficiently many points on $T$ and then quotienting by 
$\Z_3$, is presented in 
subsection \ref{subsec:Z3orbi2}, as we were unable to find 
the details in the literature, but also because it allows us to obtain the integral cohomology of
$X$ without passing through orbifold cohomology or equivariant cohomology
in subsection \ref{subsec:integralhomology}. It  thereby nicely 
prepares the ground for the study of symmetries of $X$.
\item 
A novel determination of the two  complementary
sublattices $K$  (proposition \ref{Kconstruct}) and $P$ (proposition \ref{Pconstruct}) 
of the integral cohomology lattice $H^2(X,\Z)$, where $K$ and $P$ 
are the minimal primitive sublattices containing the 
contributions from the underlying torus and from the exceptional divisors of the blow-ups respectively. 
The lattice  $H^2(X,\Z)$, which is even and self-dual, can thus  be constructed from $K$ and $P$ by using Nikulin's gluing techniques. 
The lattice $P$ is the analog of the Kummer lattice and was first constructed in \cite{be88},
using the results of \cite{shino77}, where the lattice $K$ was first determined; our derivation
is independent of \cite{shino77}.
The explicit gluing prescription is given in 
proposition \ref{gluePandK}.
\item
The description of $H^2(X,\Z)$ in terms of the lattices $K$ and $P$ by 
gluing provides an ideal set-up to determine the entire group of 
symmetries of $X$, which is $(\Z_3)^2\rtimes \Z_4$ (lemmas \ref{translationalaction} and \ref{rotationalaction})
or $(\Z_3)^2\rtimes \Z_2$ (remark \ref{generalsymmetries}), depending 
on the symmetries of the underlying torus $T$, 
and its action on $H^2(X,\Z)$. In the process, we prove that all symmetries are induced from symmetries 
of the underlying torus (proposition \ref{allsymmetries} and remark \ref{generalsymmetries}) 
and that these symmetries are uniquely
determined by their action on the lattice $P$ (corollary  \ref{F32isenough}).
Similarly to the strategy devised by two of the authors in \cite{tawe11,tawe13}, we may therefore 
use the lattice $P$ rather than the orthogonal complement of the invariant sublattice 
of $H^2(X,\Z)$ in order to track the symmetries of $X$. In the Kummer case, this was the 
crucial variation of the techniques developed by Kondo in \cite{ko98a}  that allowed
symmetry surfing in the first place; in the context of $\Z_3$-orbifold K3s we 
also find it convenient, since the lattice $P$ can be described entirely in 
terms of its root sublattice, while the orthogonal complement of the invariant
sublattice does not contain any roots.
\item 
The explicit construction of a primitive embedding of the lattice $P(-1)$ 
(obtained from $P$ by reversing the signature) in the Niemeier lattice of type $A_2^{12}$, 
called $N$ in this work (proposition \ref{embeddingexists}), and the proof that this embedding is 
unique, up to  automorphisms of $N$   (proposition \ref{primitiveembeddingisunique}). 
\item 
The proof that  the Niemeier lattice $N$ is, up to lattice automorphisms, 
the unique positive definite self-dual lattice of rank 24 which contains the lattice $P(-1)$
as a primitive sublattice (proposition \ref{primitiveNiemeierisunique}).
\item 
The proof that non-primitive embeddings of $P(-1)$ in Niemeier lattices do exist, but 
only in Niemeier lattices of type $E_6^4$ (appendix \ref{app:npembedding}). 
This allows us to illustrate the relevance of primitivity: it enables  
to track the symmetry group of $X$ in terms of its action on the Niemeier 
lattice of type $A_2^{12}$ (proposition \ref{rotationsonN}) and thereby as subgroup of the
Mathieu group $M_{12}$, while we argue that some of the so-called rotational symmetries of 
$X$ cannot be tracked on a Niemeier lattice of type $E_6^4$ 
(remarks \ref{E64} and \ref{symmetrybreaking}). 
\item 
The tracking of the symmetry group of $X$ within $M_{12}$ 
(corollary \ref{subgroupofM12}) and $M_{24}$ (proposition \ref{imageinM24})
explicitly in terms of permutations on $12$, respectively $24$ elements, 
exploiting the fact that $M_{12}$ is a subgroup of $M_{24}$.
\item 
The generation of $M_{24}$ from the symmetries of all Kummer surfaces and those of 
$X$ (theorem \ref{generatingM24}).
\end{enumerate}
Our presentation is organized as follows.
After this introductory section, the main text has two parts,   where 
section \ref{sec;Z3orbifoldK3s} is devoted to the investigation of $\Z_3$-orbifold
K3s and their symmetries in terms of the geometric orbifold construction, 
while in section \ref{sec:embedding} we focus on tracking 
these symmetries inside the Mathieu groups $M_{12}$ and $M_{24}$.

To smoothen the flow of the main text, and to allow the expert reader
a more immediate access to the new results, while keeping this work as 
self-contained as possible, additional material is found in a number of appendices. 
Appendix \ref{app:A2blowup} reviews the blow-up of
the particular singularities which arise within the standard $\Z_3$-orbifold 
construction, commonly known as $A_2$ singularities. In appendix \ref{app_ellgen},
we recall the calculation of the complex elliptic genus of $X$, which can be used
to confirm that $X$  is a K3 surface. We  present the calculation to emphasize 
the crucial role that conformal field theory has played in this area of mathematics. 
Some background on lattices, which play a fundamental role throughout our work, 
is given in appendix \ref{app:glue}, with some basic facts  assorted
in appendix \ref{subapp:background},
and a summary of Nikulin's gluing technique in appendix \ref{subapp:glue}. The Mathieu groups
$M_{12}$ and $M_{24}$, viewed as quotients by the maximal normal subgroups in the 
groups of lattice automorphisms of the Niemeier lattices of type $A_2^{12}$ and $A_1^{24}$,
are the theme of appendices \ref{subapp:niemeier} and \ref{subapp:m24}. 
Finally, appendix \ref{app:npembedding} provides the proof of the claim that
the Niemeier lattice of type $E_6^4$, up to lattice automorphisms, is unique
in allowing a non-primitive embedding of the lattice $P$ after reversal of the signature.


\section{\texorpdfstring{$\Z_3$}{TEXT}-orbifold limits of K3 surfaces}\label{sec;Z3orbifoldK3s}

This section is devoted to the discussion of the geometry of $\Z_3$-orbifold limits of K3 surfaces, 
in particular of their symmetry groups.

In subsection \ref{subsec:Z3orbi}, we review the $\Z_3$-orbifold 
construction for K3 surfaces, where 
the K3 surface is obtained by resolving all singularities in $T/\Z_3$ with $T$ an appropriate complex
two-torus. This construction is analogous to the standard Kummer construction.
That such $\Z_3$-orbifold limits of K3 exist has been known before, see for example the classical work
by Shioda and Inose \cite[Lemma 5.1]{shino77}, or Walton's discussion from a string theory viewpoint
\cite[\S3.1]{wa88}. In subsection \ref{subsec:Z3orbi2},  we offer an alternative 
construction, in which we first blow up $T$ sufficiently often, such that a $\Z_3$-quotient of the 
blown-up surface yields  the $\Z_3$-orbifold K3. This construction avoids singular intermediate
steps, which we find convenient later on. Though the existence of such a construction
is apparently known to the experts (see, for example, \cite[(2.2)]{be88}), to our knowledge the
details have not been described in the literature before.

The Torelli Theorem for K3 surfaces implies that the 
 symmetry group of $\Z_3$-orbifold K3s can be described entirely in terms of its action on the
integral second cohomology of K3 surfaces, as we shall recall at the beginning of subsection \ref{subsec:Symmetries1}.
As a preparation, we therefore discuss the integral  cohomology of $\Z_3$-orbifold K3s in 
subsection \ref{subsec:integralhomology}, in terms of the contributions that arise from the 
orbifolding construction, namely, those that descend from the underlying torus, on the one hand,
and those that arise from the blow-ups, on the other. This amounts to describing the integral second
cohomology of K3 in terms of two primitive sublattices, called $K$ and $P$,
which are glued together by means of Nikulin's gluing technique, recalled  in appendix
\ref{subapp:glue}. Here, within the integral second K3 cohomology, the lattice $K$ is the smallest 
primitive sublattice
which contains all integral classes that descend from the underlying torus, while $P$ is the
orthogonal complement of $K$. 
Both lattices have been determined before, see \cite[Lemma 5.1]{shino77} for the lattice $K$
and \cite{be88,Wendland:2000ye,Wendland:2000ry} for the lattice $P$. However, in \cite{shino77} 
Shioda and Inose use methods different from ours, while the references \cite{be88,Wendland:2000ye,Wendland:2000ry}
use the results of \cite{shino77}. Our current approach instead leads to a novel
proof of the fact that $K$ is an index $3$ sublattice of the push-forward of the $\Z_3$-invariant
second torus cohomology. 

In subsection \ref{subsec:Symmetries1} we finally obtain the symmetry group of $\Z_3$-orbifold
K3s, that is, the group of biholomorphic automorphisms of our K3 surface which leave the holomorphic volume
form and the K\"ahler class invariant. Here, we choose any (degenerate) K\"ahler class that descends from a
$\Z_3$-invariant K\"ahler class on the underlying torus. 
We show that  the symmetry group is isomorphic to $(\Z_3)^2\rtimes \Z_4$ or $(\Z_3)^2\rtimes \Z_2$,
depending on the choice of the K\"ahler class,
and we prove that all symmetries are induced by symmetries of the underlying torus.

\subsection{Constructing \texorpdfstring{$\Z_3$}{TEXT}-orbifold limits of K3 as blow-ups of quotients}\label{subsec:Z3orbi}

We denote by $T$ the complex two-torus given as a product of two elliptic curves with $\Z_3$ symmetry,
that is, $T=\C^2/L$, where the lattice $L$ is generated by 
\begin{equation}
\lambda_1 = (1,0),\;\lambda_2=(\xi,0),\;\lambda_3=(0,1),\;\lambda_4=(0,\xi); \quad
\xi=\exp(\frac{2\pi i}{3})=-\frac{1}{2}+\frac{\sqrt{3}}{2}i. \label{eq:torus}
\end{equation}
This fixes 
a complex structure on $T$  
induced by the underlying $\C^2$ such that  the standard $\Z_3$-action on $\C^2$
by unitary, orientation preserving linear
maps generated by 
\be\label{Z3action}
\C^2\longrightarrow\C^2, \qquad
(z_1,z_2)\mapsto (\xi z_1, \xi^{-1} z_2) \; ,
\ee
 descends to a holomorphic action on $T$.
This is immediate since this action maps $L$ to itself.
By the results of \cite{fu88} we may choose this presentation of the lattice $L$
underlying a complex two-torus with a holomorphic action of $\Z_3$ without loss of generality, 
since it amounts to a convenient choice of coordinates.
We observe that the standard Euclidean metric on $\C^2$ is invariant under this action, 
but so is any metric that realizes $T$ as the product of two elliptic curves with complex structure
modulus $\xi$ and arbitrary, potentially different volumes. We will come back to the discussion
of compatible metrics in proposition \ref{hyperk}. Until then, our investigations solely require us
to fix the complex structure of $T$ as stated.

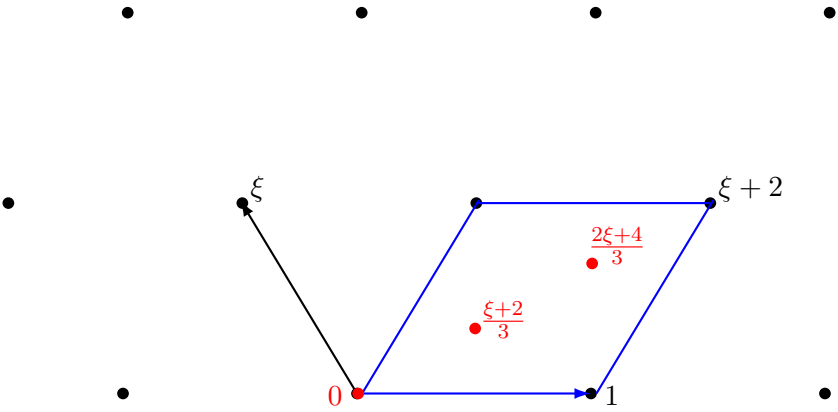
\begin{figure}[ht]
\begin{center}
\unitlength0.8em
\begin{picture}(40,20)(-15,0)
\thicklines
\multiput(-10,0)(10,0){4}{\circle*{0,5}}
\multiput(-14.9,8.2)(10,0){4}{\circle*{0,5}}
\multiput(-9.8,16.4)(10,0){4}{\circle*{0,5}}
\textcolor{blue}{\put(0,0){\vector(1,0){10}}}
\put(0,0){\vector(-3,5){5}}
\textcolor{blue}{\put(5.1,8.2){\line(1,0){10}}}
\textcolor{blue}{\put(-0.2,0){\line(3,5){5}}}
\textcolor{blue}{\put(9.4,0){\line(3,5){5}}}
\textcolor{red}{\multiput(-1.2,0)(5,2.8){3}{\circle*{0,5}}}
\put(-5.8,8.5){$\xi$}
\put(14.2,8.5){$\xi+2$}
\textcolor{red}{\put(-2.5,-0.5){$0$}}
\textcolor{red}{\put(3.6,2.8){$\frac{\xi+2}{3}$}}
\textcolor{red}{\put(7.8,6){$\frac{2\xi+4}{3}$}}
\put(8.5,-0.5){$1$}
\end{picture}
\end{center}
\caption{\small The lattice $\mbox{span}_\Z\{1,\xi\}$ with generators $1$ and $\xi=\exp(2\pi i/3)$,
lattice points $\bullet$,
a fundamental cell (blue), and the three representatives  \textcolor{red}{$\bullet$}
of fixed points under the $\Z_3$-action generated by $[z]\mapsto [\xi z]$
on $\C/\mbox{span}_\Z\{1,\xi\}$ within that cell.
We include the half-open segments $[0,1)$ from $0$ to $1$ and 
$[0,\xi+1)$ from $0$ to $\xi+1$ in the fundamental cell, but not the 
rest of the boundary.}
\label{Z3torus}
\end{figure}

Let us determine the fixed points on $T$ under the action \eqref{Z3action}. By definition,
$[(z_1,z_2)]\in T$ is a fixed point if and only if 
$(\xi z_1, \xi^{-1} z_2)-(z_1,z_2)\in L$, that is, if and only if both $z_1$ and 
$z_2$ represent a fixed point in $\C/\mbox{span}_\Z\{1,\xi\}$ 
under the $\Z_3$-action generated by $[z]\mapsto [\xi z]$.
These latter fixed points are most conveniently determined by focusing on the 
fundamental cell with vertices $0,\, 1,\, \xi+2,\, \xi+1$
in $\C$, see figure~\ref{Z3torus}. 

The map $z\mapsto\xi z$ is a rotation
by $\frac{2\pi}{3}$, such that for any $z$ in that fundamental cell,
$\xi z-z\in \mbox{span}_\Z\{1,\xi\}$ implies $\xi z-z\in \{0,\,-1,\,-2\}$.

Solving for $z$ yields precisely three distinct fixed points 
in $\C/\mbox{span}_\Z\{1,\xi\}$, namely 
$[0],\, [\frac{\xi+2}{3}]$, and $[\frac{2\xi+4}{3}]$.
Therefore there are precisely nine fixed points in $T$ under the action induced
by \eqref{Z3action}, given by $[(z_1,z_2)]$ with $z_1,z_2\in
\frac{1}{3}\mbox{span}_\Z\{\xi+2\}$.
In other words, the fixed points in $T$ are given by the elements of 
$\frac{1}{3}\widetilde L/\widetilde L$, where $\widetilde L\subset L$
is the rank-$2$ sublattice which is generated by $2\lambda_1+\lambda_2,\,
2\lambda_3+\lambda_4$. 
\bigskip

The $\Z_3$-orbifold K3s that we are studying in this work are obtained 
from the quotient $T/\Z_3$ 
through a minimal resolution of each singularity. By construction,
the singularities of $T/\Z_3$ are precisely the images of the 
$\Z_3$-fixed points in the quotient. Thus by the above, they are labelled by
$\frac{1}{3}\widetilde L/\widetilde L$. Since $\widetilde L\subset L$
is a rank-$2$ sublattice, we find $\frac{1}{3}\widetilde L/\widetilde L\cong\F_3^2$,
where $\F_3$ denotes the field with three elements. The 
fixed points are thus given by $\left[(t_1,t_2) \frac{\xi+2}{3}\right]$, 
with labels $t=(t_1,t_2)\in\F_3^2$ and where
we denote the three elements of $\F_3$ by $0,1,2$ or
equivalently by $0,1,-1$.
We observe that the label set naturally carries the structure of 
a two-dimensional affine space over $\F_3$, where this structure
is inherited from the group structure of the underlying torus $T$.

Again by construction, each of the nine singular points in 
$T/\Z_3$ is modelled by $[0]\in\C^2/\Z_3$ with the $\Z_3$-action 
generated by \eqref{Z3action}. As is explained in 
appendix \ref{app:A2blowup},
these singularities are thus of type $A_2$, and by blowing up
once in the singular point, the singularity is minimally resolved.
The exceptional divisor of the blow-up decomposes into two irreducible
components, each given by a copy of $\P^1$
that forms a $(-2)$-curve, and where the 
two components intersect transversally in one point, 
see figure~\ref{exceptionaldivisor}. The surface that we obtain 
from $T/\Z_3$ by minimally resolving all singularities is 
denoted $X=\widetilde{T/\Z_3}$, and we call it a $\mathbb{Z}_3$-\textsc{orbifold limit of a K3 surface}, 
or simply a $\mathbb{Z}_3$-\textsc{orbifold K3}.
\bigskip

Our terminology for the manifold $X$
is justified, since on the one hand, $X$  is a K3 surface as we
shall show  below. On the other hand, the notion
of \textsc{orbifold} is commonly used as a shorthand to describe the
construction of a complex manifold by quotienting another complex manifold 
by a group of symmetries and then minimally resolving all resulting singularities.
For brevity, in this work we will use the notion of orbifolding 
solely to refer to orbifolds of \textsc{complex tori}. 

Finally, 
the terminology orbifold \textsc{limit} is commonly used to remind 
of the fact that the resulting complex manifold, despite being non-singular,
still bears the traces of the singularities of the quotient. Indeed, 
the  occurrence of $(-2)$-curves is an indicator
of the fact that singularities have been resolved. Moreover, any choice 
of $\Z_3$-invariant \textsc{K\"ahler class} on the torus that underlies our 
orbifold induces a K\"ahler class on the orbifold, which however is degenerate
in the sense that it lies on the boundary of the K3 K\"ahler cone. 

We now assert 

\bprop\label{orbifoldisK3}%
The $\Z_3$-orbifold $X=\widetilde{T/\Z_3}$ is a K3 surface. 
\eprop

\bpf
By the definition of K3 surfaces,  we must show that $X$ 
is a connected, compact complex surface with trivial first cohomology, 
that is, $H^1(X,\C)=\{0\}$, and with trivial canonical bundle. 

Since $T$ is a connected, compact surface and $X$ arises from $T/\Z_3$ by replacing 
each of the nine singular points by a pair of curves that are isomorphic
to $\P^1$, it follows that $X$ is connected and compact. Moreover, 
though $T/\Z_3$ is a singular complex variety, $X$ arises by resolving
all singularities, hence altogether $X$ is a connected, compact complex 
surface. 

Using the standard
coordinates $(z_1,z_2)$ on $\C^2$, a basis of $H^1(T,\C)$ is represented
by $(dz_1,\, dz_2,\, d\overline z_1,\, d\overline z_2)$, such that none of the non-zero
classes in $H^1(T,\C)$ is invariant under the $\Z_3$-action \eqref{Z3action}.
Moreover, since the resolution of singularities in $T/\Z_3$
introduces $18$ additional copies of $\P^1$, none of which carries
a non-trivial global, closed one-form, one may conclude that indeed 
$H^1(X,\C)=\{0\}$, as claimed. 

Finally, to see that $X$ has a trivial canonical bundle, it 
suffices to show that on $X$, there exists a holomorphic $(2,0)$-form 
which vanishes nowhere. To prove the existence of such a form, note 
that $dz_1\wedge dz_2$ descends to a holomorphic $(2,0)$-form 
which vanishes nowhere on $T$. This differential form is invariant under 
our $\Z_3$-action \eqref{Z3action}, and thus it
yields a holomorphic $(2,0)$-form which vanishes nowhere
on $T/\Z_3\setminus\Delta$, with $\Delta$ denoting the set of all nine 
singularities. This differential form has a unique extension 
to the exceptional divisor under the minimal resolution of the 
singularities, as holomorphic $(2,0)$-form which vanishes nowhere; 
the explicit calculation is shown in appendix \ref{app:A2blowup}, see 
\eqref{OmegaonU0}--\eqref{OmegaonU2}. 
\epf

We have thus confirmed the claim
that $X$ is a K3 surface, a fact which can 
be expected from the calculation of the elliptic genus of $X$ using 
orbifold conformal field theory techniques.  
Indeed, recall that on the one hand, the complex elliptic 
genus $\mathcal E_Y$ of a compact, complex manifold $Y$
encodes a number of properties and invariants of the manifold 
$Y$, see for example \cite{hi66,kr90,wi88,boli00}, or for a review
\cite[\S2.4]{we14}. 
For instance, if $\mathcal E_Y$ is a weak Jacobi form of weight zero 
and index $D/2$, then $Y$ has complex dimension $D$ and a trivial canonical
bundle. Moreover, $\mathcal E_Y(\tau,z=0)$ is constant and yields the Euler
characteristic of $Y$, while the 
constant term  in the 
$q$-expansion of $\mathcal E_Y(\tau,z=1/2)$ is, up to a
prefactor $(-1)^{D/2}$, the signature of $Y$, and the constant term  in the 
$q$-expansion of $q^{D/4}\mathcal E_Y(\tau,z=\frac{\tau+1}{2})$ is, 
up to a prefactor $(-1)^{D/2}$,
the holomorphic Euler characteristic of $Y$. 

On the other hand, the construction of the complex elliptic genus of $Y$
is motivated from string theory, with the expectation that it agrees 
with the conformal field theoretic elliptic genus of every $\sigma$-model 
with target $Y$. If the target manifold is the orbifold limit of a
connected, compact complex manifold with trivial canonical bundle, 
then its elliptic genus  can be
calculated by means of orbifold techniques in conformal field theory. 
Appendix \ref{app_ellgen} illustrates these techniques for the case of 
$\widetilde{T/\Z_3}$, which is of interest in this work.
There we recall that orbifold techniques allow to recover the well-known 
expression \eqref{ellk3} for the complex elliptic genus
$\mathcal E_X$ of K3 surfaces from conformal field theory. It is a weak 
Jacobi form of weight zero and index $1$, confirming that $X$ is a complex
surface with trivial canonical bundle. Furthermore, from $\mathcal E_X$
one obtains the Euler characteristic
$24$, the signature $-16$, and the holomorphic Euler characteristic $2$
for $X$, confirming that $X$ is a K3 surface.
\bigskip
\subsection{Constructing \texorpdfstring{$\Z_3$}{TEXT}-orbifold limits of K3 as quotients of blow-ups}\label{subsec:Z3orbi2}
As mentioned in the introduction, our $\Z_3$-orbifold construction 
is similar to the Kummer construction for K3 surfaces by $\Z_2$-orbifolding.
In that construction, it is useful to observe that instead of first
taking the quotient of a complex two-torus by $\Z_2$ and then resolving
the singularities in the quotient, one may proceed in reverse order:
first blow up the fixed points of the $\Z_2$-action on the smooth torus,
then construct a continuation of the $\Z_2$-action to the resulting 
surface, and finally consider the quotient by $\Z_2$. 
One checks that the result yields back the Kummer surface, 
see e.g.~\cite[V.16]{bpv84}.

We have a similar situation\footnote{KW is grateful to Kenji Ueno for some crucial hints 
concerning this idea.} for our Kummer-like surface $X=\widetilde{T/\Z_3}$:

\bconst
To construct $X=\widetilde{T/\Z_3}$, one may alternatively first blow-up $T$
twenty-seven times, then take a quotient by $\Z_3$, and finally blow 
down nine $(-1)$-curves.
\econst

Though this is apparently known to the experts (see, for example, \cite[(2.2)]{be88}), we 
are not aware of such a construction being presented in the literature, 
yet it proves convenient for our purpose
and we therefore elaborate on it here.

First consider the local model for each fixed point under
the $\Z_3$-action, given by the origin $(z_1,z_2)=(0,0)$ in $\C^2$.
The blow-up of $\C^2$ in the origin is the closure of 
$$
\left\{ \left( (z_1,z_2), (v_1\colon v_2)\right) \in \C^2\times\P^1
\mid (z_1,z_2)\neq (0,0),\; z_1v_2=z_2v_1\right\}
$$
in $\C^2\times\P^1$,
such that the exceptional divisor is $C\subset\C^2\times\P^1$ with
$C=\{0\}\times\P^1$, an isomorphic copy of $\P^1$.
The curve $C$ intersects the strict transforms of the coordinate 
axes of $\C^2$ in the points $P_1=\left( (0,0), (1\colon 0)\right)$
and $P_2=\left( (0,0), (0\colon1)\right)$, respectively. 

\begin{figure}[ht]
\begin{center}
\unitlength1em
\begin{picture}(20,7)(3,0)
\thicklines
\put(0,0){\line(1,0){6}}
\put(0,5){\line(1,0){6}}
\put(5,-1){\line(0,1){7}}
\put(-0.5,-1){\small$\{z_1=0\}$}
\put(-0.5,5.5){\small$\{z_2=0\}$}
\put(5.3,2.3){\small$C$}
\put(3.8,4.1){\small$P_1$}
\put(3.8,0.3){\small$P_2$}
\put(15,2.5){\vector(-1,0){5}}
\put(18,0){\line(1,0){6}}
\put(18,5){\line(1,0){6}}
\put(25,1){\line(0,1){3}}
\put(23,-0.5){\line(1,1){2.5}}
\put(23,5.5){\line(1,-1){2.5}}
\put(18,-1){\small$\{z_1=0\}$}
\put(18,5.5){\small$\{z_2=0\}$}
\put(23.8,2.3){\small$C_0$}
\put(24.5,4.4){\small$C_1$}
\put(24.5,0.1){\small$C_2$}
\end{picture}
\end{center}
\caption{\small Blow-up of $\C^2$ in the origin, followed
by blow-ups in $P_1$, $P_2$.}
\label{tripleblowup}
\end{figure}
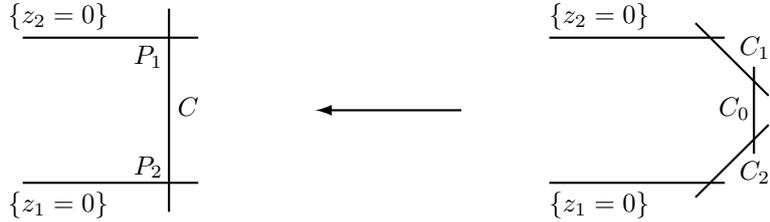

We now blow up in the points
$P_1$ and $P_2$ (see figure~\ref{tripleblowup}) 

and obtain a surface $Z$, the closure of the set
\begin{eqnarray}\label{surfaceZ}
&&\left\{ \left( (z_1,z_2)  ,(u_1\colon u_2), (v_1\colon v_2), (w_1\colon w_2)\right) \in \C^2\times \left(\P^1\right)^3 \mid \right.\\
&&\hspace*{10em}\left.(z_1,z_2)\neq (0,0),\; z_1v_2=z_2v_1,\; 
z_1u_2v_1=u_1v_2,
\; z_2v_2w_1=v_1w_2 \vphantom{\left(\P^1\right)^3 } \right\}\nonumber
\end{eqnarray}
in $\C^2\times \left(\P^1\right)^3$. Altogether we have performed
three blow-ups of points, yielding $\sigma\colon Z\longrightarrow\C^2$ with
$\sigma(z,u,v,w) = z$ for $(z,u,v,w)\in Z\subset \C^2\times \left(\P^1\right)^3$,
hence the exceptional divisor is $\sigma^{-1}(0) = C_1\cup C_0\cup C_2$ with
\begin{eqnarray*}
C_1 &=& \left\{ \left(0, (u_1\colon u_2), (1\colon 0), (1\colon 0)\right) 
\mid (u_1\colon u_2)\in\P^1\right\}, \\
C_0 &=& \left\{ \left(0, (0\colon 1), (v_1\colon v_2), (1\colon 0)\right) 
\mid (v_1\colon v_2)\in\P^1\right\}, \\
C_2 &=& \left\{ \left(0, (0\colon 1),(0\colon 1), (w_1\colon w_2)\right) 
\mid (w_1\colon w_2)\in\P^1\right\}, 
\end{eqnarray*}
yielding three isomorphic copies $C_j$, $j\in\{0,\,1,\,2\}$, of $\P^1$. Here,
$C_1$ intersects the strict transform of the $z_1$-axis of $\C^2$
under $\sigma$ transversally in $\left(0, (1\colon 0), (1\colon 0), (1\colon 0)\right)$, $C_2$ intersects the strict transform of the $z_2$-axis of $\C^2$
under $\sigma$ transversally in $\left(0, (0\colon 1),(0\colon 1), (0\colon 1)\right)$, 
while $C_0$ is the strict transform of $C$ under the
blow-up of $P_1$ and $P_2$, intersecting $C_1$ transversally in 
$\left(0, (0\colon 1), (1\colon 0), (1\colon 0)\right)$
and $C_2$ in $\left(0, (0\colon 1), (0\colon 1), (1\colon 0)\right)$.

We claim that the natural $\Z_3$-action on $\C^2$, which is generated by
\eqref{Z3action},
induces a $\Z_3$-action on $Z$ in such a way 
that the image of $C_0$ in $Z/\Z_3$
can be blown down, resulting in a  surface
isomorphic to the minimal resolution $\widetilde{S(A_2)}$
of $\C^2/\Z_3$ discussed in appendix \ref{app:A2blowup}. Indeed,  
\eqref{surfaceZ} implies that \eqref{Z3action} uniquely induces 
a biholomorphic $\Z_3$-action on $Z$ which is generated by
\be\label{actiononZ}
Z\longrightarrow Z, \qquad
\left( (z_1,z_2)  ,u, (v_1\colon v_2), w\right)  
\longmapsto
\left( (\xi z_1, \xi^{-1} z_2)  ,u, (\xi v_1\colon \xi^{-1} v_2), w\right) \; .
\ee
From \eqref{actiononZ} we observe
that the curves $C_0$, $C_1$, $C_2$, all of which are $(-1)$-curves 
since they arose from blowing up smooth points (see, for example,
\cite[p.~475]{grha78}), are each mapped
to themselves under the $\Z_3$-action. In fact, $C_1$ and $C_2$ are 
fixed pointwise, while $\Z_3$ acts non-trivially on $C_0$. The images of 
$C_1$ and $C_2$ in $Z/\Z_3$ therefore are $(-3)$-curves, while the image
of $C_0$ is still a $(-1)$-curve. By the Castelnuovo-Enriques Criterion 
(see, for example, \cite[p.~476]{grha78}), this $(-1)$-curve can
be blown down to a smooth point. Denote the resulting surface by $\widetilde S$;
by the known behaviour of the intersection form under such blow-downs  (see, for example, \cite[p.~476]{grha78}), the 
images of $C_1$ and $C_2$ descend to $(-2)$-curves $C^{(1)}$ and $C^{(2)}$ that intersect transversally in one point. 

Concretely,  by introducing 
$x_0:=z_1z_2$, $x_1:=z_1^3$,  $x_2:=z_2^3$ with $x_0^3=x_1 x_2$, from 
\eqref{surfaceZ} we find the resulting surface $\widetilde S$ as the
completion of 
\begin{eqnarray*}
&&\left\{ \left( (x_0, x_1, x_2)  ,(u_1\colon u_2), (w_1\colon w_2)\right) \in \C^3\times \left(\P^1\right)^2 \mid (x_0, x_1, x_2)\neq (0,0,0), \right.\\
&&\hspace*{10em}\left.
x_0^3 = x_1x_2,\;
x_1u_2=x_0u_1,\; 
x_0w_2=x_2w_1,\; 
u_1w_2=x_0u_2w_1
\vphantom{\left(\P^1\right)^3 } \right\}\nonumber
\end{eqnarray*}
in
$\C^3\times\left(\P^1\right)^2$. Hence
\begin{eqnarray*}
C^{(1)} &=& \left\{ \left(0, (u_1\colon u_2), (1\colon 0)\right) 
\mid (u_1\colon u_2)\in\P^1\right\}, \\
C^{(2)} &=& \left\{ \left(0, (0\colon 1), (w_1\colon w_2)\right) 
\mid (w_1\colon w_2)\in\P^1\right\},
\end{eqnarray*}
and $\widetilde S$  is indeed isomorphic to
the minimal resolution $\widetilde{S(A_2)}$ of an $A_2$-type singularity 
discussed in appendix \ref{app:A2blowup} by means of
$$
\widetilde S\longrightarrow  \widetilde{S(A_2)},\quad
\left( (x_0, x_1, x_2)  ,(u_1\colon u_2), (w_1\colon w_2)\right) 
\mapsto \left( (x_0, x_1, x_2)  ,(u_2 w_1\colon u_1 w_1\colon u_2 w_2)\right) \; ,
$$
which maps $C^{(k)}$ to $\CCC_{00}^{(k)}$ for $k\in\{1,2\}$.
The reader who is familiar with singularities of type $A_n$ with $n\in\N$
may check that the construction generalizes to all of these.

We may now implement this local construction in order to arrive 
at the promised second, very convenient description 
of our $\Z_3$-orbifold K3, 
$X=\widetilde{T/\Z_3}$. To this end, let $\widetilde T$ denote 
the $27$-fold blow-up of $T$, where in generalization of the above 
construction we first blow up each of the nine fixed points under
$\Z_3$ on $T$ and then blow up the analogs of the points $P_1$
and $P_2$ for each of them. By the above, our $\Z_3$-action on $T$
induces a $\Z_3$-action on $\widetilde T$ in such a way that 
$\widetilde T/\Z_3$ is smooth and contains nine $(-1)$-curves, 
one corresponding to each of the fixed points. Blowing all these 
$(-1)$-curves down, by what was said above, 
one recovers $\widetilde{T/\Z_3}$. 
Let $\widetilde\pi\colon \widetilde T\longrightarrow \widetilde{T/\Z_3}$ 
denote the map which first performs the $\Z_3$-quotient and then
blows down the nine $(-1)$-curves. Then altogether
we have a blowing-up and -down diagram as follows,
yielding an induced rational map $\pi\colon T\dashrightarrow X$:
\be\label{blowupanddown}
\begin{tikzcd}
\widetilde T  \arrow[d, "\widetilde\pi"] \arrow[r] &[0.5em] T \arrow[d, "\mod\Z_3"]
\arrow[dl, dashed, "{\pi}" description]
\\
X=\widetilde{T/\Z_3}\; \arrow[r] &T/\Z_3
\end{tikzcd}
\ee

\subsection{The integral cohomology of \texorpdfstring{$\Z_3$}{TEXT}-orbifold K3s}
\label{subsec:integralhomology}

As we shall explain in greater detail in subsection \ref{subsec:Symmetries1},
the classical Torelli Theorem for K3 surfaces (see \cite{pss71,burap75,lope81}) allows
to uniquely describe the symmetries of a K3 surface $X$ 
in terms of specific lattices, in particular the K3 lattice. 
We recall that the integral cohomology\footnote{The $\Z_3$-orbifold construction reviewed 
in the previous subsection is conveniently described in terms of homology. We are 
now shifting from homology to cohomology for ease of calculations, with no big 
price to pay as Poincar\'e duality ensures that the second homology and cohomology 
groups of a K3 surface are isomorphic. On cohomology we  thus use the bilinear
form which under Poincar\'e duality
is induced by the intersection form on homology. For brevity, we call it
and all its linear extensions to larger lattices or vector spaces
\textsc{intersection form} as well.}  $H^2(X,\Z)$ forms 
an even, unimodular lattice of signature $(3,19)$, called the \textsc{K3 lattice}. 
This lattice is isometric to $U^{\oplus3}\oplus E_8(-1)^{\oplus 2}$, 
with $U$ the hyperbolic lattice and $E_8(-1)$ the negative definite version 
of the root lattice of the Lie algebra $E_8$.

In this subsection, we prepare the ground for subsection \ref{subsec:Symmetries1}, 
where we give the action on the lattice $H^2(X,\Z)$ of the \textsc{symmetry group}
of the K3 surface $X=\widetilde{T/\Z_3}$ -- which according to our choice of terminology
is the group of biholomorphic automorphisms 
preserving the holomorphic volume form and the 
K\"ahler class induced by those of the underlying complex torus $T$. 
As we shall see, this group consists of the symmetries 
that are induced by symmetries
of the underlying complex torus, and as such it acts 
naturally on the torus cohomology $H^2(T,\Z)$.
In view 
of our subsequent goal, which is to track this group of symmetries in the Mathieu groups, 
we take advantage of lattice gluing techniques originating in the work of Nikulin \cite{ni80b} 
and briefly summarised in appendix \ref{app:glue}.  We use them to reconstruct the rank-$22$ 
lattice $H^2(X,\Z)$ through the gluing of 
\begin{enumerate}
 \item 
the smallest primitive sublattice $K$ of $H^2(X,\Z)$ containing all
pushforwards of $\Z_3$-invariant torus 2-forms under the rational map 
$\pi:T\dashrightarrow X$ found in \eqref{blowupanddown}; $K$ is
a rank-$4$ sublattice of $H^2(X,\Z)$;
\item  
the smallest primitive sublattice $P$ of $H^2(X,\Z)$ containing the 
Poincar\'e duals of the exceptional divisors resulting from blowing up the  
singularities of $T/\Z_3$; $P$ is a rank-$18$ lattice called the 
\textsc{Kummer-like} lattice for $X=\widetilde{T/\Z_3}$,
and it is the orthogonal complement of the lattice $K$ in $H^2(X,\Z)$.
\end{enumerate}
Since the  symmetries of $X$ that we wish to track act naturally  on $H^2(T,\Z)$
and thereby on $K$, gluing $H^2(X,\Z)$ from $K$ and $P$
will allow us to describe such symmetries in 
terms of the K3 lattice.

We thus need to determine the lattices $K$ and $P$. To do so, we first
determine the $\Z_3$-invariant part of $H^2(T,\Z)$ and its pushforward under
the degree $3$ rational map $\pi:T\dashrightarrow X$
which is defined outside the fixed points of the 
$\Z_3$-action on $T$  (lemma \ref{H2Tinv}). 
Under Poincar\'e duality, the cohomology classes that we 
determine in that lemma correspond to torus cycles in general position -- that is,
cycles that do not contain $\Z_3$-fixed points.
With $\pi_\ast:H^2(T,\Z)\longrightarrow H^2(X,\Z)$ denoting the 
linear map induced on cohomology by $\pi$,
we also determine $\pi_\ast\left(H^2(T,\Z)^{\Z_3}\right)$ in lemma \ref{H2Tinv}. 
In the classical Kummer construction, the analog of $K$ agrees 
with the pushforward of the integral torus cohomology, but for our $\Z_3$-orbifold limit 
$X$ of K3, the lattice $\pi_\ast \left(H^2(T,\Z)^{\Z_3}\right)$ is 
an index $3$ sublattice of $K$. Like in the classical Kummer construction, 
the lattice $P$ contains as a proper sublattice the rank-$18$ root lattice $R$ generated 
by the Poincar\'e duals of the exceptional divisors for the blow ups of the  
singularities of $T/\Z_3$. We proceed geometrically to obtain lattice vectors
that generate $P$ and $K$ respectively from $R$ and
$\pi_\ast \left(H^2(T,\Z)^{\Z_3}\right)$. The crucial ingredient is a list of 
special vectors  in the lattice $H^2(X,\Z)$ whose Poincar\'e dual cycles
possess unbranched $3\colon1$ covers in the blow up
$\widetilde T$ of $T$ that was constructed in \eqref{blowupanddown} (lemma \ref{gluelist}).  
The desired generators of $P$ and $K$ are then obtained as linear combinations of 
these special lattice vectors (propositions \ref{Pconstruct} and \ref{Kconstruct}).
Finally, proposition \ref{gluePandK} gives the gluing instructions to obtain $H^2(X,\Z)$
from $K$ and $P$.
\bigskip

To determine $H^2(T,\Z)^{\Z_3}$, let us first introduce some notation.
With standard coordinates $(z_1,z_2)$ as above, we identify
$\C^2\cong\R^4$ by writing $(z_1,z_2)=(x^1+ix^2, x^3+ix^4)$. The standard basis
$(e_1,\ldots,e_4)$ of $\R^4$ is identified with the coordinate basis 
$(\partial_1,\ldots,\partial_4)$ of the tangent space of $T$ at any point, yielding
the dual basis $(dx^1,\ldots,dx^4)$ for the cotangent space at any point. 
As a real manifold, $T$ can be viewed as Cartesian product of four circles, whose
fundamental group is given by the defining lattice $L$ in $T=\C^2/L$. We thus have
$L=H_1(T,\Z)$ and obtain a basis of $H^1(T,\Z)$ by taking as representatives the basis
$(\mu_1,\ldots,\mu_4)$ which is
dual to the basis $(\lambda_1,\ldots,\lambda_4)$ of $L$ given in \eqref{eq:torus}. 
We obtain
\begin{equation}\label{dualbasis}\textstyle
\mu_1=dx^1+\frac{1}{\sqrt{3}}dx^2,\quad 
\mu_2=\frac{2}{\sqrt{3}}dx^2,\quad 
\mu_3=dx^3+\frac{1}{\sqrt{3}}dx^4,\quad 
\mu_4=\frac{2}{\sqrt{3}}dx^4.
\end{equation}
Under the $\Z_3$-action \eqref{Z3action} which induces
$\lambda_1\mapsto -\lambda_1-\lambda_2\mapsto\lambda_2$ and
$\lambda_3\mapsto \lambda_4\mapsto-\lambda_3-\lambda_4$,
the differential one-forms $(\mu_1,\ldots,\mu_4)$ transform as 
\[
\mu_1\mapsto-\mu_2\mapsto\mu_2-\mu_1,\qquad\mu_3\mapsto-\mu_3+\mu_4\mapsto-\mu_4.
\]
One now finds $H^2(T,\Z)^{\Z_3}$
by a direct calculation. Moreover, on $\pi_\ast\left(H^2(T,\Z)^{\Z_3}\right)$
the quadratic form  is determined in \cite[Prop.\,1.1]{ino76}
(see also \cite[\S3]{mo84} for a detailed discussion which carries over verbatim to our
situation). 
Letting $\mu_{jk}:=\mu_j\wedge \mu_k$ when $j,k\in\{1,\ldots, 4\}$, 
we have:

\blem\label{H2Tinv}%
The lattice $H^2(T,\Z)^{\Z_3}$ is generated by the basis 
\[
\left( \mu_{13}-\mu_{24},\,\mu_{14}+\mu_{23}-\mu_{24},\,\mu_{12},\, \mu_{34}\right).
\]
Moreover, $\pi_\ast$ induces an isomorphism from $H^2(T,\Z)^{\Z_3}$ onto its image 
which multiplies the bilinear form by $3$. Hence
$\pi_\ast\left(H^2(T,\Z)^{\Z_3}\right)\cong H^2(T,\Z)^{\Z_3}(3)$, and with 
respect to the basis $\left( \pi_\ast \mu_{13}-\pi_\ast \mu_{24},\,\pi_\ast \mu_{14}+\pi_\ast \mu_{23}-\pi_\ast \mu_{24},\,
\pi_\ast \mu_{12},\, \pi_\ast\mu_{34}\right)$,  the
bilinear form on the lattice $\pi_\ast\left(H^2(T,\Z)^{\Z_3}\right)$
is given by
$$
\left(\begin{matrix}6&3&0&0\\3&6&0&0\\0&0&0&3\\0&0&3&0\end{matrix}\right)\; .
$$
\elem
\bigskip

Next, let us explain geometrically
how non-integral linear combinations of lattice vectors from
$\pi_\ast\left(H^2(T,\Z)^{\Z_3}\right)\oplus R$ can yield integral cohomology classes 
of $X$. Here, $R$ denotes the root lattice generated by the Poincar\'e duals of all
irreducible components of exceptional divisors that result from blowing up the singularities
in $T/\Z_3$. Recall from subsection \ref{subsec:Z3orbi} that there are nine fixed 
points $F_t$ of the $\Z_3$-action \eqref{Z3action} on $T$, labelled by $t\in\F^2_3$
with $\F^2_3\cong\frac{1}{3}\widetilde L/\widetilde L$. In 
$T/\Z_3$, each of these fixed points $F_t$ gives rise to a singularity of type $A_2$, 
whose resolution according to appendix \ref{app:A2blowup} yields an exceptional divisor
with two irrdeducible components $\CCC_t^{(k)}$, $k\in\{1,2\}$. The Poincar\'e dual
of $\CCC_t^{(k)}$ is denoted $E_t^{(k)}$ and by construction yields a root in $H^2(X,\Z)$.
The lattice $R$ is thus generated by the $E_t^{(k)}$ with
\begin{align*}
    \forall \, s,t\in\F^2_3, \, k,l\in\qty{1,2}\colon \qquad 
    \expval{ {E}^{(k)}_s \, , \, {E}^{(l)}_t  } = \begin{cases}
        \,\,\,\, 0  & \text{if } s\neq t \, ,\\
        \,\,\,\, 1 & \text{if } s=t, k\neq l \, , \\
        \,\,\,\, -2 & \text{if } s=t, k=l \, .
    \end{cases}
\end{align*}
Hence $R$ is a root lattice of type $A_2^9$.

Returning to our blowing up and down diagram \eqref{blowupanddown}, for 
every $t\in\F_3^2$ let $\CCC_t$ denote the image under $\widetilde\pi$ of 
the preimage $\widetilde\CCC_t$
of $F_t$ in $\widetilde T$, and let $E_t\in R$ denote
its Poincar\'e dual. Since $\widetilde\pi$ is a $3\colon 1$ cover which
is branched on the divisor $\sum_{t\in\F_3^2}\CCC_t$, we  conclude
that $\frac{1}{3}\sum_{t\in\F_3^2}E_t$ is an integral cohomology class. 
For every $t\in\F_3^2$
it follows that $\frac{1}{3}E_t$ belongs to the dual of the $A_2$-type lattice 
generated by $E_t^{(1)}$ and $E_t^{(2)}$, and in fact without loss of 
generality\footnote{Algebraically, it is natural to
choose $E^{(1)}_t$ and $E^{(2)}_t$ as simple roots of the $A_2$ root lattice that they 
generate.}
$$
\forall t\in\F_3^2\colon\qquad
E_t=E_t^{(1)}+2E_t^{(2)}\,,
$$
since $\CCC_t$ is effective and primitive.

Similarly, if $C$ is a closed $2$-cycle  in $T$ which
is mapped to itself by our $\Z_3$-action, let $\mu\in H^2(T,\Z)$ denote the Poincar\'e dual 
of $C$ and 
$S\subset \F_3^2$  the set of labels of all $\Z_3$-fixed points in $C$, and let
 $\widetilde\CCC$ denote the preimage of $C$ in $\widetilde T$.
Then $\widetilde\pi_{\mid\widetilde\CCC}$ is a 
$3\colon 1$ cover which is branched on  $\sum_{t\in S}\CCC_t$, 
thus yielding 
an unbranched covering when restricted to $\widetilde\CCC\setminus\bigcup_{t\in S}\widetilde\CCC_t$.
Hence the class $\frac{1}{3}\left( \pi_\ast\mu - \sum_{t\in S}E_t\right)$
is integral. Equivalently, one may argue that $C\longrightarrow C/\Z_3$ is a $3\colon1$ 
cover with branch points ${F}_t$, $t \in S$, such that
$\pi\left( C\setminus \{{F}_t \,\vert \, t \in S\}\,\right)
=\widetilde\pi(\widetilde\CCC)\setminus \bigcup_{t\in S}\CCC_t$ 
is the threefold of the 2-cycle on $X$ whose Poincar\'e dual  is the class
$\frac{1}{3}\left( \pi_\ast\mu - \sum_{t\in S}E_t\right)$.

Applying these observations to appropriate examples of cycles $C$, we
shall see that the following holds:

\begin{figure}
\begin{center}  
\includegraphics[width=7cm,keepaspectratio]{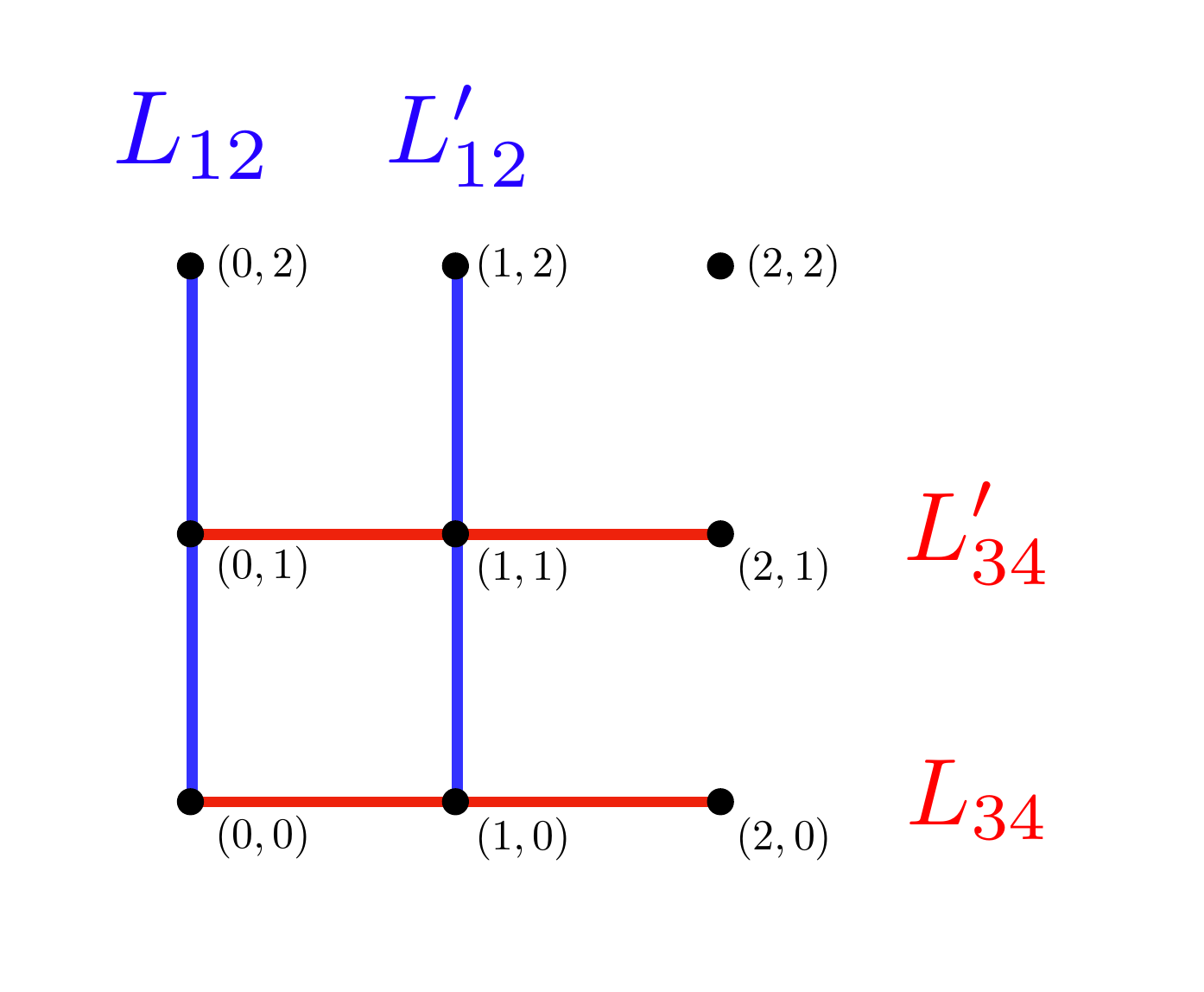}
\end{center}
\caption{The affine plane $\F_3^2$ with two pairs of affine parallel lines highlighted.}
\label{fig:affineF32}
\end{figure}

\blem\label{gluelist}%
We have four types  of two-forms 
that represent integral cohomology classes as follows:
\begin{enumerate}[label=\mbox{\rm(\roman*)}, itemsep=0pt, topsep=0pt]
\item 
Let $L_{12}\subset\F_3^2$ denote the affine line comprised of the three
points $(0,0),\, (0,1)$ and $(0,2)$, and $L_{12}^\prime\subset\F_3^2$
any affine line parallel to $L_{12}$ (see figure \ref{fig:affineF32} for illustration).
Then 
$$
\textstyle
\frac{1}{3}\pi_\ast \mu_{12} + \frac{1}{3}\sum\limits_{t\in L_{12}^\prime}E_t
$$
represents a class in $H^2(X,\Z)$.
\item 
Let $L_{34}\subset\F_3^2$ denote the affine line comprised of the three
points $(0,0),\, (1,0)$ and $(2,0)$, and $L_{34}^\prime\subset\F_3^2$
any affine line parallel to $L_{34}$.
Then 
$$
\textstyle
\frac{1}{3}\pi_\ast \mu_{34} - \frac{1}{3}\sum\limits_{t\in L_{34}^\prime}E_t
$$
represents a class in $H^2(X,\Z)$.
\item 
Let $L_{01}^\prime\subset\F_3^2$ denote an affine line 
of slope $(-1)$, that is, any affine line parallel to the line $L_{01}$
containing the three points $(0,0), (1,2)$ and $(2,1)$.
Then 
$$
\textstyle
\frac{1}{3}\pi_\ast (\mu_{12}-\mu_{34} + \mu_{14}+\mu_{23}-\mu_{24}) 
+ \frac{1}{3}\sum\limits_{t\in L_{01}^\prime}E_t
$$
represents a class in $H^2(X,\Z)$.
\item 
Let $L_{00}^\prime\subset\F_3^2$ denote an affine line 
of slope $(+1)$, that is, any affine line parallel 
to the line $L_{00}$
containing the three points $(0,0), (1,1)$ and $(2,2)$.
Then 
$$
\textstyle
\frac{1}{3}\pi_\ast (\mu_{12} -\mu_{34} + \mu_{14}+\mu_{23}-\mu_{13}) 
+ \frac{1}{3}\sum\limits_{t\in L_{00}^\prime}E_t
$$
represents a class in $H^2(X,\Z)$.
\end{enumerate}
\elem

\bpf
We apply the technique
explained above to appropriate closed 2-cycles $C$ in $T$:
\begin{enumerate}[label=\mbox{\rm(\roman*)}, itemsep=0pt, topsep=0pt]
\item 
Choose $C$ as a 2-cycle represented in $\C^2$ by
$\{z_1=c\}$ with $c\in \frac{\xi+2}{3}\Z$. 
We choose orientations such that
its Poincar\'e dual is represented by the $\Z_3$-invariant 2-form $-\mu_{12}=-\frac{2}{\sqrt{3}}( dx^1\wedge dx^2)$ 
(recall equation \eqref{dualbasis}).
If $c=0$, then the fixed points contained in $C$ are the 
$F_t$ with $t\in L_{12}$, and the claim follows. Choosing 
arbitrary $c\in \frac{\xi+2}{3}\Z$ amounts to replacing the affine line $L_{12}$
by any of its parallels $L_{12}^\prime$ in $\F_3^2$.
\item 
The argument is  analogous to that in part (i) with the cycle represented by
$\{z_1=c\}$, $c\in \frac{\xi+2}{3}\Z$, replaced by the one represented by $\{z_2=c\}$ and 
therefore $L_{12}^\prime$ replaced by $L_{34}^\prime$. That 
$\frac{1}{3}\pi_\ast \mu_{34} - \frac{1}{3}\sum_{t\in L_{34}^\prime}E_t$
is obtained rather than 
$\frac{1}{3}\pi_\ast \mu_{34} + \frac{1}{3}\sum_{t\in L_{34}^\prime}E_t$
follows from the fact that the resulting class must yield integers when
inserted into the intersection form with  the classes obtained in (i).
\item 
Choose $C$ as a 2-cycle represented by $\{z_2-\ol z_1=c\}$ with $c\in \frac{\xi+2}{3}\Z$. 
One checks that $C$ is mapped to itself by our $\Z_3$-action
and that the fixed points contained in $C$ are the 
$F_t$ with $t\in L_{01}^\prime$, an affine line in $\F_3^2$  parallel to the line
$L_{01}$.
Since on $C$, both $x^1-x^3$ and $x^2+x^4$ are constant, the
Poincar\'e dual of $C$ is, possibly up to a global sign,
\begin{eqnarray*}\textstyle
\textstyle \frac{2}{\sqrt{3}}(dx^1-dx^3)\wedge (dx^2+dx^4)
&\stackrel{\eqref{dualbasis}}{=}& \textstyle (\mu_1-\frac{1}{2}\mu_2 - \mu_3+\frac{1}{2}\mu_4) \wedge (\mu_2+\mu_4) \\
&=&
\mu_{12}-\mu_{34} + \mu_{14}+\mu_{23}-\mu_{24}\, .
\end{eqnarray*}
The claim follows when taking into account that the global sign
is already determined by the fact that the resulting  classes must yield
integers when inserted into the intersection form with the classes obtained in (i).
\item 
Choose $C$ as a 2-cycle represented by 
$\{z_2+\xi^2 \ol z_1=c\}$ with  $c\in \frac{\xi+2}{3}\Z$, and proceed as in (iii).
One checks that $C$ is mapped to itself by our $\Z_3$-action
and that the fixed points contained in $C$ are the 
$F_t$ with $t\in L_{00}^\prime$, an affine line in $\F_3^2$ parallel to the line
$L_{00}$.
Using the fact that on $C$, 
both $x^3-\frac{1}{2}x^1-\frac{\sqrt3}{2}x^2$ and $x^4-\frac{\sqrt3}{2}x^1+\frac{1}{2}x^2$ are constant, the
Poincar\'e dual of $C$ is found to be
\begin{align*}\textstyle
\textstyle\frac{2}{\sqrt3}(dx^3 - \frac{1}{2}dx^1 - \frac{\sqrt{3}}{2}dx^2)&
\textstyle
\wedge (dx^4 - \frac{\sqrt{3}}{2}dx^1+ \frac{1}{2}dx^2 )\\
&\stackrel{\eqref{dualbasis}}{=} \textstyle (\mu_3-\frac{1}{2}\mu_4 - \frac{1}{2}\mu_1-\frac{1}{2}\mu_2) \wedge (\mu_4-\mu_1+\mu_2)\\
&=
 \mu_{34}-\mu_{12}+\mu_{13}-\mu_{14}-\mu_{23}\, ,
\end{align*}
and the claim follows.
\vspace*{-1em}
\end{enumerate}
\epf

What we have shown so far allows us to determine the lattice $P$:

\bprop\label{Pconstruct}%
The lattice $P$ is generated by the root lattice $R$ together with the vectors
of the form
$$
\textstyle
\frac{1}{3} \sum\limits_{t\in L^\prime} E_t - \frac{1}{3} \sum\limits_{t\in L^{\prime\prime}} E_t\, ,
$$
where $L^\prime,\, L^{\prime\prime}$ are parallel affine lines in $\F_3^2$.\\
The lattice $P^\ast$ is generated by the root lattice $R$ together with the vectors
of the form
$$
\textstyle
\frac{1}{3} \sum\limits_{t\in L^\prime} E_t \, ,
$$
where $L^\prime$ is an affine line in $\F_3^2$.\\
The discriminant group of $P$ is isomorphic to $(\Z_3)^3$ with generators
$p_j+P$, $j\in\{1,2,3\}$, where 
\begin{eqnarray*}
p_1&:=& \textstyle\frac{1}{3} \sum\limits_{t\in L_{34}} E_t + \frac{1}{3} \sum\limits_{t\in L_{00}} E_t\, ,\\
p_2&:=& \textstyle\frac{1}{3} \sum\limits_{t\in L_{12}} E_t + \frac{1}{3} \sum\limits_{t\in L_{34}} E_t\, ,\\
p_3&:=& \textstyle\frac{1}{3} \sum\limits_{t\in L_{12}} E_t + \frac{1}{3} \sum\limits_{t\in L_{00}} E_t\, .
\end{eqnarray*}
With respect to these generators, the bilinear form $b_P$ 
with values in $\Q/\Z$ is diagonal, and the discriminant form takes values 
$q_P(p_j+P)= \frac{2}{3} \mod 2\Z$ for $j\in\{1,2,3\}$.
\eprop

\bpf
Recall first that $R\subset P$ are even lattices, such 
that $R\subset P\subset P^\ast\subset R^\ast$, where $R^\ast$ is generated by $R$
and the vectors $\frac{1}{3} E_t$ with $t\in\F_3^2$. Both $P$ and $P^\ast$ can therefore
be generated by $R$ together with appropriate vectors of the form 
$\frac{1}{3}\sum_{t\in\F_3^2} \eps_t E_t$,
where $\eps_t\in\{0,+1,-1\}$ for all $t$.

Note that the difference of any two vectors listed as type (i) in lemma \ref{gluelist} 
has the form $\frac{1}{3} \sum_{t\in L^\prime} E_t - \frac{1}{3} \sum_{t\in L^{\prime\prime}} E_t$,
where $L^\prime,\, L^{\prime\prime}$ are  affine lines in $\F_3^2$ that are parallel to
$L_{12}$. Proceeding similarly with the other types of vectors listed in lemma \ref{gluelist},
we confirm that all claimed generators of $P$ are indeed 
contained in $P$.
Note moreover that all vectors listed in lemma \ref{gluelist} have the form
$\frac{1}{3}(\pi_\ast\mu+e)$ with $\mu\in H^2(T,\Z)^{\Z_3}$ and $e\in R$.
Since $\frac{1}{3}(\pi_\ast\mu+e)\in H^2(X,\Z)$ with $K\oplus P\subset H^2(X,\Z)$,
it follows that $\frac{1}{3}\pi_\ast\mu\in K^\ast$ and $\frac{1}{3}e\in P^\ast$. In 
fact, all the vectors listed in lemma \ref{gluelist} are 
glue vectors for 
the primitive sublattices $K$ and $P$ of $H^2(X,\Z)$ according to Nikulin's gluing technique
described in appendix \ref{subapp:glue}. It in particular follows that all the 
claimed generators of $P^\ast$ are indeed contained in $P^\ast$.

Assume now that $p\in P$ with $p=\frac{1}{3}\sum_{t\in\F_3^2} \eps_t E_t$,
where $\eps_t\in\{0,+1,-1\}$ for all $t\in\F_3^2$. We need to show that $p$ is an
integral linear combination of elements of $R$ and vectors of the form 
$\frac{1}{3} \sum_{t\in L^\prime} E_t - \frac{1}{3} \sum_{t\in L^{\prime\prime}} E_t$
with parallel affine lines $L^\prime,\, L^{\prime\prime}$ in $\F_3^2$.
By adding vectors of the latter form with affine lines parallel to $L_{12}$
and  $L_{34}$, we may assume without loss of generality that 
$\eps_{00}=\eps_{10}=\eps_{01}=0$. Since we already know that for any affine line
$L^\prime$ in $\F_3^2$, we have $\frac{1}{3} \sum_{t\in L^\prime} E_t\in P^\ast$,
the intersection form of this vector with $p$ must yield an integer. In other words,
for every affine line $L^\prime$ in $\F_3^2$, we find that 
$\sum_{t\in L^\prime} \eps_t$ is divisible by $3$. Therefore $\eps_{00}=\eps_{10}=\eps_{01}=0$
implies $\eps_{02}=\eps_{20}=0$ (using $L^\prime=L_{12}$ and $L^\prime=L_{34}$), 
and similarly $\eps_t=0$ for all $t\in\F_3^2$. This completes the proof of the 
claim concerning the generators of $P$.

All the remaining claims follow similarly by a direct calculation, where it is helpful to observe
that $p_1-p_2-p_3+P=\frac{1}{3} \sum_{t\in L_{12}} E_t+P$ etc.
\epf
\bigskip

For later convenience we remark that the 
following three vectors yield a minimal set of generators of $P$ from the 
root sublattice $R$:
\begin{equation}\label{eq:v1v2v3}\textstyle
v_1:=\frac{1}{3}\sum\limits_{t \in L_{12}}E_t-\frac{1}{3}\sum\limits_{t \in (1,0)+L_{12}}E_t\,,\quad
v_2:=\frac{1}{3}\sum\limits_{t\in L_{34}}E_t-\frac{1}{3}\sum\limits_{t\in (0,1)+L_{34}}E_t\, ,\quad
v_3:=\frac{1}{3}\sum\limits_{t \in\F_3^2}E_t\, .
\end{equation}
To see this, note first of all that proposition \ref{Pconstruct} immediately implies
$v_1,\, v_2\in P$, and $v_3\in P$ follows from 
$$
\textstyle
v_3=
\underbrace{\textstyle v_1-\frac{1}{3}\Bigl( \sum\limits_{t \in (1,0)+L_{12}} E_t-\sum\limits_{t\in (2,0)+L_{12}} E_t 
\Bigr)}_{\in P}
+\underbrace{\textstyle \sum\limits_{t\in (1,0)+L_{12}}E_t}_{\in P}\, .
$$
Moreover, one checks by a direct calculation that  ($v_1+R$, $v_2+R$, $v_3+R$) is linearly independent
in $P/R$, generating $(\Z_3)^3$. It thus suffices to show that $R$ has index $3^3$ in $P$. 
By proposition \ref{Pconstruct}
we have ${\rm disc}(P)=3^3$, and  since $ A_2^\ast/A_2\cong\Z_3$,  
we have    ${\rm disc}(R)=3^9$. Hence $\vert P/ R\vert=3^3$ according to \eqref{indexdisc},
as required.
\bigskip

Proposition \ref{Pconstruct} implies that the lattice $P^\ast$ is generated by roots. It 
will prove useful to know that on $P^\ast\setminus\{0\}$, the absolute value of the quadratic form
is minimal  for these roots:

\blem\label{lemrootsinPstar}%
If $p\in P^\ast$ with $p\neq0$, then $|\langle p,p\rangle| \geq 2$.
\elem

\bpf
Assume that $p\in P^\ast$ obeys $|\langle p,p\rangle| < 2$; we need to show that $p=0$.

To  this end, we decompose $p$ as $p=\sum_{t\in\F_3^2} p_t$, where for every $t\in\F_3^2$,
$p_t$ belongs to the lattice of type $A_2^\ast$ generated by $E_t^{(2)}$ and $\frac{1}{3} E_t$. 
Assume now that for some $t\in\F_3^2$, we have $p_t\neq 0$. 
Since $|\langle p,p\rangle| < 2$ by assumption, it follows that  
$p_t$ does not belong to the $A_2$-type root lattice generated by $E_t^{(1)}$ and $E_t^{(2)}$.
On the other hand, we have $|\langle p_t,p_t\rangle |\geq \frac{2}{3}$, since in the lattice dual to 
an $A_2$-type lattice, the minimum of the absolute value of the quadratic form on non-zero
vectors is $\frac{2}{3}$ (see appendix \ref{subapp:niemeier}).

Again using our assumption $|\langle p,p\rangle| < 2$,
we may conclude that there exist $t_1,\, t_2\in\F_3^2$ such that
$$
p = \frac{\eps_1}{3} E_{t_1} + \frac{\eps_2}{3} E_{t_2} + r\, ,\qquad
\eps_1,\, \eps_2\in\{ 0, \pm1\}, \; r\in R\, .
$$
We may therefore choose parallel affine lines $L^\prime,\, L^{\prime\prime}$ in $\F_3^2$
with $L^\prime\neq L^{\prime\prime}$, $t_1\in L^\prime$ and $t_2\not\in L^\prime\cup L^{\prime\prime}$.
By proposition \ref{Pconstruct}, $\wt p:=\frac{1}{3}\left( \sum_{t\in L^\prime} E_t
-\sum_{t\in L^{\prime\prime}} E_t\right)$ yields $\widetilde p\in P$, hence 
$\langle p,\widetilde p\rangle\in\Z$. Since up to adding an integer,
$\langle p,\widetilde p\rangle$ is given by 
$\frac{\eps_1}{9} \langle E_{t_1},E_{t_1}\rangle = -\frac{2\eps_1}{3}$, we find $\eps_1=0$.
Analogously, one shows $\eps_2=0$, from which $r=0$ follows since $|\langle p,p\rangle| < 2$
by assumption. Altogether we have shown $p=0$, as desired.
\epf

\bigskip

We claim that the above also suffices to determine the lattice $K$:

\bprop\label{Kconstruct}%
A basis for the lattice $K$ is given by $(\kappa_1, \kappa_2, \kappa_3, \kappa_4)$, 
where 
$$
\begin{aligned}
    \kappa_1 &:=\textstyle \tfrac{1}{3}\pi_*(\mu_{14}+\mu_{23}-2\mu_{13}+\mu_{24}) \, , \quad  &\kappa_3 &:= \pi_*\mu_{12} \, ,\\
    \kappa_2 &:=\textstyle \tfrac{1}{3}\pi_*(-\mu_{14}-\mu_{23}-\mu_{13}+2\mu_{24}) \, , \quad &\kappa_4 &:= \pi_*\mu_{34}\, ,
\end{aligned} 
$$
such that with respect to this basis, the bilinear form on $K$ is given by
$$
\left(\begin{matrix}2&1&0&0\\1&2&0&0\\0&0&0&3\\0&0&3&0\end{matrix}\right)\; .
$$
The discriminant group of $K$ is isomorphic to $(\Z_3)^3$ with generators
$k_j+K$, $j\in\{1,2,3\}$, where 
\begin{align*}
k_1&:=\textstyle \frac{1}{3} (2\kappa_1-\kappa_2+\kappa_3+\kappa_4)
&=&\; \textstyle \frac{1}{3}\pi_\ast(\mu_{12}+\mu_{34}-\mu_{13}+\mu_{14}+\mu_{23})\, ,\\
k_2&:=\textstyle \textstyle \frac{1}{3} (\kappa_3-\kappa_4)
&=&\; \textstyle \frac{1}{3}\pi_\ast (\mu_{12}-\mu_{34})\, ,\\
k_3&:=\textstyle \textstyle \frac{1}{3} (2\kappa_1-\kappa_2-\kappa_3-\kappa_4)
&=&\; \textstyle \frac{1}{3}\pi_\ast(-\mu_{12}-\mu_{34}-\mu_{13}+\mu_{14}+\mu_{23})\, .
\end{align*}
With respect to these generators, the bilinear form $b_K$ 
with values in $\Q/\Z$ is diagonal, and the discriminant form takes values 
$q_K(k_j+K)= -\frac{2}{3} \mod 2\Z$ for $j\in\{1,2,3\}$.
\eprop

\bpf
Let us first show that $(\kappa_1, \kappa_2, \kappa_3, \kappa_4)$
indeed yields a basis of the lattice $K$, which by definition is the smallest primitive sublattice of
$H^2(X,\Z)$ that contains $\pi_\ast\left(H^2(T,\Z)^{\Z_3}\right)$. Equivalently, we may define $K$
as the orthogonal complement of $P$ in $H^2(X,\Z)$. Since 
$\kappa_1+\kappa_2=\pi_*(\mu_{24}-\mu_{13})$ and therefore 
$3\kappa_1,\kappa_1+\kappa_2, \kappa_3,\kappa_4 \in \pi_\ast\left(H^2(T,\Z)^{\Z_3}\right)$, while 
$\kappa_1\not\in \pi_\ast\left(H^2(T,\Z)^{\Z_3}\right)$, it suffices 
to show that $\kappa_1\in K$ and that $\pi_\ast\left(H^2(T,\Z)^{\Z_3}\right)$ is an index $3$ sublattice
of $K$.

To prove that $\kappa_1\in K$, consider the following vectors, which according to 
lemma \ref{gluelist} belong to $H^2(X,\Z)$ with $g_1, \ldots, g_4$ of types (i), \ldots, (iv),
and $g_4^\prime$ of type (iv)
respectively:
\begin{eqnarray*}
g_1&:=& \textstyle
\frac{1}{3}\pi_\ast \mu_{12} + \frac{1}{3}(E_{00}+E_{01}+E_{02})\, ,\\
g_2&:=&\textstyle \frac{1}{3}\pi_\ast \mu_{34} - \frac{1}{3}(E_{01}+E_{11}+E_{21})\, ,\\
g_3&:=&\textstyle \frac{1}{3}\pi_\ast(\mu_{12}-\mu_{34} + \mu_{14}+\mu_{23}-\mu_{24})  + \frac{1}{3}(E_{01}+E_{10}+E_{22})\, ,\\
g_4&:=&\textstyle \frac{1}{3}\pi_\ast(\mu_{12} -\mu_{34} + \mu_{14}+\mu_{23}-\mu_{13})  + \frac{1}{3}(E_{00}+E_{11}+E_{22})\, ,\\
g_4^\prime&:=&
\textstyle \frac{1}{3}\pi_\ast(\mu_{12} -\mu_{34} + \mu_{14}+\mu_{23}-\mu_{13})  + \frac{1}{3}(E_{02}+E_{10}+E_{21})\, .
\end{eqnarray*}
Hence the following vector belongs to $H^2(X,\Z)$,
\[
-g_1+g_2-g_3+g_4+g_4^\prime=\tfrac{1}{3}\pi_*(\mu_{14}+\mu_{23}-2\mu_{13}+\mu_{24})-E_{01}
=\kappa_1 -E_{01}\, ,
\]
implying $\kappa_1\in K$ as claimed.

To prove that $\pi_\ast\left(H^2(T,\Z)^{\Z_3}\right)$ is an index $3$ sublattice
of $K$, note that by lemma \ref{H2Tinv} we know that $\pi_\ast\left(H^2(T,\Z)^{\Z_3}\right)$
has discriminant $-3^5$, hence by \eqref{indexdisc} we need to show that $|\disc(K)|=3^3$.
Since $K$ and $P$ are orthogonal complements of one another in the even, unimodular lattice
$H^2(X,\Z)$,  Nikulin's gluing construction as discussed in appendix \ref{subapp:glue}
implies that $K$ and $P$ have isomorphic discriminant groups, that is,
$K^\ast/K\cong P^\ast/P\cong (\Z_3)^3$ by proposition \ref{Pconstruct} and 
$|\disc(K)|=3^3$ as required.
Using lemma \ref{H2Tinv}, one also confirms
that $k_j\in K^\ast$ for $j\in\{1,\, 2,\, 3\}$, and by direct calculation one checks
that the $k_j+K$ generate the discriminant group with values of the discriminant forms
as claimed.
\epf
\bigskip

The above  now implies:

\bprop\label{gluePandK}%
Consider the group homomorphism $\gamma\colon K^\ast/K\longrightarrow P^\ast/P$ which
is given by $\gamma(k_j+K)=p_j+P$ for $j\in\{1,2,3\}$. Then $\gamma$ is an isomorphism with
$q_P\circ\gamma=-q_K$,
and for any pair $(k,p)\in K^\ast\oplus P^\ast$, we have $k+p\in H^2(X,\Z)$ if and
only if $\gamma(k+K)=p+P$. In other words, $\gamma$ yields the gluing prescription
to glue $H^2(X,\Z)$ from the primitive sublattices 
$K$ and $P$, induced by the underlying geometry.
\eprop

\bpf
Note that propositions \ref{Pconstruct} and \ref{Kconstruct} imply that $\gamma$ yields an 
isomorphism between the discriminant groups of $K$ and $P$ which obeys 
$q_P\circ\gamma=-q_K$. Using lemma \ref{gluelist} one now checks that 
$k_j+p_j\in H^2(X,\Z)$ for $j\in\{1,2,3\}$, which suffices to complete
the proof of proposition \ref{gluePandK}.
\epf
\bigskip

We conclude this subsection by describing, 
in terms of the lattices $K$ and $P$, 
the transcendental lattice $T_X\subset H^2(X,\Z)$
and the K\"ahler class $\omega$ that is induced by the standard Euclidean metric on 
$T$, as well as all possible  choices of K\"ahler classes on $X$ that  descend from $\Z_3$-invariant K\"ahler
classes on the complex torus $T$.
The general facts concerning the complex geometry of K3 surfaces are 
classical, and for a summary we refer, for example, to \cite[\S3]{tawe11}.
Recall that the transcendental lattice $T_X$
is obtained as orthogonal complement in $H^2(X,\Z)$ of the Neron-Severi lattice $NS(X)$,
which in turn is obtained as $NS(X)=\Omega^\perp\cap H^2(X,\Z)$ from the two-dimensional
subspace $\Omega\subset H^2(X,\R)$ generated by the classes represented by the
real and imaginary parts of the 
holomorphic volume form of $X$. Since $H^2(X,\R)=H^2(X,\Z)\otimes_\Z \R$,
we may express all these quantities in terms of the lattices $K$ and $P$
from which we have glued $H^2(X,\Z)$. 

The holomorphic volume form of our $\Z_3$-orbifold K3 surface
$X$ is $\pi_\ast (dz_1\wedge dz_2)$, 
where $\pi_\ast$ is $\C$-linearly extended from 
$\pi_\ast\colon H^2(T,\Z)^{\Z_3}\longrightarrow H^2(X,\Z)$.
Hence $\Omega$ is generated by the two vectors
\begin{eqnarray*}
\pi_\ast (dx^1\wedge dx^3-dx^2\wedge dx^4)&\stackrel{\eqref{dualbasis}}{=}& 
\textstyle
\pi_\ast\left((\mu_1-\frac{1}{2}\mu_2)\wedge(\mu_3-\frac{1}{2}\mu_4) - \frac{3}{4}\mu_2\wedge\mu_4\right) \\
&=&
\textstyle 
\frac{1}{2} 
\pi_\ast( 2\mu_{13} - \mu_{23} - \mu_{14} - \mu_{24} ) \stackrel{\rm{}Prop.\, \ref{Kconstruct}}{=} -\frac{3}{2} \kappa_1\, ,\\
\pi_\ast (dx^1\wedge dx^4+dx^2\wedge dx^3)&\stackrel{\eqref{dualbasis}}{=}& 
\textstyle
\frac{\sqrt3}{2}\pi_\ast \left((\mu_1-\frac{1}{2}\mu_2)\wedge\mu_4 + \mu_2\wedge(\mu_3-\frac{1}{2}\mu_4)\right) \\
&=&
\textstyle 
\frac{\sqrt3}{2}
\pi_\ast ( \mu_{14} + \mu_{23} - \mu_{24} ) \stackrel{\rm{}Prop.\, \ref{Kconstruct}}{=} 
\frac{\sqrt3}{2} ( \kappa_1 - 2 \kappa_2)\, .\\
\end{eqnarray*}
It follows that $T_X$ is the lattice that is generated by $\kappa_1$ and $\kappa_2$, and
 by proposition \ref{Kconstruct} its quadratic form with respect to this basis is $\smqty(2&1\\1&2)$.
This was already proven in \cite[Lemma\,5.1]{shino77}, using different methods.

As explained at the beginning of subsection \ref{subsec:Z3orbi}, the $\Z_3$-orbifold construction is 
carried out independently of the choice of a K\"ahler class, and it yields the K3 surface that we denote by $X$.
However, as we shall recall in subsection \ref{subsec:Symmetries1}, our notion of symmetries requires the choice of a
K\"ahler class on $X$.
The standard K\"ahler form  on $X$, which is induced by the standard Euclidean metric on $\C^2$, is 
$$
\begin{array}{rclcl}
\textstyle
\omega&\sim&
\frac{i}{2}\pi_\ast (dz_1\wedge d\ol z_1+ dz_2\wedge d\ol z_2)
&=& \pi_\ast (dx^1\wedge dx^2+dx^3\wedge dx^4)\\
&\stackrel{\eqref{dualbasis}}{=}&
\textstyle
\frac{\sqrt3}{2}\pi_\ast\left((\mu_1-\frac{1}{2}\mu_2)\wedge\mu_2
+ (\mu_3-\frac{1}{2}\mu_4)\wedge\mu_4 \right)
&=& \frac{\sqrt3}{2}\pi_\ast ( \mu_{12}+\mu_{34}) \\
&\stackrel{\rm{}Prop.\, \ref{Kconstruct}}{=}&   \frac{\sqrt3}{2}( \kappa_3+\kappa_4)\, .
 \end{array}
 $$
In keeping with the notion of \textsl{$\Z_3$-orbifold limits}, for our choice of K\"ahler class on $X$
we allow any class that descends from a $\Z_3$-invariant K\"ahler class on the underlying torus. 
In order to obtain a 
(degenerate) K\"ahler class compatible with our choice of complex structure, we may choose any
$\wt\omega\in\pi_\ast\left(H^2(T,\Z)^{\Z_3}\right)$ which yields a real cohomology class in $H^{1,1}(X,\C)$ and
which encodes a positive definite Hermitean sesquilinear form. With respect to the intersection form on $H^{2}(X,\C)$,
requiring $\wt\omega\in H^{1,1}(X,\C)$ is equivalent to requiring $\wt\omega\in\left(T_X\right)^\perp\cap \pi_\ast\left(H^2(T,\Z)^{\Z_3}\right)$. 
Such $\wt\omega$ encodes a positive definite Hermitean sesquilinear form if and only if it has positive self-intersection.
In conclusion, we have shown

\bprop\label{hyperk}%
The complex structure of the $\Z_3$-orbifold K3 surface $X$ is uniquely determined
 by the two dimensional subspace $\Omega$ of $H^2(X,\R)$ with basis
$(\kappa_1,\kappa_2)$. Its standard K\"ahler form is 
$\omega:=\kappa_3+\kappa_4$, and if $\wt\omega$ represents a (degenerate) K\"ahler class on $X$ 
that descends from a $\Z_3$-invariant K\"ahler class on $T$ then $\wt\omega=V_3\kappa_3+V_4\kappa_4$
with $V_3,V_4\in\R$ and $V_3V_4>0$.
\eprop

With notations as in proposition \ref{hyperk}, $\wt\omega=V_3\kappa_3+V_4\kappa_4$ is induced by a 
K\"ahler class on $T$ which yields this complex torus as the product of two elliptic curves with complex
structure moduli $\xi$ and volumes given by $|V_3|,|V_4|$, respectively, up to a common factor that 
depends on the conventions of normalization for the respective volume forms.

\subsection{Symmetries of \texorpdfstring{$\Z_3$}{TEXT}-orbifold K3s}\label{subsec:Symmetries1}

In this subsection, we determine the symmetry group of $X=\widetilde{T/\mathbb{Z}_3}$ 
and its induced action on the integral cohomology lattice $H^2(X,\mathbb{Z})$, which as we shall recall below 
completely specifies the symmetry group. As was mentioned above, we call a map $f\colon X\longrightarrow X$ a
\textsc{symmetry of $X$} if and only if  it is a biholomorphic isometry that leaves the holomorphic 
volume form of $X$ invariant. This is in accord with the definition given in \cite[Def. 3.2.1]{tawe11}, 
and for further details on the relevant notions from complex geometry we refer the reader to that reference.
Note that here, $f$ is called an isometry if and only if it leaves invariant a choice $\wt\omega$ of 
K\"ahler form  which is induced by a $\Z_3$-invariant K\"ahler 
metric on the complex torus $T$. 
As such, $\wt\omega$  represents a degenerate K\"ahler class, that is,
a class on the boundary of the K\"ahler cone.
We  refer to the \textsc{symmetries of  $(X,\wt\omega)$}  to remind the reader of our choice of 
degenerate K\"ahler class $\wt\omega$.
The integral cohomology lattice as well as the complex structure and K\"ahler class of $X$ were
 described in the previous subsection in terms of the  sublattices $K$ and $P$. 
We shall see that the symmetry group of $X$
consists of the symmetries which are induced by the symmetries of the underlying 
torus;  hence it corresponds 
to a group of lattice automorphisms that 
map $K$ to $K$ and $P$ to $P$.\footnote{This is in contrast
to the situation for Kummer surfaces, where
there may exist symmetries that do
not map the analogs of the lattices $K$ and $P$ to themselves. 
Note however that our requirement that symmetries fix the degenerate K\"ahler
class induced from a $\Z_3$-invariant K\"ahler metric on the torus, 
rules out,  for example,  the `reflections'
investigated in \cite[Prop.~5.1]{ko98} and which were discovered by Keum \cite{ke97}: 
these reflections do not map the analogs of $K$ and $P$
to themselves, but they also do not fix the required K\"ahler class.}
\bigskip

We first clarify the role of lattice automorphisms for the description of symmetries of $(X,\omega)$. 
To simplify the exposition, we begin by focussing on the standard K\"ahler class $\omega$ of proposition \ref{hyperk},
while in remark \ref{generalsymmetries} we address the general case $\wt\omega$ with $V_3,\, V_4\in\R$, $V_3V_4>0$.
Let $f\colon X\longrightarrow X$ denote a symmetry of $X$, and $\Phi:=f_\ast$ the induced 
map on the cohomology of $X$, such that by construction $\Phi\in\Aut(H^2(X,\Z))$. 
Here and in the following, we use the same letter $\Phi$ for the linear extension of $\Phi$ to 
$H^2(X,\R)$ or any subset of this vector space.
By definition, $\Phi$ leaves invariant the cohomology classes which
specify the complex structure and K\"ahler class of $X$, that is, according to proposition \ref{hyperk},
\begin{equation}\label{Sigmadef}%
\mbox{with }\quad \Sigma:=\mbox{span}_\R \left\{ \kappa_1,\, \kappa_2, \, \kappa_3+\kappa_4\right\}, \qquad
\forall\sigma\in\Sigma, \;\Phi(\sigma)=\sigma\, .
\end{equation}
Hence $\Phi$ induces an automorphism of the negative definite lattice $\Sigma^\perp\cap H^2(X,\Z)$.
Note that by the results of the previous subsection, $\Sigma^\perp\cap H^2(X,\Z)$ is generated by $P$
and the glue vector 
\begin{equation}\label{Kaehlerglue}%
\textstyle
\frac{1}{3} (\kappa_3-\kappa_4) + \frac{1}{3} \sum\limits_{t\in L_{12}} E_t + \frac{1}{3} \sum\limits_{t\in L_{34}} E_t
=k_2+p_2\, 
\end{equation}
with $p_2$ and $k_2$ as in propositions \ref{Pconstruct} and \ref{Kconstruct}.
In particular, if $v\in \Sigma^\perp\cap H^2(X,\Z)$ and $v\not\in P$, then $|\langle v,v\rangle|>2$,
which implies
\be\label{rootsareinP}
\forall v\in \Sigma^\perp\cap H^2(X,\Z)\colon\qquad
\langle v,v\rangle = -2\quad \Longrightarrow\quad v\in P
\, .
\ee

That is, the roots in $\Sigma^\perp\cap H^2(X,\Z)$ are given by those in the root lattice
$R\subset P$ of type $A_2^9$.
Since $\Phi$ is induced by $f$, it maps roots to roots, and in addition it
has a property known as \textsc{effectiveness}, which in our setting is equivalent to the fact 
that $\Phi$ maps each of the (``simple'') roots $E_t^{(j)}$ in $R$, where $t\in\F_3^2$ and $j\in\{1,2\}$, 
to a \emph{non-negative} linear combination of such simple roots, see e.g.~\cite[\S VIII.3]{bpv84}  
for further details. We say that a lattice automorphism of 
$H^2(X,\Z)$ is \textsc{$\Z_3$-effective} if 
it has this property and it acts trivially on $\Sigma$. 
In summary, we have argued
that every symmetry of $X$ induces a $\Z_3$-effective lattice automorphism of $H^2(X,\Z)$.
The Torelli Theorem for K3 surfaces as obtained from \cite[Thms.~2.7' \& 4.3]{ni80b}
in conjunction with \cite{ku77} and \cite[Prop. VIII.3.10]{bpv84} 
states that the converse is also true:
\be\label{Torelli}
\begin{array}{rl}
\Phi\in\Aut( H^2(X,\Z) ) & \mbox{is $\Z_3$-effective }\\
\Longleftrightarrow
&\mbox{ there exists a  symmetry $f$ of  $(X,\omega)$ with } \Phi=f_\ast \; ;\\
\mbox{ then:}
& f \mbox{ is uniquely determined by } \Phi=f_\ast \, .
\end{array}
\ee
\bigskip

We next determine all symmetries of  $(X,\omega)$ that are induced by symmetries of the torus $T$.
Analogously to the notion of symmetries of K3 surfaces, a symmetry
of $T$ is a biholomorphic isometry $g\colon T\longrightarrow T$ which leaves the 
holomorphic volume form of $T$ invariant. Here, $T=\C^2/L$ is equipped with the 
standard Euclidean metric, which
is induced by the standard Euclidean metric on $\C^2$.
For any $t_0\in T$, the translation $\alpha_{t_0}\colon T\longrightarrow T$,
$t\mapsto t+t_0$ yields such a symmetry, called \textsc{translational}. 
In addition, $T$ possesses \textsc{rotational} symmetries which are induced by the
action of a finite subgroup of $\SU(2)$ on $\C^2$ descending to $T=\C^2/L$. All such symmetry
groups have been classified by Fujiki in \cite{fu88}. The resulting 
group for our specific torus $T$ is $T\rtimes\m D$, where the additive group $T$ 
acts by translational symmetries and
$\m D\subset \SU(2)$ 
is the binary dihedral group of order $12$ which is generated by $\varrho$ and $\beta$,
\be\label{eq:beta}
\varrho\colon (z_1,z_2) \longmapsto (\zeta z_1, \zeta^{-1} z_2) \mbox{ with }\zeta:=\exp(\frac{2\pi i}{6})\, ,
\quad
\beta\colon (z_1,z_2) \longmapsto (- z_2,  z_1)\,  ,
\ee
such that $\varrho^3=\beta^2=-\mbox{id}$ and $\beta=\rho\beta\rho$. Since $\rho^2$ yields the map \eqref{Z3action},
a symmetry $g$ of $T$ descends to $T/\Z_3$ if and only if $g[z]:=[gz]$ is well defined for every
$[z]\in T/\Z_3$. In other words, we must have $g\Z_3 g^{-1} = \Z_3$ for the group $\Z_3$ that is 
generated  by $\rho^2$. We thus have shown that the subgroup of $T \rtimes\m D$ whose action descends
to $T/\Z_3$ is the normalizer of $\Z_3$. A direct calculation shows that the latter is given by 
$T^{\Z_3}\rtimes\m D$, where $T^{\Z_3}$ is the group whose elements are the fixed points of the $\Z_3$ action on the torus. Since $T^{\Z_3}\cong(\Z_3)^2$ by the discussion in subsection
\ref{subsec:Z3orbi}, and because $\varrho^2$ acts trivially on 
$T/\Z_3$, while $\m D/\Z_3\cong\Z_4$ is generated by $\beta$, we find an action of 
$(\Z_3)^2\rtimes \Z_4$ on $T/\Z_3$.

Next we will show that every translational symmetry $\alpha_{t_0}$
with $t_0\in T^{\Z_3}$ induces a well-defined symmetry on the minimal resolution $X$
of $T/\Z_3$, as does $\beta$, proving
that the group of symmetries of $X$ that are induced by symmetries of $T$ is
$(\Z_3)^2\rtimes \Z_4$. According to the Torelli Theorem \eqref{Torelli} it suffices
to describe these symmetries
in terms of their induced action on $H^2(X,\Z)$, which by construction maps $K$ to $K$
and $P$ to $P$. Since every vector in $P$ is a rational linear combination of roots
$E_t^{(j)}$ with $t\in\F_3^2$ and $j\in\{1,2\}$, it suffices to state the action on 
the generators $(\kappa_1,\kappa_2,\kappa_3,\kappa_4)$ of $K$ as well as the $E_t^{(j)}$.
\bigskip

We choose generators $\alpha^1, \alpha^2$ for the action of the translational 
group $(\Z_3)^2$, 
where
\be\label{eq:gammas}
        \begin{aligned}
            \alpha^1 &: (z_1,z_2) \longmapsto (z_1+\tfrac{\xi+2}{3},z_2) \, ,\\
            \alpha^2 &: (z_1,z_2) \longmapsto (z_1,z_2+\tfrac{\xi+2}{3}) \, .
        \end{aligned} 
\ee
These maps leave the differential forms $dz_1,\, dz_2$ invariant, thus
inducing a trivial action on $H^2(T,\Z)$ and therefore on 
$\pi_\ast\left(H^2(T,\Z)^{\Z_3}\right)$ and  on $K$. 

The local description of the blow-ups and downs in \eqref{blowupanddown} 
(see appendix \ref{app:A2blowup}) shows that the $\alpha^k$ induce actions
on $\widetilde T$ and $\widetilde X$ as well, where 
$\CCC_t\mapsto\CCC_{\widetilde\alpha^k(t)}$ for $t\in\F_3^2$ with 
$\widetilde\alpha^1(t_1,t_2)=(t_1+1,t_2)$ and $\widetilde\alpha^2(t_1,t_2)=(t_1,t_2+1)$.
On cohomology, we therefore have $\alpha^k_\ast (E_t)=E_{\widetilde\alpha^k(t)}$
and thus altogether
\blem\label{translationalaction}%
The action of the translational symmetries of $(X,\omega)$ on $H^2(X,\Z)$ is generated by
$\alpha^1_\ast$ and $\alpha^1_\ast$, where 
$\alpha^1_\ast (E_{(t_1,t_2)}^{(j)})=E_{(t_1+1,t_2)}^{(j)}$,
$\alpha^2_\ast (E_{(t_1,t_2)}^{(j)})=E_{(t_1,t_2+1)}^{(j)}$ when
$(t_1,t_2)\in\F_3^2$ and $j\in\{1,2\}$, and $\alpha^1_\ast(\kappa_j)=\alpha^2_\ast(\kappa_j)=\kappa_j$
when  $j\in\{1,\ldots,4\}$.
\elem
Let us confirm that the action determined above yields lattice 
automorphisms of  $P$ and that it
is compatible with the gluing 
prescription of proposition \ref{gluePandK}. 
For the generators \eqref{eq:v1v2v3} of $P$, we have
$$
\begin{aligned}
    \alpha^1_\ast\colon& \quad v_1\longmapsto v_1-v_3 + \sum_{t\in(1,0)+L_{12}} E_t\, , \quad &&v_2\longmapsto v_2 \, , \quad &v_3\longmapsto v_3 \, , \\
    \alpha^2_\ast\colon&\quad v_2\longmapsto v_2-v_3 + \sum_{t\in(0,1)+L_{34} }E_t\, , \quad   
   &&  v_1\longmapsto v_1 \, , \quad    
    &v_3\longmapsto v_3 \, .
\end{aligned}
$$
Hence our choice of generators of $P$ from $R$ is mapped to another set of such generators,
which implies that $\alpha^1_\ast$ and $\alpha^2_\ast$ yield lattice automorphisms of $P$. 
Similarly, all the
glue vectors $k_j+p_j$ found in proposition \ref{gluePandK} are mapped to glue 
vectors, since the $k_j$ are invariant and the affine lines  that are summed
over in the $p_j$ are mapped to parallel affine lines, such that $k_j+p_j$
differs from its image by vectors of the types listed in proposition 
\ref{Pconstruct}.
\bigskip

For the rotational group $\m D/\Z_3\cong \Z_4$ which is generated by $\beta$, we proceed
similarly. The action \eqref{eq:beta} of $\beta$ 
according to \eqref{dualbasis} induces 
$(\mu_1,\mu_2,\mu_3,\mu_4)\mapsto (-\mu_3,-\mu_4,\mu_1,\mu_2)$ and thereby
fixes $\mu_{13}$ and $\mu_{24}$, while $\mu_{12}\leftrightarrow\mu_{34}$
and $\mu_{14}\leftrightarrow\mu_{23}$. Proposition \ref{Kconstruct} thus yields
$\beta_\ast\colon (\kappa_1,\kappa_2,\kappa_3,\kappa_4) \longmapsto
(\kappa_1,\kappa_2,\kappa_4,\kappa_3)$. 

The action \eqref{eq:beta} of $\beta$  moreover permutes the fixed points $T^{\Z_3}$
of $\Z_3$ and thus induces an action $\widetilde\beta$ on the labels $t\in\F_3^2$,
\begin{align}
\forall (t_1,t_2)\in \F_3^2\colon\qquad
    \widetilde\beta: (t_1,t_2) \longmapsto (-t_2,t_1) \, .
\end{align}
If $\beta$ descends to a symmetry (also denoted $\beta$) of $X$, then 
by the Torelli Theorem \eqref{Torelli}, the lattice automorphism $\beta_\ast$ is $\Z_3$-effective,
such that for $t\in\F_3^2$,  from \eqref{blowupanddown} it follows that 
$\beta_\ast$  maps the set $\{ E_t^{(1)}, E_t^{(2)} \}$
to the set $\{ E_{\widetilde\beta(t)}^{(1)}, E_{\widetilde\beta(t)}^{(2)} \}$. 
The precise map can be read off from the local analysis of the blow-up according to
appendix \ref{app:A2blowup}. If $t=(0,0)$, then
we may describe the blow-up locally by the closures of the
three coordinate patches $U_0,\, U_1,\, U_2$ given in \eqref{U0}--\eqref{U2}
in $\C^3\times\P^2$. Since for the coordinates introduced there,
$(x_0,x_1,x_2)=(z_1z_2, z_1^3, z_2^3)$
and $u_jx_k=u_kx_j$ for all $j,k\in\{0,1,2\}$, we find an
induced action of $\beta$ on $\C^3\times\P^2$ given by
$$
\beta\colon \left( (x_0,x_1,x_2), (u_0:u_1:u_2) \right)
\longmapsto \left( (-x_0,-x_2,x_1), (-u_0:-u_2:u_1) \right)\, .
$$
Thus \eqref{exceptionalcomponents} shows that $\beta$ interchanges the two irreducible
components $\CCC_0^{(j)}$ of the exceptional divisor and therefore
$\beta_\ast (E_{00}^{(1)}) = E_{00}^{(2)}$, 
$\beta_\ast (E_{00}^{(2)}) = E_{00}^{(1)}$.
For arbitrary $t_0\in T^{\Z_3}$, let 
$\widehat\beta:=\alpha_{-\beta(t_0)}\circ\beta\circ\alpha_{t_0}$,
such that $\widehat\beta(0)=0$. 
Proceeding analogously to the discussion of $\beta$, we find 
that $\widehat\beta_\ast$ interchanges 
$E_{00}^{(1)}$ with $E_{00}^{(2)}$, and therefore, 
if $t\in\F_3^2$ is the label for the fixed point $t_0$, then 
$\beta_\ast (E_{t}^{(1)}) = E_{\widetilde\beta(t)}^{(2)}$, 
$\beta_\ast (E_{t}^{(2)}) = E_{\widetilde\beta(t)}^{(1)}$.
Hence we may conclude that 
$\beta$ does indeed descend to a well-defined symmetry of $X$, and 
$\beta_\ast (E_t^{(j)}) = E_{\widetilde\beta(t)}^{(3-j)}$
when $t\in\F_3^2$ and $j\in\{1,2\}$. 
Altogether we have shown:
\blem\label{rotationalaction}%
The action of the rotational symmetries of  $(X,\omega)$ on $H^2(X,\Z)$ is generated by
$\beta_\ast$, where 
$\beta_\ast (E_{(t_1,t_2)}^{(j)})=E_{(-t_2,t_1)}^{(3-j)}$ when
$(t_1,t_2)\in\F_3^2$ and $j\in\{1,2\}$, and $\beta_\ast(\kappa_j)=\kappa_j$
if  $j\in\{1,2\}$, $\beta_\ast(\kappa_3)=\kappa_4$ and $\beta_\ast(\kappa_4)=\kappa_3$.
\elem
To confirm that lemma \ref{rotationalaction} indeed yields a lattice automorphism
of $H^2(X,\Z)$, we also check compatibility with our description of the lattice $P$
in proposition \ref{Pconstruct} and with the gluing 
prescription of proposition \ref{gluePandK}.
Since $\beta_\ast( \frac{1}{3} E_t) = -\frac{1}{3} E_{\widetilde\beta(t)} + E_{\widetilde\beta(t)}^{(1)} 
+ E_{\widetilde\beta(t)}^{(2)}$ if $t\in\F_3^2$,
for the generators \eqref{eq:v1v2v3} of $P$ from $R$, we have
\begin{align*}
    \beta_\ast(v_1) \in - v_2 + R \, , \quad \beta_\ast(v_2)\in v_1+v_3 + R \, , \quad \beta_\ast(v_3)\in - v_3 +R\, .
\end{align*}
In other words, our system of generators is mapped to another such system. Similarly, for the
glue vectors $k_j+p_j$ found in proposition \ref{gluePandK},
\begin{align}\label{betaondisc}
    \beta_\ast(k_1+p_1) \in k_1 + p_1 + P \, , \quad
    \beta_\ast(k_2+p_2) \in -k_2 - p_2 + P \, , \quad
    \beta_\ast(k_3+p_3) \in k_3 + p_3 +P \, .
\end{align}
Hence all of these glue vectors are again mapped to glue vectors.
\bigskip

So far, we have constructed the group of those symmetries of  $(X,\omega)$ which are induced
by symmetries of the torus $T$. We claim that this already yields the entire symmetry
group of  $(X,\omega)$:

\bprop\label{allsymmetries}%
 The  symmetry group of $(X,\omega)$ is $(\Z_3)^2\rtimes \Z_4$, where $(\Z_3)^2$ is
 the group of  translational symmetries, induced by translations on the underlying
 torus $T$, and  $\Z_4$ is the group of rotational symmetries, induced by the 
 action of a subgroup of $\SU(2)$ on the underlying torus $T$.
\eprop

\bpf
According to lemmas \ref{translationalaction} and \ref{rotationalaction}, the group of  symmetries
of $(X,\omega)$ which are induced by symmetries of $T$ is $(\Z_3)^2\rtimes \Z_4$. Hence 
assuming that $f$ is any symmetry of $(X,\omega)$, we need to prove that $f$ is
induced by a symmetry of $T$. 
We shall show that by composing $f$  with appropriate symmetries  induced by symmetries
of $T$, one obtains a symmetry of $(X,\omega)$ whose induced $H^2(X,\Z)$ lattice automorphism is the identity.
From this the claim is immediate by the uniqueness statement in the Torelli Theorem \eqref{Torelli},
which implies that any symmetry of $(X,\omega)$ that acts trivially on  $H^2(X,\Z)$
is the identity and hence that $f$ is induced
by symmetries of $T$.

We  denote by 
$\Phi=f_\ast$ the lattice automorphism $\Phi\in\Aut(H^2(X,\Z))$ induced by $f$.
By \eqref{Torelli}, $\Phi$ is $\Z_3$-effective. Hence $\Phi$ acts as the identity on
the three-dimensional space $\Sigma$ with basis $(\kappa_1,\, \kappa_2,\, \kappa_3+\kappa_4 )$, and therefore
it restricts to a lattice automorphims of $\Sigma^\perp\cap H^2(X,\Z)$. In particular, 
$\Phi$ permutes the roots of that lattice, which by \eqref{rootsareinP} agree with the roots of $P$.
Since by construction, every element of $P$ is a (rational) linear combination of roots, this implies
that $\Phi$ maps $P$ to itself. Furthermore, $\Phi(\kappa_j)=\kappa_j$ for $j\in\{1,2\}$ and 
$\Phi(\kappa_3+\kappa_4)=\kappa_3+\kappa_4$ implies that either 
$\Phi(\kappa_j)=\kappa_j$ for $j\in\{3,4\}$ or $\Phi\colon\kappa_3\leftrightarrow\kappa_4$.
Since by
Lemma \ref{rotationalaction}, the generator $\beta$ of the rotational group $\m D/\Z_3\cong \Z_4$
interchanges $\kappa_3$ and $\kappa_4$,
by replacing $f$ by $f\circ\beta^{-1}$
if necessary we may assume without 
loss of generality that $\Phi$ acts trivially on 
$K$.

Since $\Phi$ acts as the identity on $K$, it also acts as the identity on the 
discriminant group of $K$ and thereby, according to \eqref{autoextend}, it acts as
the identity on the discriminant group of $P$. Let $\psi:=\Phi_{|P}$ and extend
$\psi$ linearly to $\psi\in\Aut( P^\ast)$, where according to proposition \ref{Pconstruct},
the lattice $P^\ast$ is generated by the $E_t^{(j)}$ with $t\in\F_3^2$, $j\in\{1,2\}$,
along with the vectors $\frac{1}{3}\sum_{t\in L^\prime}E_t$ with $L^\prime\subset\F_3^2$ any affine line.
All these generators have length square $-2$. In fact, proposition \ref{Pconstruct}
implies that for every $v\in P^\ast\setminus P$
with $<v,v>=-2$, there exists an affine line $L^\prime\subset\F_3^2$ such that $v$
is one of the following:
\begin{equation}\label{rootsinPstar}
\textstyle
\pm\left(\frac{1}{3}\sum\limits_{t\in L^\prime}E_t-\sum\limits_{t\in L^\prime}\delta_t(E_t^{(1)}+E_t^{(2)})\,\right)\,,\quad
\delta_t\in\{0,1\} \mbox{ for } t\in L^\prime.
\end{equation}
Since $\Phi$ is $\Z_3$-effective, 
$\psi$ permutes the roots $E_t^{(j)}$ with $t\in\F_3^2$, $j\in\{1,2\}$.
Hence $\psi\in\Aut(P^\ast)$ also permutes the vectors listed in \eqref{rootsinPstar}.
More precisely, since $\psi$ induces the identity on the discriminant of $P$
and since it permutes the  $E_t^{(j)}$,
if $L^\prime\subset \F_3^2$ is an affine line,
then there exists a parallel $L^{\prime\prime}\subset \F_3^2$ 
 of $L^\prime$ such that 
$\frac{1}{3}\sum_{t\in L^\prime}E_t$ is mapped
to  $\frac{1}{3}\sum_{t\in L^{\prime\prime}}E_t$ by $\psi$.
It follows that there exists a permutation $\widetilde\psi$ of the points
of $\F_3^2$ such that $\psi(E_t) = E_{\widetilde\psi(t)}$ for every $t\in\F_3^2$.
In fact, since $\psi$ permutes the $E_t^{(j)}$, we have 
$\psi(E_t^{(j)}) = E_{\widetilde\psi(t)}^{(j)}$ for every $t\in\F_3^2$
and $j\in\{1,2\}$. Moreover, by the above $\psi$ permutes the elements of $P^\ast$
of the form $\frac{1}{3}\sum_{t\in L^\prime}E_t$ with $L^\prime\subset\F_3^2$ any affine line. Hence
$\widetilde\psi$ maps affine lines to affine lines and therefore is affine linear.

By composing the symmetry $f$ with an appropriate translational
symmetry, we may assume without loss of generality that $\widetilde\psi$ fixes $(0,0)$,
that is, $\psi(E_t^{(j)}) = E_{\widetilde\psi(t)}^{(j)}$ 
for every $t\in\F_3^2$ and $j\in\{1,2\}$, where $\widetilde\psi$ is a  linear map
on $\F_3^2$. 
Again using the fact that $\psi$ acts as the identity on
the discriminant group of $P$, it follows that $\widetilde\psi$ maps each of the lines
$L_{12},\, L_{34},\, L_{00}$ to itself, fixing the point $(0,0)$. In particular, if
$\widetilde\psi$ does not fix $L_{12}$ pointwise, then it interchanges the points
$(0,1)$ and $(0,2)$, as does $\widetilde\beta^2$ (with notations as in the discussion 
preceding lemma \ref{rotationalaction}). Hence by composing $f$ with the symmetry
that is induced by $\beta^2$ if necessary,
we may assume without loss of generality that $\widetilde\psi$ fixes the line $L_{12}$
pointwise.

Let us now show that our assumptions imply that $\widetilde\psi$ is the identity on 
$\F_3^2$. Indeed, $\widetilde\psi$ acts linearly on $\F_3^2$,
and it fixes the line $L_{12}$ pointwise. It also maps the lines $L_{34}$
and $L_{00}$ to themselves. Each of these lines
has two distinct parallels, which are either mapped to themselves or interchanged by
$\widetilde\psi$, since $\psi$ induces the identity on the discriminant
group of $P$. Since $(0,1)$ and $(0,2)$ are fixed points under $\widetilde\psi$
belonging to the two 
distinct parallels of $L_{34}$, these two lines are also mapped to themselves.
In particular, $(1,1)$ is not interchanged with $(2,2)$, which implies that the line
$L_{00}$ is also fixed pointwise. Proceeding in this manner, one checks that 
$\widetilde\psi$ is the identity.

Altogether, by composing $f$ with some symmetry
which is induced by a symmetry of $T$ we have shown that without loss of generality $\Phi$
is the identity automorphism, as desired.
\epf
\brmk\label{generalsymmetries}%
For ease of exposition, in proposition \ref{allsymmetries} and in lemmas \ref{translationalaction}, 
\ref{rotationalaction} above we restricted
our discussion to the symmetry group of $(X,\omega)$, where $\omega$ denotes the standard K\"ahler class on $X$.
According to proposition \ref{hyperk}, our notion of $\Z_3$-orbifold limits of K3 allows the choice of a more general
K\"ahler class $\wt\omega=V_3\kappa_3+V_4\kappa_4$ with $V_3,\, V_4\in\R$ and $V_3V_4>0$. We claim that the steps detailed
above in determining the symmetry group of $(X,\omega)$ apply analogously to the more general case $(X,\wt\omega)$. 
This requires to replace $\Sigma$ as introduced in \eqref{Sigmadef} by 
$$
\wt\Sigma:=\mbox{span}_\R \left\{ \kappa_1,\, \kappa_2, \, \wt\omega\right\}\, .
$$
The entire derivation yielding lemmas \ref{translationalaction}, \ref{rotationalaction} and proposition \ref{allsymmetries}
is based on the detailed knowledge of the lattice $\Sigma^\perp\cap H^2(X,\Z)$, which is now replaced by the lattice
$\wt\Sigma^\perp\cap H^2(X,\Z)$. To understand this lattice, the observation \eqref{rootsareinP} is crucial, which
continues to hold when $\Sigma$ is replaced by $\wt\Sigma$, as we shall show now.\\
If no $V\in\R\setminus\{0\}$ exists with $V\wt\omega\in H^2(X,\Z)$, then $\wt\Sigma^\perp\cap H^2(X,\Z)=P$ and 
\eqref{rootsareinP} holds trivially when $\Sigma$ is replaced by $\wt\Sigma$. If on the other hand $V\in\R\setminus\{0\}$
exists such that $V\wt\omega\in H^2(X,\Z)$, then $\wt\omega=\wt V(M_3\kappa_3+M_4\kappa_4)$ with $\wt V\in\R$ and 
$M_3,\, M_4\in\N$, and without loss of generality we may assume that $M_3$ and $M_4$ are coprime. Then the glue vector
obtained in \eqref{Kaehlerglue} is replaced by the glue vector
$$
\textstyle
\frac{1}{3} (M_4\kappa_3-M_3\kappa_4) + \frac{M_4}{3} \sum\limits_{t\in L_{12}} E_t + \frac{M_3}{3} \sum\limits_{t\in L_{34}} E_t\; .
$$
If $M_3=M_4=1$, then this is just the expression \eqref{Kaehlerglue} and $\Sigma=\wt\Sigma$, such that nothing remains to be shown.
Otherwise, $M_3M_4\geq2$, and $\wt v:=\frac{1}{3} (M_4\kappa_3-M_3\kappa_4)$ obeys 
$\left|\langle\wt v,\wt v\rangle\right|=\frac{2}{3} M_3M_4$. That \eqref{rootsareinP} 
continues to hold when $\Sigma$ is replaced by $\wt\Sigma$ is therefore immediate:
assume that $v\in \wt\Sigma^\perp\cap H^2(X,\Z)$ with $v\not\in P$. 
Then $v=M\wt v+e$ for some $M\in\Z\setminus\{0\}$ and $e\in P^\ast$.
If $M^2M_3M_4>3$, then $\left|\langle v, v\rangle\right|\geq M^2 \left|\langle\wt v,\wt v\rangle\right|>2$, hence
$\langle v, v\rangle\neq -2$. Otherwise, $M^2M_3M_4\in\{2,3\}$, which implies $M^2=1$ and $M_3=1$ or $M_4=1$. It follows that 
$e\neq0$, in fact with $\wt e=\frac{1}{3} \sum\limits_{t\in L_{12}} E_t$ or  $\wt e=\frac{1}{3} \sum\limits_{t\in L_{34}} E_t$,
we have 
$\left|\langle v, v\rangle\right|\geq \left|\langle\wt v,\wt v\rangle\right| + \left|\langle\wt e,\wt e\rangle\right| 
> \left|\langle\wt e,\wt e\rangle\right| =2$, hence
$\langle v, v\rangle\neq -2$. We have established that indeed \eqref{rootsareinP} 
continues to hold when $\Sigma$ is replaced by $\wt\Sigma$.\\
To proceed, we observe that the metric on $T$ which induces the general K\"ahler class 
$\wt\omega=V_3\kappa_3+V_4\kappa_4$ with $V_3,\, V_4\in\R$ and $V_3V_4>0$ on $X$ breaks the 
symmetry group $T\rtimes\mathcal D$ obtained with respect to the standard Euclidean metric on $T$
if and only if $V_3\neq V_4$. Namely, in this case one checks that of the symmetries listed in \eqref{eq:beta},
the symmetry $\beta$  is broken (indeed $\beta$ interchanges the two factors of $T$, which are given by 
elliptic curves of different volumes if $V_3\neq V_4$; one also checks
$\beta_\ast\wt\omega\neq\wt\omega$ since $\beta_\ast$ interchanges $\kappa_3$ and $\kappa_4$), 
while $\varrho$ remains a symmetry. This means that the symmetry group of $T$ is broken to $T\rtimes\wt{\mathcal D}$,
where $\wt{\mathcal D}\cong\Z_6$ is the cyclic group of order $6$ which is generated by $\varrho$. In particular,
$\varrho^3=\beta^2$, which acts by multiplication by $-1$ on $\C^2$, belongs to $\wt{\mathcal D}$. 
This induces an action of $(\Z_3)^2\rtimes\Z_2$ on $T/\Z_3$, as follows by the same arguments that we used to
find the $(\Z_3)^2\rtimes\Z_4$-action in the standard case. One now checks that the remaining steps leading to 
lemmas \ref{translationalaction}, \ref{rotationalaction} and proposition \ref{allsymmetries} continue to 
hold when $\Sigma$ is replaced  by $\wt\Sigma$. Altogether we obtain the following.\\
\begin{framed}\vspace*{-0.3em}
The symmetry group of $(X,\wt\omega)$ with $\wt\omega=V_3\kappa_3+V_4\kappa_4$, $V_3,\, V_4\in\R$ and $V_3V_4>0$
agrees with the symmetry group $(\Z_3)^2\rtimes\Z_4$ of $(X,\omega)$ described in proposition \ref{allsymmetries}
if $V_3=V_4$. If $V_3\neq V_4$, then the symmetry group is $(\Z_3)^2\rtimes\Z_2$, where $(\Z_3)^2$ is the group 
of translational symmetries described in lemma \ref{translationalaction}, while $\Z_2$ is the subgroup 
of the rotational symmetry group described in lemma \ref{rotationalaction} which is generated by $\beta_\ast^2$.
\end{framed}
\ermk

Inspection of the actions described in lemmas \ref{translationalaction} and
\ref{rotationalaction} also shows:

\bcor\label{F32isenough}%
The symmetries of  $(X,\wt\omega)$ with $\wt\omega$ as in proposition \ref{hyperk}
are uniquely determined by their induced actions on $\F_3^2$
as affine linear maps.
\ecor
\section{Representation of the symmetries of \texorpdfstring{$\Z_3$}{TEXT}-orbifold K3s on a Niemeier lattice}
\label{sec:embedding}

So far, we have described the integral cohomology lattice $H^2(X,\Z)$ of our K3
surface $X=\widetilde{T/\Z_3}$ in terms of the  root lattice $R$ obtained in the 
Kummer-like construction of  $X$. Here, $R$ is of type $A_2^9$ and is 
generated by the components of the exceptional divisor that arises from the resolutions 
in the Kummer-like construction of $X$.
We made use of the rank-$18$, negative definite, 
Kummer-like lattice $P\supset R$ obtained in proposition \ref{Pconstruct}, and its orthogonal 
complement $K$ obtained in proposition \ref{Kconstruct}. There, we
furthermore used Nikulin's gluing techniques to describe $H^2(X,\Z)$
in terms of the primitive sublattices $K$ and $P$. In proposition \ref{allsymmetries},
this allowed us to determine all symmetries of $X$,  and to conveniently
describe them in terms of automorphisms of $H^2(X,\Z)$ in lemmas \ref{translationalaction} and 
\ref{rotationalaction}. 

In order to track such symmetries, it is natural to view them as permutation groups,
acting, for example, on each set of vectors of equal length in the lattice $H^2(X,\Z)$. However,
since this lattice has signature $(3,19)$,  every non-zero lattice vector 
is a member of an infinite  set of lattice vectors of the same length. 
It is therefore more convenient
to represent our symmetries on a lattice of definite signature. To do so, inspired
by Kondo's seminal paper \cite{ko98a}, and since corollary \ref{F32isenough} implies
that each of our symmetries is uniquely determined by its action on $P$, 
in subsection \ref{subsec:embeddingP} we search for primitive embeddings of the lattice $P(-1)$ in 
some Niemeier lattice.  
It turns out that the unique possibility is the Niemeier lattice $N$ of type $A_2^{12}$
(proposition \ref{primitiveNiemeierisunique}), and that the embedding is unique up to
lattice automorphisms (proposition \ref{primitiveembeddingisunique}).
The lattice $N$ is introduced in appendix \ref{subapp:niemeier} and may be generated by the 
set \eqref{Ngen}. An explicit primitive embedding of the Kummer-like 
lattice $P(-1)$ into $N$ is given in proposition \ref{embeddingexists}.
In remarks \ref{E64} and \ref{symmetrybreaking} we also discuss the only existing 
\emph{non-primitive}
embeddings of $P(-1)$ into a Niemeier lattice of different type, which can only be
of type $E_6^4$.

In subsection \ref{subsec:symmetriesinM12}, we argue that our embedding of $P(-1)$
into the Niemeier lattice $N$ allows us to represent the symmetry group of $X$
on $N$. This is a variation of the methods introduced by Kondo in \cite{ko98a}, 
and we discuss the advantages and disadvantages of this variation; in fact, as we show in the
introduction to subsection \ref{subsec:symmetriesinM12}, our results prove that $N$ is
one of the lattices that may be used to represent the symmetries of $X$
on a Niemeier lattice by Kondo's original method. That this should be the case 
was already claimed in 
\cite[p.29, bottom]{chha15}, though without a proof. To our knowledge,
our above-mentioned uniqueness result for the Niemeier lattice $N$ to allow a primitive
embedding of $P(-1)$ is also new. We arrive at an explicit representation of
the symmetry group of $(X,\omega)$ on the lattice $N$ in propositions \ref{translationsonN} 
and \ref{rotationsonN}. 
Since the symmetry group of $(X,\wt\omega)$, with the general choice of $\wt\omega$ obtained in proposition \ref{hyperk},
is a subgroup of the symmetry group of $(X,\omega)$ by remark \ref{generalsymmetries}, 
this representation of the symmetry group of $(X,\omega)$ is sufficient to represent
all symmetry groups of $\Z_3$-orbifold K3s. We may therefore restrict our attention to the symmetry group of
$(X,\omega)$ henceforth, which for simplicity we refer to as the symmetry group of $\Z_3$-orbifold K3 surfaces $X$.
We argue that the projection of $\Aut(N)$ to the Mathieu
group $M_{12}$ described in appendix \ref{subapp:niemeier} yields an injective
image of the symmetry group within $M_{12}$, explicitly given in corollary 
\ref{subgroupofM12}.

That the Mathieu group
$M_{12}$ is a subgroup of the Mathieu group $M_{24}$ allows us to realize the symmetry
group of $\Z_3$-orbifold K3 surfaces as subgroup of the largest Mathieu group
in subsection \ref{subsec:symmetriesinM24}. This implies that one may represent our 
symmetries on the Niemeier lattice of type $A_1^{24}$; for some background 
information on this lattice and its relation to the largest Mathieu group
we refer to appendix \ref{subapp:m24}.
We explicitly express our symmetry group as subgroup of $M_{24}$ in proposition 
\ref{imageinM24}, and in theorem \ref{generatingM24} we show that this group
together with the overarching symmetry group of all Kummer surfaces
generates the entire group $M_{24}$.
\subsection{Embeddings of the Kummer-like lattice in Niemeier lattices}\label{subsec:embeddingP}

In the following, the Niemeier lattice $N$ of type $A_2^{12}$  plays a fundamental role. 
It is an even, self-dual lattice of signature $(24,0)$ which contains a root 
sublattice $\wt R$ of type $A_2^{12}$;
$\wt R$ has index $3^6$ in $N$. More details can be found 
in appendix \ref{subapp:niemeier}, where we also introduce our notations.
In particular, our choice of simple roots in $\wt R$ is denoted $\wt E_j^{(k)}$ with 
$j\in\{1,\ldots,12\}$ labelling the irreducible $A_2$-type sublattices of $\wt R$
and $k\in\{1,2\}$, and \eqref{gluecodewords} gives a minimal choice of 
generators $(w_1,\ldots,w_6)$ of $N$ from $\wt R$, using the glue matrix \eqref{matrixcode}.
Let us show that the Kummer-like lattice $P$ which was constructed in proposition 
\ref{Pconstruct} admits a primitive embedding $P(-1)\hookrightarrow N$:

\bprop\label{embeddingexists}%
The $\Q$-linear extension $\iota$ of the map 
$$
\mbox{for } \ell\in\{1,2\},\quad
\iota\colon
\left\{
\begin{array}{rcrrcrrcr}
E_{00}^{(\ell)}&\longmapsto&\wt E_6^{(\ell)}\, ,&
E_{01}^{(\ell)}&\longmapsto&\wt E_4^{(\ell)}\, ,&
E_{02}^{(\ell)}&\longmapsto&\wt E_7^{(\ell)}\, ,\\[0.3em]
E_{10}^{(\ell)}&\longmapsto&-\wt E_5^{(\ell)}\, ,&
E_{11}^{(\ell)}&\longmapsto&-\wt E_3^{(\ell)}\, ,&
E_{12}^{(\ell)}&\longmapsto&-\wt E_2^{(\ell)}\, ,\\[0.3em]
E_{20}^{(\ell)}&\longmapsto&\wt E_8^{(\ell)}\, ,&
E_{21}^{(\ell)}&\longmapsto&\wt E_1^{(\ell)}\, ,&
E_{22}^{(\ell)}&\longmapsto&\wt E_9^{(\ell)}\, ,
\end{array}
\right.
$$
yields a primitive embedding  $\iota\colon P(-1)\hookrightarrow N$ with 
image 
\[
\widetilde{P}
:={\rm span}_\Z\left\{\widetilde{E}^{(\ell)}_j,\; j\in \{1,\ldots,9\},\, \ell\in \{1,2\};\quad 
w_1,\, w_2,\, w_3\right\}\, .
\]
\eprop

\bpf
Recall from equation \eqref{eq:v1v2v3} that the lattice $P$ is 
obtained from its root sublattice $R$ by including the three additional generators
$v_1,\, v_2,\, v_3$. By direct calculation one confirms that
\[
\iota(v_1)= w_1,\quad \iota(v_2)= w_2,\quad 
\iota(v_3)=
w_1+w_2+w_3+\widetilde{E}_1-\widetilde{E}_2-\widetilde{E}_3\,,
\]
which implies that $\iota$ gives an isometry from $P(-1)$ to $\widetilde P$.
Using \eqref{matrixcode} we moreover identify
$\widetilde P$ as the orthogonal complement of the sublattice 
$\bigoplus\limits_{j=10}^{12} A_{2,j}$ in $N$; this implies that the embedding 
$\iota\colon P(-1)\hookrightarrow N$ is primitive.
\epf

Since $\wt P$ is a primitive sublattice of $N$, we may use Nikulin's gluing
techniques as discussed in appendix \ref{subapp:glue} to describe $N$ in terms of $\wt P$
and its orthogonal complement; for later convenience, let us read off the 
corresponding gluing isomorphism explicitly:

\blem\label{gluinginN}%
For $\wt P=\iota(P(-1))$ as in proposition \ref{embeddingexists}, let 
$\widetilde K:=\widetilde P^\perp\cap N$. Using the $\Q$-linear extension of $\iota$
and the generators $p_1,\, p_2,\, p_3$ of $P^\ast$ from $P$ that were determined in 
proposition \ref{Pconstruct},
also let $\widetilde p_j:=\iota(p_j)$, $j\in\{1,2,3\}$. 
Then $\wt K=\bigoplus\limits_{j=10}^{12} A_{2,j}$, and
the gluing isomorphism $\widetilde\gamma:\widetilde P^\ast/\wt P\longrightarrow\wt K^\ast/\wt K$
can be given by 
$$
\textstyle
\wt\gamma( \wt p_1+\wt P) = \frac{1}{3}\wt E_{10}+\wt K,\quad
\wt\gamma( \wt p_2+\wt P) = -\frac{1}{3}\wt E_{12}+\wt K,\quad
\wt\gamma( \wt p_3+\wt P) = -\frac{1}{3}\wt E_{11}+\wt K.
$$
\elem

\bpf
By the explicit data given in proposition \ref{embeddingexists}, it is immediate that
$\wt K$ contains the root lattice $\bigoplus\limits_{j=10}^{12} A_{2,j}$.
To prove equality, consider  
$\wt k=\sum\limits_{\ell=1}^6 c_\ell w_\ell+\wt r$ with $c_1,\ldots,c_6\in\Z$
and $\wt r\in\wt R$ such that 
$\wt k\in\wt K$; we must 
show that all coefficients $c_\ell$ are multiples of $3$. 
With $w_\ell = \frac{1}{3} \sum\limits_{j=1}^{12} a_{\ell j} \widetilde E_j$ as in
\eqref{gluecodewords}, observe that 
$\wt k\perp A_{2,j}$ for all $j\in\{1,\ldots,9\}$ implies $\sum\limits_{\ell=1}^6 c_\ell a_{\ell j}\in 3\Z$. 
Modulo $3\Z$ the coefficients $a_{\ell j}\in\{0,\pm1\}$ yield elements $\ol a_{\ell j}$ of the finite field $\F_3$. 
Now the claim is immediate from 
inspection of the glue matrix \eqref{matrixcode},  which shows that the $6\times6$ block
given by the entries $\ol a_{ij}\in\F_3$ with $i,j\in\{1,\ldots,6\}$ forms an invertible
$6\times6$-matrix over $\F_3$.

To construct the gluing isomorphism $\wt\gamma:\wt P^\ast/\wt P\longrightarrow\wt K^\ast/\wt K$,
we need to find $\wt k_1,\, \wt k_2,\, \wt k_3\in\wt K^\ast$ such that $\wt p_j+\wt k_j\in N$
for $j\in\{1,2,3\}$; then $\wt\gamma(\wt p_j+\wt P)=\wt k_j+\wt K$
for $j\in\{1,2,3\}$ uniquely determines $\wt\gamma$.
To this end one checks by direct calculation that
$\wt p_1+\frac{1}{3} \wt E_{10}\in w_2+w_3+w_4+\wt R$, which implies that 
we may choose $\wt k_1=\frac{1}{3} \wt E_{10}$; similarly,
$\wt p_2-\frac{1}{3} \wt E_{12}\in w_1+w_2-w_6+\wt R$
and $\wt p_3-\frac{1}{3} \wt E_{11}\in w_1+w_3-w_5+\wt R$
proves the remaining claims.
\epf

Since in the above construction, we have made many choices, it is natural to ask about uniqueness.
Fixing $N$, we first obtain

\bprop\label{primitiveembeddingisunique}%
The primitive embedding of $P(-1)$ in the Niemeier lattice $N$ is unique
up to automorphisms of $N$.
\eprop

\bpf
Assume that 
$\wh P\subset N$ is the image of $P(-1)$ under another primitive embedding 
$\wh\iota\colon P(-1)\hookrightarrow N$, such that 
$\varphi:=\wh\iota\circ\iota^{-1}\colon \wt P\longrightarrow \wh P$ is a lattice isomorphism.
Our claim is equivalent to the existence of some
lattice automorphism $\Phi\in\Aut(N)$ which on $\wt P$ restricts to $\varphi$.

To construct such $\Phi$, first set $\wh K:=\wh P^\perp\cap N$. By the same 
arguments as in the proof of lemma \ref{gluinginN}, we may conclude that $\wh K$
is a root lattice of type $A_2^3$. In fact, again using the generators $p_1,\, p_2,\, p_3$
of $P^\ast$ from $P$ that we found in proposition \ref{Pconstruct}, let
$\wh p_j:=\wh\iota(p_j)$ for $j\in\{1,2,3\}$, 
$\wh E_t:=\wh\iota( E_t )$, 
$\wh E_t^{(\ell)}:=\wh\iota\left( E_t^{\smash{(\ell)}} \right)$ if $t\in\F_3^2$, $\ell\in\{1,2\}$.
Then according to proposition \ref{Pconstruct}, with respect to the generators
$\left(\wh p_1+\wh P,\wh p_2+\wh P,\wh p_3+\wh P\right)$ of $\wh P^\ast/\wh P$ the 
induced bilinear form is diagonal and the discriminant form $q_{\wh P}$ takes value 
$-\frac{2}{3}\mod2\Z$. Hence
the gluing isomorphism $\wh\gamma\colon \wh P^\ast/\wh P\longrightarrow\wh K^\ast/\wh K$
for $N$ from $\wh P$ and $\wh K$ maps $\left\{\wh p_1+\wh P,\wh p_2+\wh P,\wh p_3+\wh P\right\}$
to a generating set of $\wh K^\ast/\wh K$ on which the bilinear form $b_{\wh K}$ is diagonal 
and the discriminant form $q_{\wh K}$ takes value $\frac{2}{3}\mod2\Z$. By explicit calculation
of the discriminant form on $\wh K^\ast/\wh K$ this implies that 
for each $j\in\{1,2,3\}$ we have
$\wh\gamma(\wh p_j+\wh P)=\frac{1}{3} \wh E_{n_j}+\wh K$ or 
$\wh\gamma(\wh p_j+\wh P)=-\frac{1}{3} \wh E_{n_j}+\wh K$, 
where $n_1, n_2, n_3\in\{1,\ldots,12\}$ 
with $\wh K=\bigoplus\limits_{j=1}^3 A_{2,n_j}$.

Now let $\ol\varphi\colon \wt P^\ast/\wt P\longrightarrow \wh P^\ast/\wh P$
denote the isomorphism which is induced by $\varphi$. Note that by construction,
$\ol\varphi(\wt p_j+\wt P)=\wh p_j+\wh P$ if $j\in\{1,2,3\}$.
Let 
$\ol\psi\colon \wt K^\ast/\wt K\longrightarrow \wh K^\ast/\wh K$
such that $\wh\gamma\circ\ol\varphi=\ol\psi\circ\wt\gamma$, where 
$\wt\gamma$ is the gluing isomorphism that we determined in  lemma \ref{gluinginN}.
To ease notation below, we set $(m_1,m_2,m_3):=(10,12,11)$ and observe
that by the above, $\ol\psi(\frac{1}{3}\wt E_{m_j}+\wt K)=\pm\frac{1}{3}\wh E_{n_j}+\wh K$
for $j\in\{1,2,3\}$.
Since $\wt K\cong\wh K$ are root lattices of type $A_2^3$, we may 
explicitly construct a lattice isomorphism $\psi$ between $\wt K$ and $\wh K$ which induces $\ol\psi$.
Indeed, we may take $\psi(\wt E_{m_j}^{(\ell)}):=\wh E_{n_j}^{(\ell)}$ if 
$\ol\psi(\frac{1}{3}\wt E_{m_j}+\wt K)=\frac{1}{3}\wh E_{n_j}+\wh K$ and 
$\psi(\wt E_{m_j}^{(\ell)}):=\wh E_{n_j}^{(3-\ell)}$ otherwise. 

By $\Q$-linearly extending $\varphi$ and $\psi$ to $\wt P^\ast$ and $\wt K^\ast$
(denoting the resulting maps by $\varphi$, $\psi$, as well), we obtain a lattice isomorphism
$\Phi\colon \wt P^\ast\oplus\wt K^\ast\longrightarrow \wh P^\ast\oplus\wh K^\ast$
with $\Phi(\wt p+\wt k):=\varphi(\wt p)+\psi(\wt k)$ if $\wt p\in\wt P^\ast$ and 
$\wt k\in\wt K^\ast$. That $\Phi$ 
restricts to a lattice automorphism of $N$ with the desired properties 
now follows from a straightforward generalization of \eqref{autoextend}.

In detail, Nikulin's gluing results as discussed in appendix \ref{subapp:glue} imply
\begin{eqnarray*}
N
&=& \left\{\wt p+\wt k\in \wt P^\ast\oplus\wt K^\ast \mid \wt\gamma(\wt p+\wt P)=\wt k+\wt K\right\}\\
&=& \left\{\wh p+\wh k\in \wh P^\ast\oplus\wh K^\ast \mid \wh\gamma(\wh p+\wh P)=\wh k+\wh K\right\}\, .   
\end{eqnarray*}
Because  $\wh\gamma\circ\ol\varphi=\ol\psi\circ\wt\gamma$
by construction, for $\wt p+\wt k\in \wt P^\ast\oplus\wt K^\ast$
we have $\wt p+\wt k\in N$ if and only if $\varphi(\wt p)+\psi(\wt k)\in N$; in other words,
$\Phi$ restricts to a lattice automorphism of $N$, and $\Phi_{|\wt P}=\varphi$ by construction.
This concludes the proof.
\epf

Next, let us show  that  the lattice $N$,
up to lattice automorphisms, is also uniquely determined by the requirement that it
is even, self-dual of signature $(24,0)$, and that it admits a
primitive embedding of $P(-1)$:

\bprop\label{primitiveNiemeierisunique}%
Assume that $\widehat N$ is a Niemeier lattice which allows a primitive embedding of
$P(-1)$. Then $\widehat N$ is a Niemeier lattice of type $A_2^{12}$.
\eprop

\bpf
Denote by $\widehat E_t$, respectively $\widehat E_t^{(j)}$, with $t\in\F_3^2$
and $j\in\{1,2\}$ the images of the $E_t$ and the $E_t^{(j)}$ under 
the embedding of $P(-1)$ into $\widehat N$, and by $\widehat P$ the image of $P(-1)$.
Now assume that $\widehat r\in\widehat N$ is a root in $\widehat N$.
Then $\widehat r$ decomposes as
$\widehat r =\widehat p+\widehat k$
with $\widehat p\in \widehat P\otimes_\Z\Q$ and $\widehat k\perp\widehat P$.
Since $\widehat N$ is an integral lattice, we have $\widehat p\in \widehat P^\ast$. 
If $\wh p\neq0$, then according to lemma
\ref{lemrootsinPstar} it follows that $\langle\widehat p,\wh p\rangle\geq2$, such that 
our assumption $\langle\widehat r,\wh r\rangle=2$ implies $\wh k= 0$. In other words, 
$\widehat r=\widehat k$ or  $\widehat r=\widehat p$.
Since by our assumptions, $\widehat P$ is a primitive sublattice of $\widehat N$,
if $\widehat r=\widehat p$, then $\widehat p\in\widehat P$ follows.

Altogether, we have 
shown that the root lattice of $\widehat N$ decomposes as an orthogonal direct
sum isomorphic to $\bigoplus\limits_{j=1}^9 A_{2,j}\bigoplus \widehat K$, where $\widehat K$
is a root lattice of rank $6$. Comparing to the list of root lattices for Niemeier lattices 
(see, for example, \cite[Table 16.1]{conway1998sphere}), 
we may conclude that $\widehat K$ is of type $A_2^{3}$ and thus 
$\widehat N$ is of type $A_2^{12}$.
\epf
\brmk\label{E64}%
While proposition \ref{primitiveNiemeierisunique} completes the discussion of 
\emph{primitive} embeddings of $P(-1)$ into Niemeier lattices, we shed some light 
on \emph{non-primitive} embeddings into Niemeier lattices in appendix \ref{app:npembedding}.
There, we prove that up to isometry, the unique Niemeier lattice that allows
a non-primitive embedding of $P(-1)$ is of type $E_6^4$.
The fact that $P(-1)$ can be embedded in the Niemeier lattice $\wh N$ of type 
$E_6^4$ does not come as a surprise, since our K3 surface $X$ is elliptically
fibered over $\C\P^1$ with three singular fibers of Kodaira type $IV^\ast$
(see for example \cite[\S4]{nawe01}). Every singular fiber of type $IV^\ast$
decomposes into seven irreducible components, each given by a $(-2)$-curve, and 
arranged in a pattern according to the extended Dynkin diagram of $E_6$.
In remark \ref{symmetrybreaking} below, we shall see that the non-primitive embedding 
of $P(-1)$ into the Niemeier lattice $\wh N$ of type $E_6^4$
does not allow to express our K3 symmetry $\beta$ in terms of a lattice automorphism
of $\wh N$, while the translational symmetries as well as $\beta^2$ do induce well-defined automorphisms
of $\wh N$. This is in accord with the fact that the translational symmetries and $\beta^2$ respect
this elliptic fibration of $X$, but $\beta$ does not.
\ermk


\subsection{Tracking the symmetries of \texorpdfstring{$\Z_3$}{TEXT}-orbifold K3s within the 
Mathieu group \texorpdfstring{$M_{12}$}{TEXT}}\label{subsec:symmetriesinM12}

As mentioned at the beginning of this section,
in order to track the symmetries of the $\Z_3$-orbifold K3 $X$, it is convenient to represent 
them on a Niemeier lattice, since this allows us to view them as permutation groups 
of finitely many lattice vectors. This is the idea underlying Kondo's seminal paper \cite{ko98a}, 
where it is shown that every symmetry group $G$ of a K3 surface $Y$ can be 
faithfully represented on some Niemeier lattice. In more detail, using the induced action
of $G$ on the K3 lattice $H^2(Y,\Z)$,  let $M_G$ denote the orthogonal complement
of $H^2(Y,\Z)^G$ in $H^2(Y,\Z)$. In \cite[Lemma 5]{ko98a}, Kondo shows that
there exists a Niemeier lattice which allows a primitive embedding of $M_G(-1)$
and such that the action of $G$ can be extended to a lattice automorphism of 
the Niemeier lattice.
\bigskip

In the following, we will be using a variation of this construction. 
We now restrict our attention to our $\Z_3$-orbifold K3 $X$, and
the symmetry group $G$ of  $(X,\omega)$, i.e.~$G\cong(\Z_3)^2\rtimes\Z_4$
by proposition \ref{allsymmetries}. Since $H^2(X,\Z)^G$ contains the lattice
that generates $\Sigma$ (see lemma \ref{H2Tinv} and \eqref{Sigmadef}), it follows that 
$M_G\subset \Sigma^\perp\cap H^2(X,\Z)$. To arrive at \eqref{rootsareinP}, we have seen 
that $\Sigma^\perp\cap H^2(X,\Z)$ is
the lattice generated by $P$ and $k_2+p_2$,
where we continue to
use the notations introduced in propositions \ref{Pconstruct} and \ref{Kconstruct}.
Motivated by corollary \ref{F32isenough}, which also implies that the action of $G$ is uniquely 
determined by its action on $P$, we find it more convenient to focus on the lattice $P$ rather
than the lattice $M_G$.
According to propositions \ref{primitiveembeddingisunique} and \ref{primitiveNiemeierisunique}, 
for the lattice $P$,
the Niemeier lattice $N$ of type $A_2^{12}$ is unique in allowing a primitive embedding
$\iota\colon P(-1)\hookrightarrow N$, and  this embedding is unique up to lattice 
automorphisms. One checks that $N$ also allows a primitive embedding of 
$(\Sigma^\perp\cap H^2(X,\Z))(-1)$ which on $P(-1)$ restricts to $\iota$, 
implying that $M_G(-1)$ embeds primitively in $N$.
The induced action of $G$
on the discriminant of $M_G$ is trivial by construction. 
According to Nikulin's gluing technique \eqref{autoextend},  this means that the 
action on $M_G(-1)$ can be extended  to $N$ by a trivial action on the orthogonal 
complement of $M_G(-1)$.

Our focus on the lattice $P$ instead of $M_G$ is
a variation of Kondo's arguments which two of the authors have already successfully applied 
similarly in the Kummer case \cite{tawe11,tawe12,tawe13}.  
We find it convenient, because it achieves the same purpose of 
tracking symmetries of $X$ as lattice automorphisms on a Niemeier lattice, but it
allows us to base our calculations on root
systems, while $M_G$ according to \cite[Theorem 4.3]{ni80b} cannot contain any roots.
Two of the authors have also previously mentioned  that one can 
naturally avoid the confusing transition  
from rank $22$ for the lattice $H^2(X,\Z)$ to rank $24$ for the lattice $N$ 
by representing $G$ on the total integral cohomology 
$H^\ast(X,\Z)$ of $X$, which has rank $24$, rather
than on $H^2(X,\Z)$ (see \cite{tawe11}); $G$ acts trivially on the orthogonal complement of $H^2(X,\Z)$
in $H^\ast(X,\Z)$, making this modification natural indeed.
\bigskip

When in Kondo's construction,
one uses $M_G(-1)\hookrightarrow N$ and extends the  automorphisms of $M_G(-1)$ that 
are induced
by symmetries of $X$ trivially to the orthogonal complement of $M_G(-1)$ in
$N$, this ensures that one obtains an action of $G$ on $N$.
As already follows from what is said in Kondo's work \cite[Proof of Theorem 4]{ko98a}, 
the projection from $\Aut(N)$ to $M_{12}$ induces an 
injective group homomorphism $G\hookrightarrow M_{12}$. 
Indeed, when retracing our steps from $\Aut(N)$ to $\Aut( H^2(X,\Z) )$,
any $\wt\Phi\in \Aut(N)$ which acts trivially on the orthogonal complement of 
$M_G(-1)$  induces 
a unique automorphism $\Phi$ of $H^2(X,\Z)$ that on $M_G$ agrees with $\wt\Phi$ and
that acts trivially on the orthogonal complement of $M_G$.
Recall from  appendix \ref{subapp:niemeier}
that $\Aut(N)$ possesses a normal subgroup $G_0\times G_1$, 
where $G_0$
is the Weyl group of the  root sublattice of $N$ of type $A_2^{12}$, and $G_1\cong\Z_2$ is generated by 
the involution which acts by multiplication by $-1$ on $N$; moreover, the quotient group 
$\Aut(N)/(G_0\times G_1)$ is isomorphic to $M_{12}$. 
If $\wt\Phi\in G_0\times G_1$, then $\Phi$ is $\Z_3$-effective if
and only if it is the identity, by the definition of the notion of $\Z_3$-effectiveness. 
More generally, by the same line of thought,
if a symmetry of $X$ induces $\Phi\in\Aut(H^2(X,\Z))$ and
$\wt\Phi\in\Aut(N)$, then the $\Z_3$-effectiveness of $\Phi$ implies that
no two elements of the 
$(G_0\times G_1)$-orbit of $\wt\Phi$ can be induced by different symmetries of $X$.
This is a crucial step in Kondo's proof of the fact that 
any symmetry group $H$ of a K3 surface can be viewed 
as a subgroup of $M_{12}$ if $M_{H}(-1)$ can be primitively embedded in the Niemeier lattice
of type $A_2^{12}$, see \cite[Proof of Theorem 4]{ko98a}.

Below, by way of the variation of Kondo's techniques mentioned above,
we construct an action of the symmetry group $G$ of $X$
on $N$ which restricts to the representation induced on $\wt P$ 
by the action of the  symmetries of $X$ on $H^2(X,\Z)$. The induced action of each
element of $G$ is uniquely defined up to a contribution from 
the Weyl group $S_3^3$ of the lattice $\wt K$, a root lattice of type $A_2^3$.
In contrast to the setting in Kondo's construction, the induced action on the discriminant 
group of $\wt P$ is not always trivial, meaning that we cannot always choose the identity
to consistently extend our lattice automorphisms to $\wt K$. 
But we will see that in our setting it is straightforward to make choices which indeed yield a
$G$-action, resulting in an injective group homomorphism $G\hookrightarrow\Aut(N)$.
Moreover, the choices made when  inducing an automorphism of $N$  
from the action of a symmetry of $X$ on $H^2(X,\Z)$ take
place solely within the group $G_0\times G_1$. Hence  the projection from 
$\Aut(N)$ to $M_{12}$ is independent of these choices.
\bigskip

We begin by finding a representation of the translational symmetry
group of $X$ on the lattice $N$:
\bprop\label{translationsonN}%
The translational symmetry group of $X$ acts faithfully on the Niemeier lattice $N$,
such that the resulting subgroup of $\Aut(N)$ is isomorphic to $(\Z_3)^2$ with
generators $\wt\alpha_\ast^1$, $\wt\alpha_\ast^2$,
which are induced by the symmetries
$\alpha^1,\, \alpha^2$ given in \eqref{eq:gammas}, and which
are obtained by $\Q$-linearly extending the following maps:
\begin{eqnarray*}
\mbox{for } \ell\in\{1,2\},\\
\wt\alpha_\ast^1\colon
&& \hspace*{-1em}\left\{
\begin{array}{rcrrcrrcr}
\wt E_1^{(\ell)}&\longmapsto&\wt E_4^{(\ell)}\, ,&
\wt E_2^{(\ell)}&\longmapsto&-\wt E_9^{(\ell)}\, ,&
\wt E_3^{(\ell)}&\longmapsto&-\wt E_1^{(\ell)}\, ,\\[0.3em]
\wt E_4^{(\ell)}&\longmapsto&-\wt E_3^{(\ell)}\, ,&
\wt E_5^{(\ell)}&\longmapsto&-\wt E_8^{(\ell)}\, ,&
\wt E_6^{(\ell)}&\longmapsto&-\wt E_5^{(\ell)}\, ,\\[0.3em]
\wt E_7^{(\ell)}&\longmapsto&-\wt E_2^{(\ell)}\, ,&
\wt E_8^{(\ell)}&\longmapsto&\wt E_6^{(\ell)}\, ,&
\wt E_9^{(\ell)}&\longmapsto&\wt E_7^{(\ell)}\, ,\\[0.3em]
\wt E_{10}^{(\ell)}&\longmapsto&\wt E_{10}^{(\ell)}\, ,&
\wt E_{11}^{(\ell)}&\longmapsto&\wt E_{11}^{(\ell)}\, ,&
\wt E_{12}^{(\ell)}&\longmapsto&\wt E_{12}^{(\ell)}\, ;\\
\end{array}
\right.\\[1em]
\wt\alpha_\ast^2\colon
&&\hspace*{-1em}\left\{
\begin{array}{rcrrcrrcr}
\wt E_1^{(\ell)}&\longmapsto&\hphantom{-}\wt E_9^{(\ell)}\, ,&
\wt E_2^{(\ell)}&\longmapsto&\hphantom{-}\wt E_5^{(\ell)}\, ,&
\wt E_3^{(\ell)}&\longmapsto&\hphantom{-}\wt E_2^{(\ell)}\, ,\\[0.3em]
\wt E_4^{(\ell)}&\longmapsto&\wt E_7^{(\ell)}\, ,&
\wt E_5^{(\ell)}&\longmapsto&\wt E_3^{(\ell)}\, ,&
\wt E_6^{(\ell)}&\longmapsto&\wt E_4^{(\ell)}\, ,\\[0.3em]
\wt E_7^{(\ell)}&\longmapsto&\wt E_6^{(\ell)}\, ,&
\wt E_8^{(\ell)}&\longmapsto&\wt E_1^{(\ell)}\, ,&
\wt E_9^{(\ell)}&\longmapsto&\wt E_8^{(\ell)}\, ,\\[0.3em]
\wt E_{10}^{(\ell)}&\longmapsto&\wt E_{10}^{(\ell)}\, ,&
\wt E_{11}^{(\ell)}&\longmapsto&\wt E_{11}^{(\ell)}\, ,&
\wt E_{12}^{(\ell)}&\longmapsto&\wt E_{12}^{(\ell)}\, .\\
\end{array}
\right.
\end{eqnarray*}
\eprop

\bpf
With $\iota\colon P(-1)\hookrightarrow N$ as in proposition \ref{embeddingexists}
and $\alpha_\ast^k$, $k\in\{1,2\}$, as in lemma \ref{translationalaction}, we require 
$\wt\alpha_\ast^k\circ\iota=\iota\circ\alpha_\ast^k$. This uniquely
determines the action of the $\wt\alpha_\ast^k$ on $\wt P$. By a direct
calculation one checks that this proves all claims made for the $\wt\alpha_\ast^k$-images of the 
$\wt E_j^{(\ell)}$ with $j\in\{1,\ldots, 9\}$ and $\ell\in\{1,2\}$ in the statement
of the proposition.

Since the translational symmetries induce a trivial action on the 
discriminant group $\wt P^\ast/\wt P$, we may extend these lattice automorphisms
trivially to the orthogonal complement $\wt K$ of $\wt P$, in accord with the claims. 
The $\wt\alpha_\ast^k$ thus yield well-defined lattice automorphisms
that are induced by the respective translational symmetries.

As a consistency check, one may convince oneself of the fact that 
the $\Q$-linear extensions of these maps indeed map the glue vectors $(w_1,\ldots,w_6)$
to another choice of glue vectors via
$$
\wt\alpha_\ast^1 \left(\begin{matrix}w_1\\w_2\\w_3\\w_4\\w_5\\w_6\end{matrix}\right)
= \left(\begin{matrix}-w_2-w_3\\w_2\\w_1+w_2-w_3\\w_4-w_2+w_3\\w_3+w_5\\w_6-w_2-w_3\end{matrix}\right)\,;
\qquad
\wt\alpha_\ast^2 \left(\begin{matrix}w_1\\w_2\\w_3\\w_4\\w_5\\w_6\end{matrix}\right)
= \left(\begin{matrix}w_1\\-w_1-w_3\\w_1+w_2-w_3\\w_4-w_3\\w_1-w_3+w_5\\w_6-w_1-w_3\end{matrix}\right)\,.
$$
Since the $\wt\alpha_\ast^k$ act trivially on the orthogonal complement $\wt K$ of $\wt P$, 
it follows that they define a faithful representation of the translational symmetry group of $X$ on 
the lattice $N$.
\epf

We may now extend the action of the translational symmetry group to
an action of the full symmetry group of $X$ on the lattice $N$:

\bprop\label{rotationsonN}%
The symmetry group of $(X,\omega)$ acts faithfully on the Niemeier lattice $N$.
The resulting subgroup of $\Aut(N)$, which is isomorphic to $(\Z_3)^2\rtimes\Z_4$,
can be generated by the automorphisms
$\wt\alpha_\ast^1,\, \wt\alpha_\ast^2$ obtained in proposition \ref{translationsonN},
along with the automorphism $\wt\beta_\ast$ which is induced by the rotational symmetry
$\beta$ of \eqref{eq:beta} and which arises 
by $\Q$-linearly extending the following map:
$$
\begin{array}{rl}
\mbox{for } \ell\in\{1,2\},\\[0.5em]
\wt\beta_\ast\colon\qquad
&\hspace{-2em} \left\{
\begin{array}{rcrrcrrcr}
\wt E_1^{(\ell)}&\longmapsto&\wt E_9^{(3-\ell)}\, ,&
\wt E_2^{(\ell)}&\longmapsto&\wt E_3^{(3-\ell)}\, ,&
\wt E_3^{(\ell)}&\longmapsto&-\wt E_1^{(3-\ell)}\, ,\\[0.3em]
\wt E_4^{(\ell)}&\longmapsto&\wt E_8^{(3-\ell)}\, ,&
\wt E_5^{(\ell)}&\longmapsto&-\wt E_4^{(3-\ell)}\, ,&
\wt E_6^{(\ell)}&\longmapsto&\wt E_6^{(3-\ell)}\, ,\\[0.3em]
\wt E_7^{(\ell)}&\longmapsto&-\wt E_5^{(3-\ell)}\, ,&
\wt E_8^{(\ell)}&\longmapsto&\wt E_7^{(3-\ell)}\, ,&
\wt E_9^{(\ell)}&\longmapsto&-\wt E_2^{(3-\ell)}\, ,\\[0.3em]
\wt E_{10}^{(\ell)}&\longmapsto&\wt E_{10}^{(\ell)}\hphantom{-}\, ,&
\wt E_{11}^{(\ell)}&\longmapsto&\wt E_{11}^{(\ell)}\hphantom{-}\, ,&
\wt E_{12}^{(\ell)}&\longmapsto&\wt E_{12}^{(3-\ell)}\, .\\
\end{array}
\right.
\end{array}
$$
\eprop

\bpf
Analogously to the proof of  proposition \ref{translationsonN}, the requirement
$\wt\beta_\ast\circ\iota=\iota\circ\beta_\ast$ with $\beta_\ast$ 
as in lemma \ref{rotationalaction} uniquely
determines the action of $\wt\beta_\ast$ on $\wt P$. By a direct
calculation one checks that this proves all claims made for the 
$\wt\beta_\ast(E_j^{(\ell)})$ with $j\in\{1,\ldots, 9\}$ and $\ell\in\{1,2\}$ in the statement
of the proposition.

From \eqref{betaondisc} we conclude that on the 
discriminant group $\wt P^\ast/\wt P$, the automorphism $\wt\beta_\ast$
induces an isomorphism $\ol\beta$ which
acts by multiplication by $-1$ on the generator $\wt p_2+\wt P$
of $\wt P^\ast/\wt P$,  while
it acts trivially on the generators $\wt p_j+\wt P$ with $j\in\{1,3\}$.
Let $\ol\psi$ denote the automorphism 
of the discriminant group $\wt K^\ast/\wt K$ of the orthogonal
complement $\wt K$ of $\wt P$ which obeys 
$\wt\gamma\circ\ol\beta=\ol\psi\circ\wt\gamma$, where $\wt\gamma$ 
is the gluing isomorphism given in lemma \ref{gluinginN}.
To consistently extend $\wt\beta_\ast$ to the entire lattice $N$, by
\eqref{autoextend} we must find an automorphism $\wt\psi$ of $\wt K$ 
which induces $\ol\psi$.
One checks that by the above, $\ol\psi$ 
acts by multiplication by $-1$ on the generator $\frac{1}{3}\wt E_{12}+\wt K$
of $\wt K^\ast/\wt K$,  while
it acts trivially on the generators $\frac{1}{3}\wt E_{j}+\wt P$ with $j\in\{10,11\}$.
Thus the map $\wt\psi\colon (\wt E_{10}^{(\ell)},\, \wt E_{11}^{(\ell)},\,
\wt E_{12}^{(\ell)})\longmapsto (\wt E_{10}^{(\ell)},\, \wt E_{11}^{(\ell)},\,
\wt E_{12}^{(3-\ell)})$ for $\ell\in\{1,2\}$ is found to $\Z$-linearly extend 
to an automorphism 
of $\wt K$ which induces $\ol\psi$. In summary, $\wt\beta_\ast\in\Aut(N)$ 
is well-defined, and it is induced by the symmetry $\beta$ of $X$.

As a consistency check, one may convince oneself of the fact that 
$\wt\beta_\ast$ indeed maps the  glue vectors $(w_1,\ldots,w_6)$
to another choice of glue vectors via
$$
\wt\beta_\ast \left(\begin{matrix}w_1\\w_2\\w_3\\w_4\\w_5\\w_6\end{matrix}\right)
= \left(\begin{matrix}-w_2\\-w_1+w_2+w_3\\-w_2+w_3\\w_3+w_4\\w_5-w_2+w_3\\w_2+w_3-w_6\end{matrix}\right)\,.
$$

From our claim, it is immediate that $\wt\beta_\ast$ has the same order $4$ as $\beta$;
since the translational symmetries act by $\wt\alpha_\ast^k$ and thus trivially on $\wt K$, 
it follows that altogether we have obtained a faithful representation of the  symmetry group 
of $X$ on the lattice $N$.
\epf
\brmk\label{symmetrybreaking}%
By remark \ref{E64}, there exists an embedding of $P(-1)$ into the Niemeier lattice $\wh N$
of type $E_6^4$, which however is not primitive. It is therefore a priori not clear whether 
or not a given symmetry of $X$ induces a well-defined lattice automorphism of $\wh N$. By 
what is said in appendix \ref{app:npembedding}, the smallest primitive sublattice $\wh P^\ast\cap\wh N$
of $\wh N$ which contains
the image $\wh P$ of $P(-1)$ is obtained from $\wh P$ by including one additional generator, 
denoted $\wh r$  in appendix \ref{app:npembedding}. Every lattice automorphism of 
$H^2(X,\Z)$ that is induced by a symmetry of $X$ induces a lattice automorphism of $\wh P$
which can be $\Q$-linearly extended to $\wh P\otimes_\Z\Q$. This yields an
automorphism of $\wh P^\ast\cap\wh N$ if and only if $\wh r$
is mapped into $\wh P^\ast\cap\wh N$. For translational symmetries, one 
checks that this is the case and that consistent continuations to automorphisms of $\wh N$ exist. 
However, one also checks  that the automorphism $\wh\beta_\ast$ that is induced by $\beta$ yields
$\langle \wh\beta_\ast\wh r,\wh r\rangle\not\in\Z$, which implies 
$\wh\beta_\ast\wh r\not\in \wh N$, while $\beta_\ast^2\wh r-\wh r\in\wh P$.  
In summary, the embedding of $P(-1)$ into the Niemeier lattice of type $E_6^4$ breaks the symmetry
group $(\Z_3)^2\rtimes\Z_4$ of $(X,\omega)$ to $(\Z_3)^2\rtimes\Z_2$, which is the
symmetry group of $(X,\wt\omega)$ with $\wt\omega=V_3\kappa_3+V_4\kappa_4$, $V_3,\, V_4\in\R$, $V_3V_4>0$ and $V_3\neq V_4$. 
As was explained in remark \ref{E64}, 
this is consistent with what one might expect from studying the elliptic fibration of
our K3 surface $X$ which is induced from the decomposition of the underlying complex torus
$T$ into a product of two $\Z_3$-symmetric elliptic curves.
\ermk
\bigskip

Let us return to our aim of tracking the symmetries of $X$ in the Mathieu group $M_{12}$.
As was mentioned in appendix \ref{subapp:niemeier}, we view $M_{12}$ as 
the quotient group of the automorphism 
group of the extended ternary
Golay code $\CCC_{12}\subset\F_3^{12}$ by its centre $\Z_2$.
$M_{12}$ acts by permutations on the components
of the codewords in $\F_3^{12}$. 
Generators of the code $\CCC_{12}$ are obtained from the 
rows of the matrix $\m G$ of \eqref{matrixcode} by interpreting them 
as vectors in $\F_3^{12}$; each permutation in $M_{12}$ can be lifted 
to an automorphisms of $\m C_{12}$ by accompanying it with appropriate sign flips
``$+\leftrightarrow-$'' on the $\F_3$ components in $\F_3^{12}$.
However, it is important to keep in mind that the automorphism group of $\m C_{12}$ 
is a \emph{non-split} extension of $M_{12}$, so that the lift of $M_{12}$ to 
the automorphism group of $\m C_{12}$ cannot be achieved by means of a 
group homomorphism. 
In propositions \ref{translationsonN} and \ref{rotationsonN}, we have already 
expressed the symmetry group of $X$ as a subgroup of
$\Aut(N)$. As explained at the beginning of this subsection, there is a projection 
$\Aut(N)\longrightarrow M_{12}$ which on that subgroup restricts to an injective map,
and which in addition yields a result that is independent of the choices that were
made along the way. For any $\wt\Phi\in\Aut(N)$, the image under the
projection to $M_{12}$ is the permutation that $\wt\Phi$ induces 
on the twelve irreducible root sublattices of type $A_2$ in $N$, denoted
$A_{2,j}$, $j\in\{1,\ldots,12\}$. We may therefore directly
read off the relevant permutations from the results of propositions 
\ref{translationsonN} and \ref{rotationsonN}:
\bcor\label{subgroupofM12}%
The above construction leads to an injective group homomorphism $\Theta$ mapping 
the symmetry group of $X$ into the Mathieu group $M_{12}$. The image is generated by
the following three permutations:
\begin{eqnarray*}
\Theta(\alpha^1)=   (1, 4, 3) (2, 9, 7)(5, 8, 6)\,,\qquad
\Theta(\alpha^2)&=& ( 1, 9,8) ( 2, 5, 3) (4, 7,6)\,,\\
\Theta(\beta)&=& (1, 9, 2, 3)(4, 8, 7, 5) \,.\qquad
\end{eqnarray*}
\ecor

\subsection{Tracking the symmetries of \texorpdfstring{$\Z_3$}{TEXT}-orbifold K3s within the 
Mathieu group \texorpdfstring{$M_{24}$}{TEXT}}\label{subsec:symmetriesinM24}%
According to \cite[\S10.2.3, Theorem 15]{conway1998sphere}, there is yet another construction
of the Mathieu group $M_{12}$, namely as subgroup of the Mathieu group $M_{24}$.
The latter  is the automorphism group of the extended binary Golay code $\m C_{24}\subset\F_2^{24}$,
as we briefly recall in appendix \ref{subapp:m24}, where we also state our conventions for the
description of $\m C_{24}$ and $M_{24}$. 
Consider a \textsc{duum} in $\m C_{24}$, that is, 
a pair of complementary dodecads.
Then the subgroup of $M_{24}$ given by those permutations that
either fix the two complementary dodecads or interchange them
is a maximal subgroup isomorphic to
$M_{12}\rtimes \Z_2$, and the Mathieu group $M_{12}$ arises as the normal subgroup of this group
which fixes each of the dodecads. Since $M_{24}$ acts transitively on the dodecads in
$\m C_{24}$, we may obtain $M_{12}$ as the subgroup of $M_{24}$ which fixes any choice of dodecad.
Because we have already realized the symmetry group of $\Z_3$-orbifold K3s as a subgroup
of $M_{12}$, viewing $M_{12}$ as a subgroup of $M_{24}$ allows us to track the symmetries of $\Z_3$-orbifold 
K3s within the Mathieu group $M_{24}$. We will do so explicitly in this subsection, in a form that is compatible
with the study of symmetries of $\Z_2$-orbifold K3s within the Mathieu group $M_{24}$ which two
of us have carried out previously \cite{tawe11,tawe12,tawe13}.
\clearpage
\bprop\label{M12inM24}%
Consider the dodecad in $\m C_{24}$ whose non-zero entries have labels in 
$\m D:=\{1,2,3,4,6,8,9,12,13,16,18,23\}$.
Assign the labels $\{1,\ldots,12\}$ to the entries of $\m D$
as indicated in the following MOG-array:
$$
\begin{array}{|cc|cc|cc|}
\hline
{8}&\hphantom{\bullet}&{9}&\hphantom{\bullet}&5&\hphantom{\bullet}\\
\hphantom{\bullet}&7&\hphantom{\bullet}&{11}&\hphantom{\bullet}&3\\
\hline
\hphantom{\bullet}&10&\hphantom{\bullet}&{6}&\hphantom{\bullet}&2\\
\hphantom{\bullet}&12&\hphantom{\bullet}&{4}&\hphantom{\bullet}&1\\
\hline \end{array}\, ,
$$
that is, replace the labels in $\m D$ according
to 
\be\label{relabelD}
(12, 8, 18, 13, 2, 16, 3, 23, 1, 6,4,9)
\quad \longmapsto \quad 
(1, \ldots, 12)\, .
\ee
Then the  
subgroup of $M_{24}$ which stabilizes the set $\m D$
yields the Mathieu group $M_{12}$ according to our conventions.
In other words, under \eqref{relabelD}, the stabilizer group of the dodecad corresponding to $\m D$
becomes a permutation group on $\{1,\ldots,12\}$,
 which is generated by $A, B, C, D$ as listed in \eqref{generatingM12}; vice versa,  
the action of every $\sigma\in M_{12}$ induces a permutation of the
twelve labels in $\m D$
which allows an accompanying permutation of the
complement dodecad with labels in $\{5,7,10,11,14,15,17,19,20,21,22,24\}$
that extends $\sigma$ to an element of $M_{24}$.
\eprop
\bpf
Using the MOG as described in appendix \ref{subapp:m24}, one first
checks that the set $\m D$ indeed corresponds to a dodecad in $\m C_{24}$. In fact,
by what was said above, we might work with any dodecad; 
however to realize elements of $M_{12}$ in terms of $M_{24}$, we must then find a
relabelling of the twelve non-zero entries of our dodecad
by $\{1, \ldots, 12\}$, such that the permutation of the set $\{1, \ldots, 12\}$
by each $\sigma\in M_{12}$ induces a permutation of the
twelve labels in the chosen dodecad which in turn allows an accompanying permutation of the
complement dodecad that extends $\sigma$ to an element of $M_{24}$. This is a non-trivial 
step which we perform with the help of the results presented by Conway and Sloane in
\cite{conway1998sphere}, where several approaches are described. Here, for our 
choice of dodecad we follow the suggestion of \cite[\S11.17, Figure 11.37]{conway1998sphere},
namely
$$
\begin{array}{|cc|cc|cc|}
\hline
{\bullet}&\hphantom{\bullet}&{\bullet}&\hphantom{\bullet}&\bullet&\hphantom{\bullet}\\
\hphantom{\bullet}&\bullet&\hphantom{\bullet}&{\bullet}&\hphantom{\bullet}&\bullet\\
\hline
\hphantom{\bullet}&\bullet&\hphantom{\bullet}&{\bullet}&\hphantom{\bullet}&\bullet\\
\hphantom{\bullet}&\bullet&\hphantom{\bullet}&{\bullet}&\hphantom{\bullet}&\bullet\\
\hline \end{array}
\quad\mbox{ with ``shuffle numbering'' }\qquad
\begin{array}{|cc|cc|cc|}
\hline
{0}&\hphantom{\bullet}&{1}&\hphantom{\bullet}&2&\hphantom{\bullet}\\
\hphantom{\bullet}&3&\hphantom{\bullet}&{4}&\hphantom{\bullet}&5\\
\hline
\hphantom{\bullet}&6&\hphantom{\bullet}&{7}&\hphantom{\bullet}&8\\
\hphantom{\bullet}&9&\hphantom{\bullet}&{10}&\hphantom{\bullet}&11\\
\hline \end{array}\; .
$$
In \cite[\S11.18 (page 328)]{conway1998sphere} we are given a dictionary 
between the ``shuffle numbering'' and the ``numbering $\mod 11$'',
which in turn we have already translated into our own labeling by $\{1,\ldots,12\}$ in 
\eqref{relabelCSToBTUWZ}, as follows:
$$
\begin{array}{||l||c|c|c|c|c|c|c|c|c|c|c|c||}
\hline\hline
\mbox{shuffle numbering}&0&1&2&3&4&5&6&7&8&9&10&11\\
\hline
\mbox{numbering $\mod 11$}&\infty&1&9&3&4&5&0&8&6&2&10&7\\
\hline
\mbox{our numbering}&8&9&5&7&11&3&10&6&2&12&4&1\\
\hline\hline
\end{array}\, .
$$
Applying this dictionary to the above MOG with ``shuffle numbering'' yields 
the MOG-labelling claimed in the proposition. By \eqref{ConwayMOG}, we obtain the 
relabelling of the entries of $\m D$ as claimed. The remaining
claims now follow from what is said in \cite[\S11.17--11.18]{conway1998sphere},
and they are immediately confirmed by GAP.
\epf

Proposition \ref{M12inM24} shows that there exists an embedding
of $M_{12}\hookrightarrow M_{24}$ such that the induced action of each $\sigma\in M_{12}$
on the labels of our chosen dodecad $\m D$ is obtained by means of \eqref{relabelD}. 
Hence one might calculate the image of each of the generators $A,B,C,D$ of $M_{12}$
within $M_{24}$, express $\sigma$ in terms of these generators and thereby obtain
the image of $\sigma$ in $M_{24}$. However, our main interest is in the symmetries
of $\Z_3$-orbifold K3s; we therefore calculate their images in $M_{24}$ without
the use of $A,B,C,D$:

\bprop\label{imageinM24}%
Let $\vartheta$ denote the group homomorphism that maps the symmetry group of 
our $\Z_3$-orbifold K3 surface $(X,\omega)$ into the Mathieu group $M_{24}$ and which arises by 
composing the group homomorphism $\Theta$ of corollary \ref{subgroupofM12} with the embedding
$M_{12}\hookrightarrow M_{24}$ of proposition \ref{M12inM24}. Then the image of the symmetry
group of $(X,\omega)$ in $M_{24}$ is
generated by 
\begin{eqnarray*}
\vartheta(\alpha^1)&=&   (1,3,8)(2,23,16) (5,22,14)(10,20,21)(11,17,24)(12,13,18)\,,\\
\vartheta(\alpha^2)&=& (1,23,12)(2,18,8) (3,16,13)(5,14,22)(7,15,19)(10,20,21)\,,\\
\vartheta(\beta)&=& (1,8,18,12)(2,13,23,3) (5,10)(7,11,15,24)(14,21,22,20)(17,19) \,.\qquad
\end{eqnarray*}
\eprop

\bpf
By corollary \ref{subgroupofM12}, $\vartheta(\alpha^1)$, $\vartheta(\alpha^2)$, 
$\vartheta(\beta)$ generate the desired subgroup of $M_{24}$. Hence we must prove 
that these symmetries are described by the permutations as claimed.
For $\sigma\in\{\alpha^1, \alpha^2, \beta\}$ one first checks by a direct calculation
that the claimed permutation $\vartheta(\sigma)$ maps the dodecad $\m D$ 
that was used in proposition \ref{M12inM24} to itself by permuting $\m D$ according to the 
image of $\Theta(\sigma)\in M_{12}$ under \eqref{relabelD}.
The general strategy by which we confirm the claimed action of $\vartheta(\sigma)$ on 
the entire set $\{1,\ldots, 24\}$ uses the special octads listed in \eqref{domino}: in every step, 
the action on a subset $\m D^\prime\subset \{1,\ldots, 24\}$ is already known, where 
$\m D\subset\m D^\prime$. We choose an octad from the list \eqref{domino},
such that the labels with non-zero entries overlap with $\m D^\prime$ in at least
$5$ and at most $7$ entries. Then the image of this octad under $\vartheta(\sigma)$ is an octad
of which at least $5$ of the labels of non-zero entries are known; these uniquely determine the image octad
according to \cite[\S11.6]{conway1998sphere}. By making appropriate choices from \eqref{domino},
in each case we obtain sufficiently many restrictions on the possible action of 
$\vartheta(\sigma)$ to uniquely determine the image.

Let us illustrate the procedure by determining the image of $\alpha^1$ in $M_{24}$. We begin 
with the octad with label $10$ in \eqref{domino}; according to \eqref{ConwayMOG}, for this octad
the list of labels of non-zero entries overlaps with $\m D$ 
in $\{3,4,6,9,13,16\}$ which is mapped to $\{2,4,6,8,9,18\}$ under $\vartheta(\alpha^1)$;
the unique octad with non-zero entries labelled by this latter list carries the number $5$ in \eqref{domino},
and the remaining non-zero entries have labels $\{11,17\}$. We conclude that the remaining labels 
with non-zero entries of the octad number $10$, that is, $\{11,24\}$,
are mapped to  $\{11,17\}$ under $\vartheta(\alpha^1)$.
The calculation is more efficiently depicted by the following diagrams, where we include the MOG-labels
from \eqref{domino} and record the image of each octad under $\vartheta(\alpha^1)$; star entries 
$\star$ signify entries whose image under $\vartheta(\alpha^1)$ has not yet been determined:
$$
\begin{array}{|cc|cc|cc|}
\hline
\hphantom{\bullet}&{\star}&\hphantom{\bullet}&{\star}&\hphantom{\bullet}&\hphantom{\bullet}\\
\hphantom{\bullet}&{{\bullet}}&\hphantom{\bullet}&{\bullet}&\hphantom{\bullet}&\hphantom{\bullet}\\
\hline
\hphantom{\bullet}&{\bullet}&\hphantom{\bullet}&{\bullet}&\hphantom{\bullet}&\hphantom{\bullet}\\
\hphantom{\bullet}&{\bullet}&\hphantom{\bullet}&{\bullet}&\hphantom{\bullet}&\hphantom{\bullet}\\
\hline \end{array}_{\;10}
\longrightarrow\quad
\begin{array}{|cc|cc|cc|}
\hline
\hphantom{\bullet}&\hphantom{\bullet}&\hphantom{\bullet}&{\star}&{\bullet}&\hphantom{\bullet}\\
\hphantom{\bullet}&\hphantom{\bullet}&\hphantom{\bullet}&{\bullet}&\hphantom{\bullet}&{\bullet}\\
\hline
\hphantom{\bullet}&{\bullet}&\hphantom{\bullet}&\hphantom{\bullet}&{\star}&{\bullet}\\
\hphantom{\bullet}&{\bullet}&\hphantom{\bullet}&\hphantom{\bullet}&\hphantom{\bullet}&\hphantom{\bullet}\\
\hline \end{array}_{\;5}
\longrightarrow\quad
\begin{array}{|cc|cc|cc|}
\hline
{\bullet}&{\star}&{\bullet}&\hphantom{\bullet}&\hphantom{\bullet}&\hphantom{\bullet}\\
\hphantom{\bullet}&\hphantom{\bullet}&\hphantom{\bullet}&{\bullet}&\hphantom{\bullet}&\hphantom{\bullet}\\
\hline
\hphantom{\bullet}&{\bullet}&\hphantom{\bullet}&\hphantom{\bullet}&{\star}&\hphantom{\bullet}\\
\hphantom{\bullet}&{\bullet}&\hphantom{\bullet}&\hphantom{\bullet}&\hphantom{\bullet}&{\bullet}\\
\hline \end{array}_{\;19}\; .
$$
By \eqref{ConwayMOG}, and because $\vartheta(\alpha^1)$ has order $3$, we conclude that 
$\vartheta(\alpha^1)$ yields $\{11,24\}\longrightarrow\{11,17\}\longrightarrow\{17,24\}\longrightarrow\{11,24\}$,
hence $11\mapsto17\mapsto24\mapsto11$.

Analogously, we find
$$
\begin{array}{|cc|cc|cc|}
\hline
\hphantom{\bullet}&\hphantom{\bullet}&\hphantom{\bullet}&{{\bullet}}&\hphantom{{\bullet}}&{\star}\\
\hphantom{\bullet}&\hphantom{\bullet}&\hphantom{\bullet}&{\bullet}&\hphantom{\bullet}&{\bullet}\\
\hline
\hphantom{\bullet}&\hphantom{\bullet}&\hphantom{\bullet}&{\bullet}&\hphantom{\bullet}&{\bullet}\\
\hphantom{\bullet}&\hphantom{\bullet}&\hphantom{\bullet}&{\bullet}&\hphantom{\bullet}&{\bullet}\\
\hline \end{array}_{\;2}
\longrightarrow\quad
\begin{array}{|cc|cc|cc|}
\hline
\hphantom{\bullet}&\hphantom{\bullet}&{\bullet}&\hphantom{\bullet}&{\bullet}&\hphantom{\bullet}\\
\hphantom{\bullet}&\hphantom{\bullet}&\hphantom{\bullet}&{\bullet}&\hphantom{\bullet}&{\bullet}\\
\hline
\hphantom{\bullet}&\hphantom{\bullet}&{\star}&\hphantom{\bullet}&{{\bullet}}&\hphantom{\bullet}\\
\hphantom{\bullet}&\hphantom{\bullet}&\hphantom{\bullet}&{\bullet}&\hphantom{\bullet}&{\bullet}\\
\hline \end{array}_{\;4}
\longrightarrow\quad
\begin{array}{|cc|cc|cc|}
\hline
{\bullet}&{{\bullet}}&\hphantom{\bl{\bullet}}&\hphantom{\bullet}&\hphantom{\bullet}&\hphantom{\bullet}\\
\hphantom{\bullet}&{\bullet}&\hphantom{\bullet}&{\bullet}&\hphantom{\bullet}&{\bullet}\\
\hline
\hphantom{\bullet}&\hphantom{\bullet}&\hphantom{\bullet}&\hphantom{\bullet}&\hphantom{\bullet}&\hphantom{\bullet}\\
{\star}&\hphantom{\bullet}&\hphantom{\bullet}&{\bullet}&\hphantom{\bullet}&{\bullet}\\
\hline \end{array}_{\;23} \; ,
$$
which allows us to conclude that $\vartheta(\alpha^1)$ yields
$22\mapsto14\mapsto5\mapsto22$.
From
$$
\begin{array}{|cc|cc|cc|}
\hline
{\bullet}&{\bullet}&{\bullet}&\hphantom{\bullet}&\hphantom{\bullet}&\hphantom{\bullet}\\
\hphantom{\bullet}&{\bullet}&\hphantom{\bullet}&\hphantom{\bullet}&{\star}&\hphantom{\bullet}\\
\hline
\hphantom{\bullet}&{\bullet}&\hphantom{\bullet}&\hphantom{\bullet}&\hphantom{\bullet}&{\bullet}\\
\hphantom{\bullet}&\hphantom{\bullet}&\hphantom{\bullet}&{\bullet}&\hphantom{\bullet}&\hphantom{\bullet}\\
\hline \end{array}_{\;20}%
\longrightarrow\quad
\begin{array}{|cc|cc|cc|}
\hline
\hphantom{\bullet}&\hphantom{\bullet}&{\bullet}&{\bullet}&\hphantom{\bullet}&\hphantom{\bullet}\\
\hphantom{\bullet}&{\bullet}&{\star}&\hphantom{\bullet}&\hphantom{\bullet}&{\bullet}\\
\hline
\hphantom{\bullet}&{\bullet}&\hphantom{\bullet}&{\bullet}&\hphantom{\bullet}&{\bullet}\\
\hphantom{\bullet}&\hphantom{\bullet}&\hphantom{\bullet}&\hphantom{\bullet}&\hphantom{\bullet}&\hphantom{\bullet}\\
\hline \end{array}_{\;7}
\longrightarrow\quad
\begin{array}{|cc|cc|cc|}
\hline
\hphantom{\bullet}&\hphantom{\bullet}&{\bullet}&\hphantom{\bullet}&{\bullet}&\hphantom{\bullet}\\
\hphantom{\bullet}&{\bullet}&\hphantom{\bullet}&\hphantom{\bullet}&\hphantom{\bullet}&\hphantom{\bullet}\\
\hline
\hphantom{\bullet}&{\bullet}&\hphantom{\bullet}&\hphantom{\bullet}&{\bullet}&{\bullet}\\
\hphantom{\bullet}&\hphantom{\bullet}&{\star}&\hphantom{\bullet}&\hphantom{\bullet}&{\bullet}\\
\hline \end{array}_{\;6} \; 
$$
we obtain that $\vartheta(\alpha^1)$ maps $10\mapsto20\mapsto21\mapsto10$, and finally
\begin{eqnarray*}
\begin{array}{|cc|cc|cc|}
\hline
\hphantom{\bullet}&\hphantom{\bullet}&\hphantom{\bullet}&\hphantom{\bullet}&{\bullet}&{\bullet}\\
\hphantom{\bullet}&\hphantom{\bullet}&\hphantom{\bullet}&\hphantom{\bullet}&{\bullet}&{\bullet}\\
\hline
\hphantom{\bullet}&\hphantom{\bullet}&\hphantom{\bullet}&\hphantom{\bullet}&{\bullet}&{\bullet}\\
\hphantom{\bullet}&\hphantom{\bullet}&\hphantom{\bullet}&\hphantom{\bullet}&{\star}&{\bullet}\\
\hline \end{array}_{\;1}
&\longrightarrow&
\begin{array}{|cc|cc|cc|}
\hline
{\bullet}&{\bullet}&{\bullet}&\hphantom{\bullet}&\hphantom{\bullet}&\hphantom{\bullet}\\
\hphantom{\bullet}&\hphantom{\bullet}&{\bullet}&\hphantom{\bullet}&\hphantom{\bullet}&\hphantom{\bullet}\\
\hline
\hphantom{\bullet}&\hphantom{\bullet}&{\bullet}&\hphantom{\bullet}&\hphantom{\bullet}&\hphantom{\bullet}\\
\hphantom{\bullet}&\hphantom{\bullet}&\hphantom{\bullet}&{\bullet}&{\star}&{\bullet}\\
\hline \end{array}_{\;18}
\; ,\\
&&\begin{array}{|cc|cc|cc|}
\hline
\hphantom{\bullet}&\hphantom{\bullet}&{\bullet}&\hphantom{\bullet}&\hphantom{\bullet}&\hphantom{\bullet}\\
\hphantom{\bullet}&{\bullet}&\hphantom{\bullet}&\hphantom{\bullet}&{\bullet}&\hphantom{\bullet}\\
\hline
{\star}&{\bullet}&\hphantom{\bullet}&{\bullet}&\hphantom{\bullet}&\hphantom{\bullet}\\
\hphantom{\bullet}&{\bullet}&\hphantom{\bullet}&\hphantom{\bullet}&\hphantom{\bullet}&{\bullet}\\
\hline \end{array}_{\;13}
\longrightarrow
\begin{array}{|cc|cc|cc|}
\hline
\hphantom{\bullet}&\hphantom{\bullet}&\hphantom{\bullet}&\hphantom{\bullet}&{\bullet}&\hphantom{\bullet}\\
\hphantom{\bullet}&{\bullet}&{\bullet}&\hphantom{\bullet}&\hphantom{\bullet}&\hphantom{\bullet}\\
\hline
{\star}&{\bullet}&\hphantom{\bullet}&\hphantom{\bullet}&\hphantom{\bullet}&{\bullet}\\
\hphantom{\bullet}&{\bullet}&\hphantom{\bullet}&{\bullet}&\hphantom{\bullet}&\hphantom{\bullet}\\
\hline \end{array}_{\;12}
\end{eqnarray*}
shows that $7$ and $15$, and thereby also $19$,
are fixed by $\vartheta(\alpha^1)$. Altogether, we have proven the
claim for the permutation $\vartheta(\alpha^1)$.

To prove the claims for $\vartheta(\alpha^2)$ and $\vartheta(\beta)$, one proceeds analogously. 
Specifically, for $\vartheta(\alpha^2)$, 
we have used the MOGs   in \eqref{domino} labelled
$10$ (which is fixed, implying that elements $11$ and $24$ are fixed), 
$2$ (which maps to octad number $11$ and then on to octad number $17$, 
yielding the permutation $(22,5,14)$, 
$19$ (which is fixed, implying that the element $17$
is fixed), 
$28$ (which maps to octad number $16$ and then on to octad number $14$, 
yielding the permutation $(15,19,7)$) 
and $1$ (mapping to octads numbers 
$29$ and then $27$, yielding the permutation $(10,20,21)$).
For $\vartheta(\beta)$, we have used the MOGs in \eqref{domino} labelled $10$ (which 
maps to octad number $25$ and then back to octad number $10$) and $2$ 
(which maps to octad number $24$, then $8$ and then on to octad number $3$). 
These MOGs and their images under powers of $\vartheta(\beta)$ yield the permutations 
$(11,15,24,7)$ and $(22,20,14,21)$. 
We also used the octads labelled $19$ (which maps to number 
$15$, then $5$ and $22$, yielding the permutation $(17,19)$)
and $28$ (mapping to numbers $9$, $21$ and then $26$, yielding the permutation $(5,10)$).
\epf

We have carefully chosen our conventions for the realization of $M_{12}$ and $M_{24}$ in accord with
those used by two of the authors in their previous work \cite{tawe11,tawe12,tawe13}. This allows us to combine the 
image of the symmetry group of our $\Z_3$-orbifold K3 surface $X$ in $M_{24}$ with the \textsc{octad subgroup}
$(\Z_2)^4\rtimes A_8$ of $M_{24}$ which was obtained by symmetry surfing all $\Z_2$-orbifold K3s in \cite{tawe13},
that is, by combining all symmetry groups of Kummer surfaces that are induced by symmetry groups of the underlying
complex tori. 
While symmetry surfing the moduli space of Kummer surfaces allowed us to naturally
combine symmetry groups from distinct points in moduli space through a common polarization, we
emphasize that the choices made above to  realize our symmetry group of $\Z_3$-orbifold K3s
as a subgroup of $M_{24}$ remain, at this stage,  ad hoc. We leave it for future work to establish choices that are motivated
by geometry or conformal field theory. However, as a proof of concept we remark that as expected,
our newly constructed symmetries extend the maximal subgroup $(\Z_2)^4\rtimes A_8$ 
of $M_{24}$ obtained from symmetry surfing to the entire group $M_{24}$:

\bthm\label{generatingM24}%
The subgroup of type $(\Z_2)^4\rtimes A_8$ of $M_{24}$ which is obtained as combined symmetry group
of all Kummer K3s in \cite[\S4]{tawe13},
together with the image of the symmetry group of $\Z_3$-orbifold K3 surfaces under the map $\vartheta$
obtained in proposition \ref{imageinM24} generate the Mathieu group $M_{24}$.
\ethm

\bpf
While the claim may be confirmed by GAP by using proposition \ref{imageinM24} above as well as the permutations listed in 
\cite[(2.1),(3.1),(3.2),(3.9)]{tawe13},
we find the following argument more elegant and more rewarding.

Recall from \cite[\S2]{tawe13} and \cite[(2.14)]{tawe11} that the subgroup of type $(\Z_2)^4\rtimes A_8$ 
of $M_{24}$ which is obtained as combined symmetry group of all Kummer K3s in \cite[\S4]{tawe13} is the 
subgroup of $M_{24}$ which stabilizes the octad with non-zero entries labelled by $\{3,5,6,9,15,19,23,24\}$,
that is, the octad  with label $28$ in \eqref{domino}. This is a maximal subgroup of $M_{24}$.
Hence any permutation in $M_{24}$ which does not stabilize this octad will suffice to generate $M_{24}$ from 
this group. Proposition \ref{imageinM24} shows that each of the generators $\vartheta(\alpha^1)$,
$\vartheta(\alpha^2)$, $\vartheta(\beta)$ mixes the set $\{3,5,6,9,15,19,23,24\}$ with its complement
in $\{1,\ldots,24\}$, thus proving the claim.
\epf

\acknowledgments
This project was suggested by two of us as one of the research themes for the 2020 Workshop 
\textsl{Women in Mathematical Physics}, Banff International Research Station,
Banff, Canada. We are indebted to the organizers, Ana Ros Camacho and Nezhla Aghaei,
for their initiative, but in particular for their
unfaltering dedication to run this event despite the pandemic --- as an online 
workshop (WoMaP I, September 21--22, 2020), followed by the in-person meeting WoMaP II at the BIRS, Banff, Canada, August 13 -- 18, 2023.
We are also grateful to Katrina Barron and Gaywalee Yamskulna,
who complemented the organising committee for WoMaP II,
and to the Banff International Research Station,
Banff, Canada, for their hospitality during this event,
where some of us met for the first time in person. 
We thank Lyda Urresta for her participation in early stages of the project.
KW is grateful to Kenji Ueno for some crucial hints 
concerning the construction of $\Z_3$-orbifold limits of K3 by  
blowing up a complex torus before quotienting it.
IGZ wishes to thank G.~Aldazabal, E.~Andr\'es, A.~Font, K.~Narain and N.~Scheithauer for helpful discussions.
We thank an anonymous referee for their careful reading of our manuscript and for their helpful comments.

AT and KW would like to thank the Isaac Newton Institute for Mathematical Sciences,
Cambridge, UK, for support and hospitality during the programme 
\textsl{New connections in number theory and physics},
where work on this paper was undertaken; this work was supported by EPSRC grant no EP/R014604/1.
KW thanks the Galileo Galilei Institute for Theoretical Physics at Firenze, Italy,
for its hospitality during the programme 
\textsl{BPS Dynamics and Quantum Mathematics},
and the INFN for partial support during the completion of part of this work.
IGZ thanks ICTP and MITP for support and hospitality during different stages of this work. 
KB's work was supported in part by the Simons Collaboration on Confinement and QCD Strings and the Simons Collaboration on Celestial Holography.

The free open-source mathematics software system  
SAGE\footnote{W.A. Stein et~al., 
\emph{{S}age {M}athematics {S}oftware \mbox{\rm(}{V}ersion 4.6.2}\mbox{\rm)}}, 
The Sage Development Team, 2010, \url{http://www.sagemath.org} and in particular its 
component GAP4\footnote{The GAP~Group, \emph{GAP -- Groups, Algorithms, and Programming,   
Version 4.14.0; 2024.
  \url{https://www.gap-system.org}}}
  were used,
 as well as Wolfram Mathematica,
  in checking several assertions made in the group 
theoretical component of our work.

\appendix
\section{Blowing up \texorpdfstring{$A_2$}{TEXT} singularities}\label{app:A2blowup}
This work investigates properties of $\Z_3$-orbifold K3 surfaces,
which can be obtained as minimal resolutions of $\Z_3$ quotients of 
certain tori, where all singularities in the quotient are of a type called
$A_2$. In this appendix, we recall some properties of these particular singularities,
and we showcase the minimal resolution of singularities of type $A_2$.
Though these techniques are well-known to algebraic geometers
\cite{du34,ar66,hi66}, 
presenting the details 
of the calculation makes our work self-contained, and it also serves us to introduce the notations that are used in the main part of the text.
\bigskip

Consider the following surface $S(A_2)$ in $\C^3$:
$$
S(A_2) := \left\{(x_0, x_1, x_{2})\in \C^3 \mid x_0^3= x_1 x_2\right\} \; .
$$
One checks that this surface is singular at the origin $x=(0,0,0)$ and 
smooth everywhere else; in fact, the origin is a 
\textsc{singularity of type $A_2$}, 
providing a local model of the standard $\Z_3$-orbifold singularity $[0]\in\C^2/\Z_3$, 
as can be seen as follows. First recall the standard $\Z_3$-action on $\C^2$ by
linear unitary, orientation preserving maps, which is generated by \eqref{Z3action}. 
The origin $(z_1,z_2)=(0,0)$ is the unique fixed point
of this action, and it descends to a singular point in $\C^2/\Z_3$.
One has an isomorphism of complex varieties
$$
\C^2/\Z_3 \longrightarrow S(A_2), \qquad
\left[(z_1,z_2)\right]\mapsto (z_1 z_2, z_1^3, z_2^3)\; ,
$$
which maps $[(z_1,z_2)]=[0]$ to $x=(0,0,0)$, justifying our claim
that $S(A_2)$ is a model of $\C^2/\Z_3$.

\begin{figure}[ht]
\begin{center}
\unitlength1em
\hfill
\begin{picture}(5,7)
\thicklines
\put(-2,6){(a)}
\qbezier(0,0)(1,4)(5,5)
\put(3,5.5){$\CCC_{00}^{(1)}$}
\qbezier(5,0)(1,1)(0,5)
\put(3,1){$\CCC_{00}^{(2)}$}
\end{picture}
\hfill
\begin{picture}(5,7)
\thicklines
\put(-2,6){(b)}
\put(1,3){\circle*{0.4}}
\put(3,3){\circle*{0.4}}
\put(1,3){\line(1,0){2}}
\end{picture}
\hspace*{\fill}
\end{center}
\caption{\small (a) Two copies $\CCC_{00}^{(1)}$, $\CCC_{00}^{(2)}$
of $\P^1$, depicted as (complex) lines,
which intersect transversally in one point.
(b) The Dynkin diagram of type $A_2$.}
\label{exceptionaldivisor}
\end{figure}
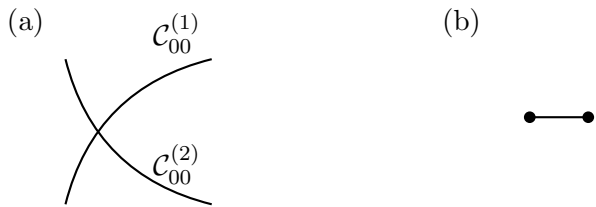

In order to calculate a minimal resolution of an $A_2$-type singularity, the
model $S(A_2)$ is particularly convenient. Indeed, the blow up $\widetilde{S(A_2)}$ 
of $S(A_2)$ in the singular point $x=(0,0,0)$, by definition, is the closure of the 
surface
$$
S^0(A_2)
:=\left\{ (x,u) \in \C^3\times \P^2 \mid 
x\in S(A_2),\; x\neq0,\; x_ju_k=x_ku_j \;\forall j,k\in\{0,\,1,\,2\} \right\}
$$
in $\C^3\times \P^2$, where it is useful to note that 
$S^0(A_2)\cong S(A_2)\setminus\{0\}$ and thus 
$S^0(A_2) \cong \left(\C^2\setminus\{0\}\right)/\Z_3$. 
To exhibit the geometry of $\widetilde{S(A_2)}$,
we consider $U_j$ with $j\in\{0,\,1,\,2\}$, the 
intersection of $\widetilde{S(A_2)}$ with the $j^{\rm th}$ standard
affine coordinate chart 
$$
\left\{ (x,u) \in \C^3\times \P^2 \mid u=(u_0\colon u_1\colon u_2), u_j\neq 0\right\}\; .
$$
It is the closure of the intersection of $S^0(A_2)$ with the $j^{\rm th}$ standard
affine coordinate chart in $\C^3\times \P^2$:
\begin{eqnarray}
U_0 &=& \left\{ (x,u) \in \C^3\times \P^2 \mid u_0=1, \; x_1 = u_1 x_0, \; 
x_2= u_2x_0, \; x_0=u_1u_2 \right\} \label{U0}\\
&& \qquad = \left\{ \left( (u_1u_2, u_1^2u_2, u_1u_2^2), (1\colon u_1\colon u_2)\right) \in \C^3\times \P^2 \mid  (u_1, u_2)\in\C^2\right\} ,
\nonumber\\
U_1 &=& \left\{ (x,u) \in \C^3\times \P^2 \mid u_1=1, \; x_0 = u_0 x_1, \; 
x_2= u_2x_1, \; u_2=u_0^3x_1 \right\} 
\label{U1}\\
&& \qquad = \left\{ \left( (u_0x_1, x_1, u_0^3x_1^2), (u_0\colon 1\colon u_0^3x_1)\right) \in \C^3\times \P^2 \mid  (u_0, x_1)\in\C^2\right\}  \; ,\nonumber\\
U_2 &=& \left\{ (x,u) \in \C^3\times \P^2 \mid u_2=1, \; x_0 = u_0 x_2, \; 
x_1= u_1x_2, \; u_1=u_0^3x_2 \right\} 
\label{U2}\\
&& \qquad = \left\{ \left( (u_0x_2, u_0^3x_2^2, x_2), (u_0\colon u_0^3x_2\colon 1)\right) \in \C^3\times \P^2 \mid  (u_0, x_2)\in\C^2\right\}\; .
\nonumber
\end{eqnarray}
The above gives a smooth parametrization for each $U_j$; thus each $U_j$ is smooth,
and we may confirm that $\widetilde{S(A_2)}$ is a smooth surface.
We may now also read off the exceptional divisor
of the blow-up, $\CCC_{00}:=\left\{(x,u)\in \widetilde{S(A_2)} \mid x=0\right\}$,
where we find $\CCC_{00}=\CCC_{00}^{(1)}\cup \CCC_{00}^{(2)}$ with
\be\label{exceptionalcomponents}
\begin{array}{rcl}
\CCC_{00}^{(1)}&:=&\left\{\left(0, (u_0\colon u_1\colon 0)\right) \in \C^3\times \P^2 
\mid (u_0\colon u_1)\in\P^1\right\}, \\[5pt]
\CCC_{00}^{(2)}&:=&\left\{\left(0, (u_0\colon 0\colon u_2)\right) \in \C^3\times \P^2 
\mid (u_0\colon u_2)\in\P^1\right\}\; .
\end{array}
\ee
In other words, $\CCC_{00}$ decomposes into two irreducible components, 
$\CCC_{00}^{(1)}$ and $\CCC_{00}^{(2)}$, where each of these components is isomorphic
to $\P^1$, and these two irreducible components of $\CCC_{00}$ have one
intersection point $\left( 0, (1\colon0\colon0)\right)$, where they
intersect transversally (see figure~\ref{exceptionaldivisor}(a)). 
One checks that each  $\CCC_{00}^{(k)}$ is a $(-2)$-curve,
that is, it has self-intersection number $-2$.
\bigskip

The form of the exceptional divisor in the minimal resolution of $S(A_2)$
which we have obtained above also explains why the singular point in $S(A_2)$
is called of type $A_2$. One uses a graphical depiction of the exceptional divisor, 
by introducing a node corresponding to each irreducible component of the
exceptional divisor and by introducing an edge between two nodes if 
and only if the corresponding curves intersect. 
Following this 
instruction, the exceptional divisor 
$\CCC_{00}=\CCC_{00}^{(1)}\cup \CCC_{00}^{(2)}$ is depicted by a graph 
with two nodes, and with one edge, connecting these two nodes
(see figure~\ref{exceptionaldivisor}(b)). 
This is the standard Dynkin diagram of type $A_2$ -- thus the terminology.
\bigskip

We finally claim that the resolution $\widetilde{S(A_2)}$ constructed
above is in fact a so--called \textsc{crepant resolution}, which is 
equivalent to the claim that on $\widetilde{S(A_2)}$, there exists
a holomorphic $(2,0)$-form $\Omega$ which vanishes nowhere. Indeed, the
existence of such a form $\Omega$
can be checked by direct calculation. To this end, note that 
$dz^1\wedge dz^2$ is a holomorphic $(2,0)$-form on $\C^2$ which vanishes
nowhere. It descends to a holomorphic $(2,0)$-form on 
$S^0(A_2) \cong \left(\C^2\setminus\{0\}\right)/\Z_3$ since 
it is invariant under the $\Z_3$-action \eqref{Z3action}. One now checks
that the resulting form may be extended to a holomorphic 
$(2,0)$-form $\Omega$ on $\widetilde{S(A_2)}$ which vanishes nowhere.
Indeed, $x_0=z_1z_2$, $x_1=z_1^3$, $x_2=z_2^3$ implies 
$$
dx_1\wedge dx_2 = 9 x_0^2\; dz_1\wedge dz_2  \; ,
$$
such that on $U_0$, where by \eqref{U0}
we use $(u_1,u_2)$ as coordinates 
with $x_0=u_1u_2,\, x_1=u_1^2u_2,\; x_2=u_1u_2^2$, we find that 
the unique holomorphic extension of $dz^1\wedge dz^2$ is expressed by 
\be\label{OmegaonU0}
\Omega_{\mid U_0}= \frac{1}{3} du_1\wedge du_2\, ,
\ee
which indeed never vanishes. Similarly, on $U_1$ 
by \eqref{U1} we use 
$(u_0, x_1)$ as coordinates with $x_0=u_0x_1$,
$x_2=u_0^3x_1^2$, and thus 
the unique holomorphic extension of $dz^1\wedge dz^2$ is expressed by 
\be\label{OmegaonU1}
\Omega_{\mid U_1}
=\frac{1}{3} dx_1\wedge du_0\, ,
\ee
which also never vanishes. Finally, on $U_2$ by \eqref{U2} we use 
$(u_0, x_2)$ as coordinates with $x_0=u_0x_2$,
$x_1=u_0^3x_2^2$, and thus 
the unique holomorphic extension of $dz^1\wedge dz^2$ is expressed by 
\be\label{OmegaonU2}
\Omega_{\mid U_2}
=\frac{1}{3} du_0\wedge dx_2\, ,
\ee
which never vanishes either. Altogether, we have constructed $\Omega$
as claimed.
\section{Elliptic genus of \texorpdfstring{$\Z_3$}{TEXT}-orbifolds}\label{app_ellgen}

In this appendix we compute the elliptic genus of the K3 surface 
constructed from the $T/\Z_3$ quotient using orbifold conformal 
field theory techniques. To this end we compute the elliptic genus 
of the non-linear $\sigma$-model on $T/\Z_3$. The elliptic genus is 
an invariant of complex manifolds and so is the same for all K3 $\sigma$-models. 
It is 
given by \cite{Eguchi:1988vra}
\be\label{ellk3}
\mathcal E_{\rm K3}(\tau,z)=8\bigg[\Big(\frac{\vartheta_2(\tau,z)}{\vartheta_2(\tau,0)}\Big)^2+\Big(\frac{\vartheta_3(\tau,z)}{\vartheta_3(\tau,0)}\Big)^2+\Big(\frac{\vartheta_4(\tau,z)}{\vartheta_4(\tau,0)}\Big)^2 \bigg]
\ee
where $\tau\in\mathbb H$, $z\in\mathbb C$, and $\vartheta_i$ are the Jacobi theta functions
\be
 \vartheta_{(\alpha,\beta)} (\tau,z)= \sum_{n=-\infty}^{\infty} e^{i\pi\tau(n+\alpha)^2 + 2 i \pi (n+\alpha)(z+\beta)}
\label{sumform}
\ee
with the conventions
\be\label{theta1234}
\begin{array}{rcrrcl}
\vartheta_1 (\tau,z)&:=&{-}\vartheta_{(\frac{1}{2},\frac{1}{2})} (\tau,z)\ ,\qquad
&\vartheta_2 (\tau,z)&:=&\vartheta_{(\frac{1}{2},0)}(\tau,z)\ ,\\
\vartheta_3 (\tau,z)&:=&\vartheta_{(0,0)}(\tau,z)\ ,\qquad
&\vartheta_4 (\tau,z)&:=&\vartheta_{(0,\frac{1}{2})} (\tau,z)\ ,
\end{array}
\ee
and $\vartheta_i(\tau):=\vartheta_i(\tau,0)$. In terms of the quantities $q:=e^{2\pi i\tau}$ and $y:=e^{2\pi i z}$, 
the elliptic genus \eqref{ellk3} is expanded as
\be\label{ellk3ser}
\mathcal E_{\rm K3}(\tau,z)
=20 + 2y^{-1} + 2 y + q\bigg({20}y^{-2} - {128}y^{-1}+ 216 - 128 y + 20 y^2\bigg) + O(q^2)\ .
\ee
The elliptic genus of a K3 $\sigma$-model is defined as \cite{Witten:1986bf,Eguchi:1988vra}
\be\label{ellcft}
\mathcal E(\tau,z)
={\rm tr}_{{\rm R\widetilde R}}\Big(q^{L_0-\frac14}y^{J_0}(-1)^F\ol q^{\widetilde L_0-\frac14}(-1)^{\widetilde F}\Big)
\ee
where $L_0$ and $J_0$ give respectively the holomorphic conformal dimension 
and $R$-charge of the states (the latter is under the $\mathfrak{su}(2)$ $R$-symmetry algebra) 
and $F$ is the fermion number operator. The tilde denotes the anti-holomorphic sector 
and the subscript R refers to the Ramond sector. We next compute the elliptic 
genus \eqref{ellcft} for the $T/\Z_3$ $\sigma$-model. We verify that the 
result agrees with eq. \eqref{ellk3}, i.e.\ the elliptic genus of K3 surfaces, as expected. This shows that the constructed $T/\Z_3$ conformal field theory is indeed a K3 $\sigma$-model.

Prior to the action of the orbifold, we have a $\sigma$-model on $T$, which is a theory of four real free bosons and four real free fermions. The elliptic genus of this conformal field theory is given by
$$
\mathcal E_{T}(\tau,z)
=\bigg(\frac{\vartheta_1(\tau,z)}{\eta(\tau)^3}\bigg)^2 \bigg(\frac{\vartheta_1(\ol\tau,\ol z)}{\eta(\ol\tau)^3}\bigg)^2Z_{\Gamma}(\tau)\Big|_{\ol z=0}
$$
where $Z_{\Gamma}$ is the contribution from the momentum and winding modes. 
$\mathcal E_{T}$ vanishes because $\vartheta_1(\ol\tau,\ol z)|_{\ol z=0}=0$.

The action of the generator $g$ of the $\Z_3$ orbifold group is defined by multiplying one pair of the 
bosonic fields by an overall phase $\xi=e^{\frac{2\pi i}3}$ and the other 
pair by an overall phase $\xi^{-1}=e^{-\frac{2\pi i}3}$ --- see eq. \eqref{Z3action}. 
The two pairs of fermionic fields
are treated similarly. 
Given this action of $\Z_3$ on the free bosons and fermions of the toroidal theory, one may deduce
a consistent $\Z_3$-action on all fields of the sigma model which preserves the $\mathcal N=(4,4)$ symmetry
as well as the spectral flow operators. The order $3$ symmetries in $\Z_3$ are of type 3A in Volpato's 
classification \cite{vo14} of $\mathcal N=(4,4)$-preserving symmetry groups of non-linear sigma models on two-tori.
The contribution to the elliptic genus from the
bosonic oscillator modes transforms as
\be\label{oz3b}
\frac1{\eta(\tau)^4}\;\frac1{\eta(\bar\tau)^4}\quad\longrightarrow\quad
9\,\frac{\eta(\tau)^2}{\vartheta_1\big(\tau,\frac{1}3\big)\;\vartheta_1\big(\tau,-\frac{1}3\big)}\;
\frac{\eta(\bar\tau)^2}{\vartheta_1\big(\bar\tau,\frac{1}3\big)\;\vartheta_1\big(\ol\tau,-\frac{1}3\big)}\, ,
\ee
and that from the fermionic oscillator modes transforms as
\be\label{oz3f}
\frac{\vartheta_1(\tau,z)^2}{\eta(\tau)^2}\;\frac{\vartheta_1(\bar\tau,\bar z)^2}{\eta(\bar\tau)^2}\;\longrightarrow\;\frac{\vartheta_1\big(\tau,z+\frac{1}3\big)\;\vartheta_1\big(\tau,z-\frac{1}3\big)}{\eta(\tau)^2}\;
\frac{\vartheta_1\big(\bar\tau,\bar z+\frac{1}3\big)\;\vartheta_1\big(\bar\tau,\bar z-\frac{1}3\big)}{\eta(\bar\tau)^2}\ .
\ee
We refer the interested reader to references \cite{Wendland:2000ye,Font:2005td} for more details on constructions of toroidal orbifold conformal field theories.
We choose the convention where $(g^\ell,g^m)$, $\ell,m\in\{0,1,2\}$, denotes the $m^{\rm th}$ projected sector of the $\ell^{\rm th}$ twisted sector and $\ell=0$ corresponds to the untwisted sector. Using eqs. \eqref{oz3b} and \eqref{oz3f}, the elliptic genus of the $(1,g^m)$ untwisted sector of the $T/\Z_3$ theory is found to be
\be\label{ellz3unt}
\mathcal E_{(1,g^m)}(\tau,z)
=9\,\frac{\vartheta_1(\tau,z+\frac{m}3)\;\vartheta_1(\tau,z-\frac{m}3)}
{\vartheta_1(\tau,\frac{m}3)\;\vartheta_1(\tau,-\frac{m}3)}\ ,\qquad m\in\{1,2\}\ .
\ee

We next consider the twisted sector of the theory. The elliptic genus in the twisted sector can be obtained as follows. First perform a modular $S$ transformation on the elliptic genus of the untwisted sector, which maps the $(1,g^{\ell})$ sector to the $(g^{\ell},1)$ sector:
\be\label{ellz3gl1}
\mathcal E_{(g^\ell,1)}(\tau,z)
=9\,\frac{\vartheta_1(\tau,z+\frac{\ell\tau}3)\;\vartheta_1(\tau,z-\frac{\ell\tau}3)}
{\vartheta_1\,(\tau,\frac{\ell\tau}3)\;\vartheta_1(\tau,-\frac{\ell\tau}3)}\ .
\ee
The coefficient $9$ agrees with  
the number of singular points of the orbifold $T/\Z_3$. This is then followed by modular $T$ transformations which map to the projected sectors $(g^{\ell},g^m)$. 
Since the order of $\Z_3$ is a prime number, applying repeated $T$ transformations in each 
twisted sector maps $(g^\ell,1)$ to $(g^\ell,g^m)$, $m\in\{1,2\}$, in that sector. We find
\be\label{ellz3glm}
\mathcal E_{(g^\ell,g^m)}(\tau,z)
=9\,\frac{\vartheta_1(\tau,z+\frac{\ell\tau}3+\frac m3)\;\vartheta_1(\tau,z-\frac{\ell\tau}3-\frac m3)}
{\vartheta_1(\tau,\frac{\ell\tau}3+\frac m3)\;\vartheta_1(\tau,-\frac{\ell\tau}3-\frac m3)}\ .
\ee

All in all, the elliptic genus of the $T/\Z_3$ orbifold is given by
\be\label{ellt4z3}
\mathcal E_{T/\Z_3}(\tau,z)
={\textstyle\frac13}\!\!\!\sum_{\;\;\ell,m\in\Z/3\Z}\!\!\!\mathcal E_{(g^\ell,g^m)}(\tau,z)
=3\!\!\!\!\!\sum_{(\ell,m)\ne(0,0)}\!\!\!\!\frac{\vartheta_1(\tau,z+\frac{\ell\tau}3+\frac m3)\;\vartheta_1(\tau,z-\frac{\ell\tau}3-\frac m3)}
{\vartheta_1(\tau,\frac{\ell\tau}3+\frac m3)\;\vartheta_1(\tau,-\frac{\ell\tau}3-\frac m3)}\ .
\ee
This expression has the same $q$-expansion as 
\eqref{ellk3ser}. In fact, the elliptic genus of K3 is a weak Jacobi form of weight 0 and index 1, 
of which up to scaling there is a unique one (referred to as $\phi_{0,1}$). Hence, in checking that \eqref{ellt4z3} is a weak Jacobi form of weight $0$ and index $1$ and comparing the Fourier expansion of \eqref{ellt4z3} with \eqref{ellk3ser}, it is sufficient to match one coefficient to show that the two quantities are equal. An identity relating Jacobi theta functions is obtained from the equality of eqs. \eqref{ellt4z3} and \eqref{ellk3}.

\section{Lattices}\label{app:glue}%
Lattices play an important role in our work. 
We first recall some standard 
lattice terminology in appendix \ref{subapp:background},
and then  in appendix \ref{subapp:glue} we
proceed with a brief description of a lattice gluing 
technique due to Nikulin that we use repeatedly. 
In appendix \ref{subapp:niemeier}, we discuss the Niemeier
lattice of type $A_2^{12}$, which is a particular lattice that we use to
describe symmetries of $\Z_3$-orbifold K3 surfaces, including its 
automorphism group and the role that the Mathieu group $M_{12}$ plays
for these lattice automorphisms. Finally, appendix \ref{subapp:m24}
briefly recalls the Niemeier lattice of type $A_1^{24}$ along with 
the role that the Mathieu group $M_{24}$ plays
for the  automorphisms of this lattice.

\subsection{General background}\label{subapp:background}
Consider a $d$-dimensional real vector space $V$ with  scalar product 
$\langle\cdot,\cdot\rangle$ of  signature $(\gamma_+,\gamma_-)$.
If not stated otherwise, we equip $V=\R^d$ with the standard Euclidean
scalar product of signature $(d,0).$\\

A \textsc{lattice} $\Gamma$ in $V$ is a free $\Z$-module $\Gamma\subset V$ together
with the symmetric bilinear form induced by $\langle\cdot,\cdot\rangle$
on $\Gamma$. Its \textsc{discriminant} $\disc(\Gamma)$  
is the determinant of the matrix of the bilinear form with respect to any choice of basis of $\Gamma$.  
The lattice $\Gamma$ is \textsc{non-degenerate} if $\disc(\Gamma)\neq0$, and it is 
\textsc{unimodular} if $|\!\disc(\Gamma)|=1$.  
By $\Gamma(N)$ with $N$ an integer we denote the same $\Z$-module as $\Gamma$, 
but with bilinear form rescaled by a factor $N$.

By this definition, we can view every sublattice $\Lambda\subset\Gamma$ of $\Gamma\subset V$ 
as a lattice in the vector space $\Lambda\otimes_{\Z}\R\subset V$ with the induced bilinear form.
The \textsc{dual} of a lattice $\Gamma$ in $V$ is $\Gamma^\ast:=\Hom(\Gamma,\Z)$.
If $V=\Gamma\otimes_\Z\R$, then $\Gamma^\ast$ is naturally identified with 
$\left\{ v\in V\mid \langle v,\gamma\rangle\in\Z\; \forall\gamma\in\Gamma\right\}$
by means of the map $v\mapsto\langle v,\cdot\rangle$. With the induced bilinear form,
we therefore may view $\Gamma^\ast$ as a lattice in $V$.
\bigskip

\begin{figure}
\begin{center}
\includegraphics[width=12cm,keepaspectratio]{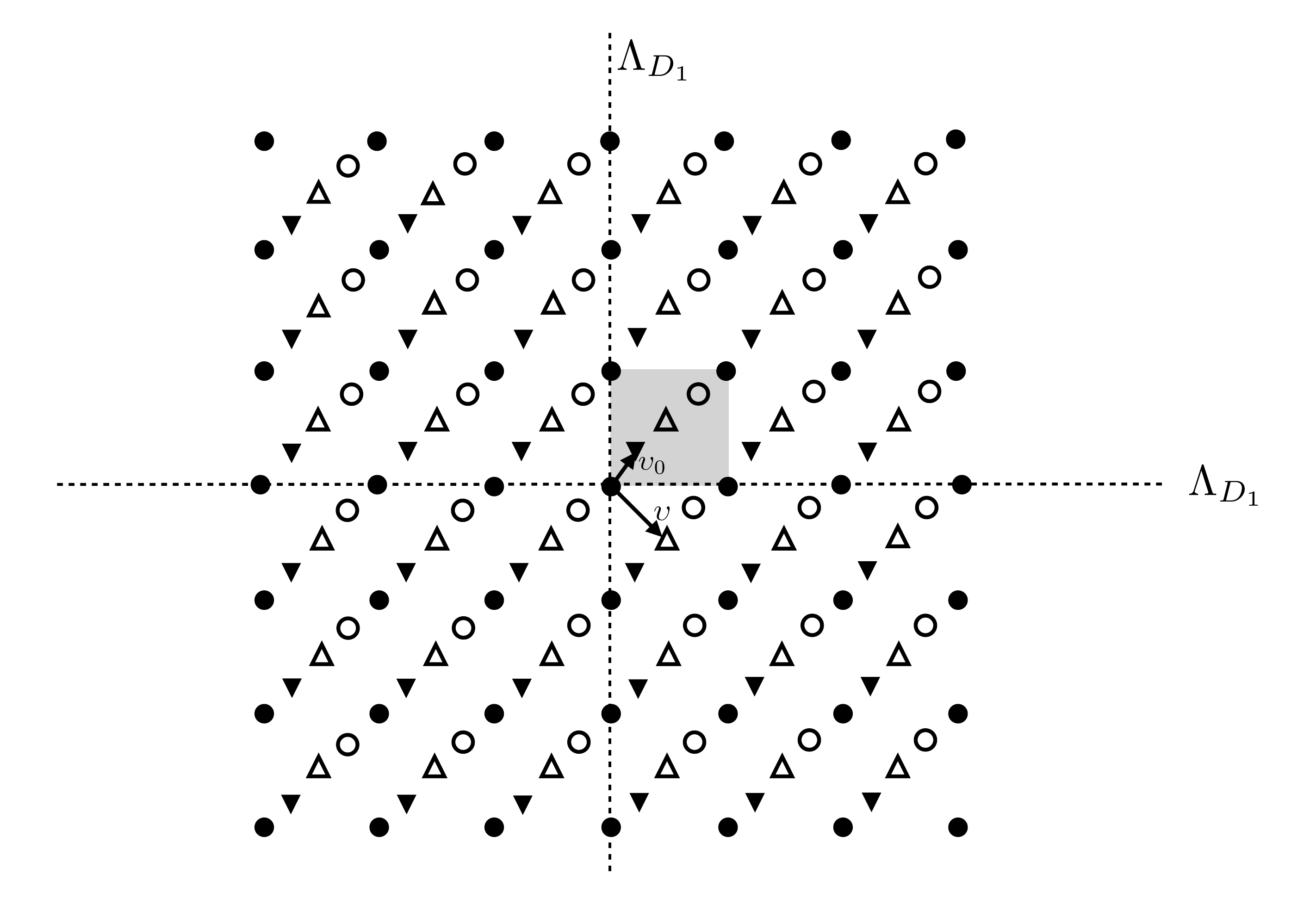}\,\, 
\end{center}
\caption{\small  A small window on the even Lorentzian self-dual lattice $\Gamma^{1,1}$ constructed in terms of two orthogonal copies of the one-dimensional lattice $\Lambda_{D_1}$. The lattice points represented by dots, empty and filled triangles, and circles correspond to vectors $(v_1,v_2)\in ([i],[i])$ for $ [i]=[0], [i]=[v], [i]=[s]$ and $[i]=[c]$ respectively. The vectors $\upsilon_0$ and $\upsilon$ generate $\Gamma^{1,1}$.\,
}
\label{fig:Gamma11}
\end{figure}

A lattice $\Gamma$ in $V$ is \textsc{integral} if its bilinear form has 
values in $\Z$ only. 
It is \textsc{even} if the associated quadratic form has values in $2\Z$ only.  An integral, unimodular lattice is \textsc{self-dual}.
An example of self-dual lattice is provided by the 2-dimensional even Lorentzian lattice $\Gamma^{1,1}$ (see figure \ref{fig:Gamma11}), 
which may be defined in terms of the one-dimensional lattice 
$\Lambda_{D_1}=[0]:=2\Z$, whose dual in $\R$ is
$\Lambda_{D_1}^\ast=\{\hf(4n+m),\, n\in \Z,\, m\in \{-1,0,1,2\}\,\} $. Let 
$[v]:=1+\Lambda_{D_1}, [s]:=\hf+\Lambda_{D_1}$ and $[c]:=-\hf+\Lambda_{D_1}$ be the 
non-trivial cosets\footnote{This notation stems from the fact that the semi-simple Lie 
algebras $D_n=\mathfrak{so}(2n), n>1$, have an associated root lattice $\Lambda_{D_n}$ with one vector and two spinor cosets. $\Lambda_{D_1}$ is defined by extending to $n=1$ the definition of the root lattice $\Lambda_{D_n}$. Note however that there are no vectors of norm $\sqrt{2}$ in $\Lambda_{D_1}$, hence no  roots.}  of $\Lambda_{D_1}$ in $\Lambda_{D_1}^\ast$. 
\bigskip

We set
\begin{align}\label{Gamma11latt}
\Gamma^{1,1}&:=\{(v_1,v_2)\in \R^{1,1}\mid v_1, v_2\in [i], i\in \{0,v,s,c\}\,\} \nonumber \\
&\textstyle
=\{(2n+\hf m,2n^\prime+\hf m)\mid n,n^\prime \in \Z, m\in \{-1,0,1,2\}\}.
\end{align}
With respect to the Lorentzian metric $\rm{diag}(1,-1)$, this lattice is  even, hence also integral. 
It is also unimodular. Indeed, 
a fundamental cell in $\Gamma^{1,1}$ is generated by the two vectors 
\[
\textstyle
\upsilon_0=(\hf, \hf)\,\,{\rm and}\,\,\upsilon=(1,-1),
\]
with Gram matrix $g=\left(\begin{matrix}0&1\\1&0\end{matrix}\right)$ so that
$
{\rm Vol}(\Gamma^{1,1})=\sqrt{|{\rm det}g|}=1.
$
Hence, being integral and unimodular, $\Gamma^{1,1}$ is self-dual. It is also known as the 
\textsc{hyperbolic lattice $U$}, and belongs to a large class of even Lorentzian self-dual lattices called Englert-Neveu lattices in the physics literature \cite{en85,Lerche1989lattices}.\bigskip

If $\Gamma$ is a non-degenerate even lattice in $V$, then there is a natural embedding $\Gamma\hookrightarrow\Gamma^\ast$ 
of finite index by means 
of the bilinear form on $\Gamma$. Thus there is an induced $\Q$-valued
symmetric bilinear form on $\Gamma^\ast$. The 
\textsc{discriminant form}  of $\Gamma$
is the map $q_\Gamma: \Gamma^\ast/\Gamma\longrightarrow\Q/2\Z$ together 
with the symmetric bilinear form $b_\Gamma$ on the \textsc{discriminant group}
$\Gamma^\ast/\Gamma$ with values in $\Q/\Z$, 
induced by the quadratic form and symmetric bilinear form on $\Gamma^\ast$, respectively.
\bigskip

A sublattice $\Lambda\subset\Gamma$ of a lattice $\Gamma$ 
in $V$ is a \textsc{primitive sublattice}
if $\Gamma/\Lambda$ is free.\bigskip

Note that
$\Lambda$ is a primitive sublattice of $\Gamma$ if and only if
$\Lambda=\left(\Lambda\otimes_{\Z}\Q\right)\cap\Gamma$. 
Moreover,  for every non-degenerate integral lattice $\Gamma$ we have
$|\disc(\Gamma)| = |\Gamma^\ast/ \Gamma|$, 
since by definition, $|\disc(\Gamma)|$ is the square of the volume
of any fundamental cell of $\Gamma$ which in turn equals the squared inverse of the
volume of any fundamental cell of $\Gamma^\ast$. Hence  such
$\Gamma$ is unimodular if and only if $\Gamma=\Gamma^\ast$. 
Moreover, it follows that for any sublattice $\Lambda$ of a lattice $\Gamma$,
if $\rk(\Lambda)=\rk(\Gamma)$, then 
\be\label{indexdisc}
|\disc(\Lambda)| = |\disc(\Gamma)| \cdot |\Gamma/\Lambda|^2\; .
\ee

\subsection{Lattice gluing technique}
\label{subapp:glue}
In situations where one has control over a primitive sublattice $\Lambda\subset\Gamma$ of an even unimodular lattice $\Gamma$ and wishes to deduce properties of the lattice $\Gamma$ from those of $\Lambda$,  
 \textsc{gluing techniques} as well as \textsc{criteria for primitivity of sublattices} that were developed by Nikulin in \cite{ni80b,ni80} prove tremendously useful, 
see also \cite{mo84} and \cite[Chapter 4, Section 3]{conway1998sphere}. 
We now recall these techniques and provide a simple example of application.\bigskip

If $\Gamma$ is an even unimodular lattice and $\Lambda\subset\Gamma$ 
is a primitive sublattice, 
then according to \cite[Prop.~1.1]{ni80b} the discriminant forms $q_\Lambda$ and $q_{\m V}$ 
of $\Lambda$ and its orthogonal complement $\m V:=\Lambda^\perp\cap\Gamma$ obey
$q_\Lambda=-q_{\m V}$. Moreover,
let $\Gamma$ denote an even unimodular lattice and $\Lambda\subset\Gamma$ 
a non-degenerate primitive sublattice. Then 
the embedding $\Lambda\hookrightarrow\Gamma$ is uniquely determined by an isomorphism 
$\gamma:\Lambda^\ast/\Lambda\longrightarrow {\m V}^\ast/{\m V}$, where
$\m V=\Lambda^\perp\cap\Gamma$ and 
the discriminant forms  
obey $q_\Lambda=-q_{\m V}\circ\gamma$. Moreover,
\be\label{glueconstruction}
\Gamma =\left\{ \,(\lambda,v)\in \Lambda^\ast\oplus {\m V}^\ast\mid
\gamma([\lambda])=[v]\, \right\},
\ee
where for $L$ a non-degenerate even lattice,
$[l]$ denotes the projection of $l\in L^\ast$  to $L^\ast/L$. \bigskip

According to \cite[\S4.3]{eb02}, the \textsc{glue vectors} 
for the lattices $\Lambda$ 
and $\m V$ are a set of 
vectors $(\lambda, v)\in \left( \Lambda^\ast\oplus\m V^\ast\right)\cap\Gamma$
which together with $\Lambda\oplus\m V$ generate $\Gamma$. For 
a glue vector $(\lambda, v)$, we may replace $\lambda\in\Lambda^\ast$ by any 
other representative in the coset $[\lambda]$, and analogously for $v$.
 \bigskip

Note that \eqref{glueconstruction} allows us to describe $\Gamma$ entirely by means of its 
sublattices $\Lambda$ and $\Lambda^\perp\cap\Gamma$ along with the isomorphism 
$\gamma$.\bigskip

We illustrate the lattice gluing technique by reconstructing the even self-dual lattice $\Gamma^{1,1}$ \eqref{Gamma11latt}
from the sublattice
\[
\textstyle
\Lambda:=\left\{n(\upsilon_0+\upsilon), n\in \Z\right\}
=\left\{n(\frac{3}{2},-\frac{1}{2}), n\in \Z\right\},
\]
which is primitive in $\Gamma^{1,1}$. Indeed,
\[
\Lambda\otimes_{\Z} \Q=\{r(\upsilon_0+\upsilon), r\in \Q\}
\]
and $(\Lambda\otimes_{\Z} \Q)\cap \Gamma^{1,1}=\Lambda$. 
The dual lattice  of $\Lambda$ in $\Lambda\otimes_\Z\R$
is $\Lambda^\ast=\{\frac{n}{2}(\upsilon_0+\upsilon) \mid n\in \Z\}$, so 
that the discriminant group of $\Lambda$ is the set $\Lambda^\ast/\Lambda=\{ \,[0]_+, [\hf]_+\}$ with
$[0]_+=\Lambda$ and $[\hf]_+=\hf (\upsilon_0+\upsilon)+\Lambda$ together with the addition law.
Its orthogonal complement in $\Gamma^{1,1}$ is  $\m V:=\{n(\upsilon_0-\upsilon)\mid n \in \Z\}$ with dual
$\m V^\ast=\{\frac{n}{2}(\upsilon_0-\upsilon) \mid n\in \Z\,\}$ in $\m V\otimes_\Z\R$
and discriminant group $\m V^\ast/\m V=\{\,[0]_-, [\hf]_-\,\}$, where 
$[0]_-=\m V$ and $[\hf]_-=\hf (\upsilon_0-\upsilon)+\m V.$
\begin{figure}[ht]
\begin{center}
\includegraphics[width=10cm,keepaspectratio]{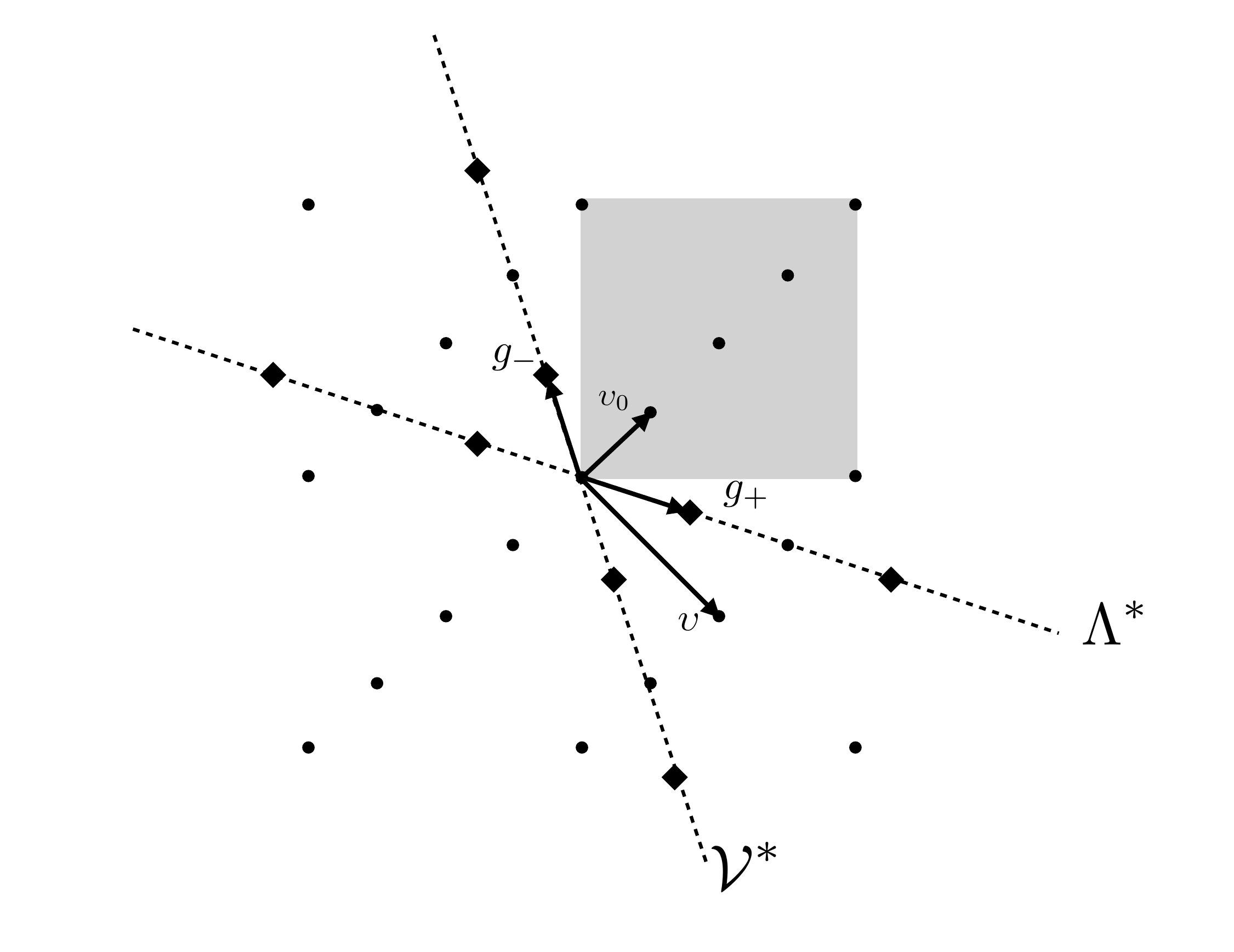}
\end{center}
\caption{\small  Reconstruction of $\Gamma^{1,1}$ from the one-dimensional 
lattice $\Lambda$ primitively embedded in $\Gamma^{1,1}$ and its orthogonal 
complement $\m V$ in $\Gamma^{1,1}$. The chosen glue vector 
is $g_++g_-$, where $g_+=\hf(\upsilon_0+\upsilon)\in\Lambda^\ast$, 
and $g_-=\hf(\upsilon_0-\upsilon)\in\m V^\ast$. All points of the 
lattice $\Gamma^{1,1}$  (represented as dots) 
are obtained by taking integer linear combinations of the 
 form $n_+g_++n_-g_-$, $n_+\equiv n_-$ (mod 2). }
\label{fig:glue}
\end{figure}

The discriminant forms of interest are
\begin{align*}
    &q_{\Lambda}: \Lambda^\ast/\Lambda \longrightarrow \Q/2\Z,\qquad 
    x+\Lambda\mapsto q_{\Lambda}([x])=\langle x,x \rangle\,\,{\rm mod}\,2\Z;\\
    &q_{\m V}: \m V^\ast/\m V \longrightarrow \Q/2\Z,\qquad 
    x+\m V\mapsto q_{\m V}([x])=\langle x,x\rangle\,\,{\rm mod}\,2\Z,
\end{align*}
where $\langle \cdot,\cdot\rangle$ is  given by the 
 induced bilinear form 
on $\Lambda^\ast$ and $\m V^\ast$ respectively. 
For instance, $q_{\Lambda}\left([\hf]_+ \right)=\langle\,\hf(\upsilon_0+\upsilon),\hf(\upsilon_0+\upsilon)\,\rangle
=\frac{1}{4}\cdot 2\langle\, \upsilon_0, \upsilon \, \rangle=\hf$. 
A choice of coset representatives  
for $\Lambda$ in $\Lambda^\ast$ is $\{ (0,0), \hf(\upsilon_0+\upsilon)\}$ 
while for $\m V$ we may take $\{ (0,0), \hf(\upsilon_0-\upsilon)\}$.
Note that the off-diagonal elements of the associated bilinear forms on the discriminant groups take values in
$\Q/\Z$.
\bigskip

The gluing isomorphism $\gamma: \Lambda^\ast/\Lambda\rightarrow \m V^\ast/\m V$ defined by
\[
\textstyle
\gamma(\,[0]_+)=[0]_-\,\,{\rm and}\,\,\gamma([\hf]_+)=[\hf]_-
\]
ensures that $q_{\Lambda}=-q_{\m V}\circ \gamma$. The lattice $\Gamma^{1,1}$ is thus reconstructed as
\[
\Gamma^{1,1}\cong\left\{ \,(\lambda,v)\in \Lambda^\ast\oplus {\m V}^\ast\mid
\gamma([\lambda])=[v]\, \right\}, 
\]
where $\lambda=\frac{n_+}{2}(\upsilon_0+\upsilon)$ and $v=\frac{n_-}{2}(\upsilon_0-\upsilon)$ 
for integers $n_+$ and $n_-$ obey $\gamma([\lambda])=[v]$ if $n_+$ and $n_-$ have the 
same parity (see figure \ref{fig:glue}).

Note that our original definition \eqref{Gamma11latt} of the hyperbolic lattice  can be viewed as yet another
construction of $\Gamma^{1,1}$ by gluing, namely from two $D_1$ root lattices $\Lambda\cong\m V\cong\Lambda_{D_1}$, 
where $\Lambda=\left\{ (2n,0) \mid n\in\Z\right\}$, $\m V=\left\{ (0,2n^\prime) \mid n^\prime\in\Z\right\}$,
$\Lambda^\ast/\Lambda\cong \m V^\ast/\m V\cong\Z_4$.
\bigskip

The above lattice gluing technique yields a convenient criterion
to decide whether or not an automorphism of a primitive 
sublattice $\Lambda\subset\Gamma$ of an even 
unimodular lattice $\Gamma$ can be extended to a lattice
automorphism\footnote{Here, a lattice automorphism of a lattice $L$ is an isometry of $L$
which possesses a linear extension to an automorphism of the ambient vector space.
The automorphism group of $L$ is denoted $\Aut(L)$.} of $\Gamma$. 
As before, denote by $\m V$ the orthogonal complement
$\m V =\Lambda^\perp\cap\Gamma$
of $\Lambda$ in $\Gamma$, and let $\gamma:\Lambda^\ast/\Lambda\longrightarrow {\m V}^\ast/{\m V}$
denote the gluing isomorphism such that \eqref{glueconstruction} holds. 
Moreover, let $\varphi\in\Aut(\Lambda)$, and denote by $\ol\varphi$ the induced 
automorphism of the discriminant group $\Lambda^\ast/\Lambda$. 
We claim
\be\label{autoextend}
\Phi\in\Aut(\Gamma) \mbox{ exists with } \Phi_{|\Lambda}= \varphi
\quad\Longleftrightarrow\quad
\psi\in\Aut(\m V) \mbox{ exists with } \gamma\circ\ol\varphi=\ol\psi\circ\gamma \; 
\ee
where $\ol\psi$ denotes the automorphism of the discriminant group $\m V^\ast/\m V$ 
induced by $\psi\in\Aut(\m V)$. Here and in the following, for ease of notation,
we denote linear extensions of lattice automorphisms to larger lattices of same rank by the 
same letter.

To prove \eqref{autoextend}, recall that by definition, lattice automorphisms
respect the bilinear form, such that 
$\Phi\in\Aut(\Gamma)$ with $\Phi_{|\Lambda}= \varphi$
exists if and only if $\psi\in\Aut(\m V)$ exists with the property
that $\Phi(\lambda,v):=(\varphi(\lambda),\psi(v))$ for all 
$(\lambda,v)\in\Lambda^\ast\oplus\m V^\ast$ yields $\Phi_{|\Gamma}\in\Aut(\Gamma)$.
By \eqref{glueconstruction}, equivalently we have that 
$\psi\in\Aut(\m V)$ exists with the property that 
$$
\forall (\lambda,v)\in\Lambda^\ast\oplus\m V^\ast\colon\qquad
\gamma([\lambda])=[v]
\quad\Longleftrightarrow\quad
\gamma( \underbrace{[\varphi(\lambda)]}_{\ol\varphi( [\lambda])} ) 
= \underbrace{[ \psi(v)]}_{\ol\psi([v])}\; .
$$
Equivalently, 
$\psi\in\Aut(\m V)$ exists with $\gamma\circ\ol\varphi=\ol\psi\circ\gamma$,
and we have proven \eqref{autoextend}.

\subsection{The Niemeier lattice of type \texorpdfstring{$A_2^{12}$}{TEXT} and the Mathieu group \texorpdfstring{$M_{12}$}{TEXT}}
\label{subapp:niemeier}
The groups of lattice automorphisms for indefinite lattices of 
rank at least $3$ are infinite (see e.g. \cite{vi75}), while
for lattices with definite signature, as permutation groups
on every set of vectors of equal length, they are always finite.
Therefore, it can be advantageous to use positive definite lattices
when describing symmetry groups. Here, the Niemeier lattices
are of particular interest.

A \textsc{Niemeier lattice} is an even unimodular positive definite 
lattice of rank $24$. According to Niemeier's Theorem \cite[Satz 8.3]{ni73},
up to isometries, there are $24$ such lattices, and each
of them is uniquely determined by its \textsc{root lattice}, that is, by
the sublattice generated by all lattice vectors of length $\sqrt2$. The 
Niemeier lattices may therefore be
labelled by the Dynkin diagrams of their root lattices of
rank $24$ or, in the case of the Leech lattice, by $0$.

In this note, we are interested in the Niemeier lattice of type $A_2^{12}$, 
which we denote by $N$.
The root lattice of $N$, which we call $\widetilde{R}$,
is given by the direct sum of twelve pairwise 
perpendicular sublattices of type $A_2$, labelled as $A_{2,j}$, $j\in \qty{1, \dots , 12}$. That is,
$$ 
\widetilde{R}:=\bigoplus_{j=1}^{12}A_{2,j}\, . 
$$
 We choose simple roots in $A_{2,j}$, denoted by $\widetilde E^{(l)}_j$, $l\in\qty{1,2}$, such that
\begin{align*}
    \forall \, j,k\in\qty{1,\dots,12}, \, l,m\in\qty{1,2} \qquad \expval{ \widetilde{E}^{(l)}_j \, , \, \widetilde{E}^{(m)}_k  } = \begin{cases}
        \,\,\,\, 0  & \text{if } j\neq k \, ,\\
        -1 & \text{if } j=k, l\neq m \, , \\
        \,\,\,\, 2 & \text{if } j=k, l=m \, .
    \end{cases}
\end{align*}
Moreover, let $\widetilde E_j:=\widetilde E^{(1)}_j+2\widetilde{E}^{(2)}_j$ 
for $j\in\qty{1,\dots, 12}$; then for the sublattice $A_{2,j}$, the dual lattice
$A_{2,j}^\ast$  may be generated by 
the $\widetilde E^{(l)}_j$ with $l\in\qty{1,2}$
and $\frac{1}{3} \widetilde E_j$. One checks that $\left( \wt E_j^{(2)}, \frac{1}{3} \widetilde E_j\right)$
forms a basis of the dual lattice $A_{2,j}^\ast$, and that the shortest non-zero vectors 
in this lattice have squared length $\frac{2}{3}$.

The Niemeier lattice $N$ is obtained from the components $A_{2,j}$ of 
the root lattice $\widetilde{R}$ by including a set of additional generators (also often called 
glue vectors, as in \cite{conway1998sphere,eb02} for instance)
of the form
\begin{align*}
    y
    =\textstyle
    \frac{1}{3}\sum\limits_{j=1}^{12}y_j\, \widetilde E_j\,, 
\end{align*}
where each $y_j$,\, $j\in \qty{1, \dots , 12}$, corresponds to a 
representative $\frac{1}{3}y_j\, \widetilde E_j$
of an element 
of the discriminant group $A_{2,j}^*/A_{2,j}$ of order 3. 
As a convenient shorthand notation we also introduce
$$
\textstyle
[0]:=A_{2,j},\quad
  [+]:=\frac{1}{3}\widetilde E_j + A_{2,j}\quad {\rm and} \quad 
    [-]:=-\frac{1}{3}\widetilde E_j + A_{2,j}.
    $$
By the above, $\widetilde R^\ast/\widetilde R\cong \left( A_2^\ast/A_2\right)^{12}\cong\F_3^{12}$,
and since $\widetilde R\subset N\subset\widetilde R^\ast$, we may view the image of $N$
under the projection to $\widetilde R^\ast/\widetilde R$ as a subgroup of $\F_3^{12}$.
This group is called the \textsc{glue code}, and we denote the image of $y\in N$ in the glue code 
by $[y]=[y_1,\,\ldots,\,y_{12}]$, where  
$y_j$ is the projection of the component of $y$ in $A_{2,j}^\ast$ to $\F_3\cong A_{2,j}^\ast/A_{2,j}$.
According to \cite[\S18.4.II]{conway1998sphere}, the glue code of 
$N$ is the \textsc{extended ternary Golay code} 
$\m{C}_{12}$, 
a $6$-dimensional vector subspace of $\F_3^{12}$ whose elements are 
called codewords.\footnote{The code $\m{C}_{12}$ is often referred to 
as the linear code $[12, 6, 6]_3$, where $12$ is its length, 
i.e. the number of coordinates describing the ambient vector 
space $\F_3^{12}$, the first $6$ is its dimension as a vector space 
and the second $6$ is the minimal non-trivial weight of its codewords 
(the weight of a codeword is the number of its non-zero components). 
The subscript $3$ signifies that one works over the finite field $\F_3$.} 
The vector space $\m C_{12}$ contains $3^6=729$ codewords, which can be 
generated by the images $[w_i]\in \F_3^{12}$ 
of the six generators of $N$ from $\widetilde R$ given by
\begin{equation}\label{gluecodewords}
\textstyle
    w_i = \frac{1}{3} \sum\limits_{j=1}^{12} a_{ij} \widetilde E_j,
\end{equation}
where $a_{ij}\in\qty{0,\pm 1}$ with $i\in\qty{1,\dots, 6}$, $j\in\qty{1,\dots, 12}$ are the entries of the generating matrix $\m G$:
\begin{equation}\label{matrixcode}
\m G=\left[
\begin{array}{cccccc|cccccc}
0 & + & + & + & + & + & + & 0 & 0 & 0 & 0 & 0 \\
- & 0 & + & - & - & + & 0 & + & 0 & 0 & 0 & 0 \\
- & + & 0 & + & - & - & 0 & 0 & + & 0 & 0 & 0 \\
- & - & + & 0 & + & - & 0 & 0 & 0 & + & 0 & 0 \\
- & - & - & + & 0 & + & 0 & 0 & 0 & 0 & + & 0 \\
- & + & - & - & + & 0 & 0 & 0 & 0 & 0 & 0 & + \\
\end{array}
\right] \, .
\end{equation}
The generating matrix $\m G$ stated in \eqref{matrixcode} is obtained from the description 
given in \cite[\S10.1.5]{conway1998sphere} 
as follows. We use the basis vectors of $\m C_{12}$
denoted $\w_\infty$, $\w_1$, $\w_3$, $\w_4$, $\w_5$, $\w_9$ in 
\cite[\S10.1.5]{conway1998sphere}, each of which has $12$ entries from $\F_3$ 
that are labelled by $(\infty,0,1,\ldots,10)$. 
We relabel these entries according to 
\be\label{relabelCSToBTUWZ}
(\infty,0,1,\ldots,10)\quad\longmapsto\quad(8,10,9,12,7,11,3,2,1,6,5,4)\,,
\ee
such that, for example, $\w_1=(+,-,+,-,-,-,+,+,+,-,+,-)$ after relabelling (and reordering in 
numerical order $(1,2,\ldots,12)$)
becomes $(+,+,+,-,+,-,-,+,+,-,-,-)$.
One verifies that the resulting basis vectors of a copy of $\m C_{12}$ are linear 
combinations of the row vectors of the generating matrix $\m G$ given in \eqref{matrixcode}.
These row vectors serve as a convenient basis of our representation of $\m C_{12}$.

Altogether, the Niemeier lattice $N$ can thus be conveniently generated by the set
\begin{align}\label{Ngen}
    \qty{ \widetilde E_j^{(l)} \, , j\in\qty{1,\dots 12} \, , l\in\qty{1,2}  ; \; w_1,\dots w_6 } \, .
\end{align}

Now we move on to describe the group of lattice 
automorphisms of the Niemeier lattice $N$ following \cite[p.100]{conway1998sphere}. 
This group, which we call $\Aut(N)$, is
described in terms of three groups, denoted $G_0$, $G_1$, $G_2$, where
$\Aut(N)\cong \left(G_0\times G_1 \right) \cdot G_2$, that is, $G_0\times G_1\lhd \Aut(N)$
and $\Aut(N)$ is a non-split extension of $G_0\times G_1$ by $G_2$. Here, 
\begin{itemize}
\item 
$G_{0}\cong(S_3)^{12}$ is the direct product of the Weyl groups $S_3$ of each $A_2$ component. 
In other words, $G_0$ is the normal subgroup of $\Aut(N)$
generated by the reflections in the hyperplanes whose normal vectors are simple roots in any one of the $A_2$ components.
\item 
$G_1\cong\Z_2$ is the group whose generator acts by multiplication by $(-1)$ on every lattice vector.
It is the centre of $\Aut(N)$.
\item 
$G_2\cong M_{12}$ is the  group  of  permutations  of the twelve 
components $A_{2,j}$, $j\in\{1,\ldots,12\}$, of the root lattice $\wt R$ induced by
$\Aut(N)$. 
\end{itemize}
According to \cite[\S2.8.5]{conway1998sphere}, 
the automorphism group $\Aut(\m C_{12})$
of the extended ternary Golay code $\m C_{12}$
is isomorphic to $G_1\cdot M_{12}$,
where the non-trivial element of $G_1$ acts by $+\longleftrightarrow-$ in each component of the 
codewords.
By the above, $G_1\cdot M_{12}$ is the group of automorphisms of $N$ modulo the normal subgroup generated by 
reflections in the hyperplanes given by the roots.
It is one of the 23 umbral groups in \cite{cheng2014umbral}.

According to \cite[\S10.1.7]{conway1998sphere}, the Mathieu group $M_{12}\cong \Aut(\m C_{12})/\Z_2$
may be represented faithfully as permutation group which acts on the twelve entries of the code
words of $\m C_{12}\subset\F_3^{12}$, keeping in mind that lifting elements from $M_{12}$ to 
$\Aut(\m C_{12})$ amounts to introducing appropriate sign flips ``$+\leftrightarrow-$''
in some of the components and that this cannot be done by means of a group homomorphism.
Implementing the relabelling obtained in \eqref{relabelCSToBTUWZ}, the generators of 
$M_{12}$ that are introduced in \cite[\S10.1.7]{conway1998sphere} are given by the permutations
\be
\label{generatingM12}
\left.
\begin{array}{rclrcl}
A&=& (1,6,5,4,10,9,12,7,11,3,2)\, ,&
B&=& (1,4,6,12,2)(3,11,9,7,5)\, ,\\
C&=& (1,7)(2,5)(3,12)(4,9)(6,11)(8,10)\, ,&
D&=& (1,2)(3,5)(4,12)(7,11)\, .
\end{array}\right\}
\ee
One checks that $A,\, B,\, D$ lift to automorphisms of $\m C_{12}$ without
any additional sign flips, while $C$ must be accompanied by sign flips
$$
\F_3^{12}\longrightarrow \F_3^{12}\, , \;\;
(x_1,\ldots,x_{12})
\mapsto\pm (x_1,x_2,-x_3,x_4,-x_5,x_6,-x_7,x_8,-x_9,-x_{10},-x_{11},x_{12})\, ,
$$
in accord with \cite[\S10.1.7]{conway1998sphere}.

\subsection{The Niemeier lattice of type \texorpdfstring{$A_1^{24}$}{TEXT} and the Mathieu group \texorpdfstring{$M_{24}$}{TEXT}}
\label{subapp:m24}

In this appendix, we briefly recall a few standard facts about the Niemeier lattice
$\ch{N}$ of type $A_1^{24}$ and its automorphism group. 
For a more detailed discussion, see \cite{conway1998sphere}, or see \cite[Appendix A]{tawe11}
for a review; we continue to use the conventions that two of the authors
have chosen in the latter reference.

The root lattice of $\ch N$ is 
given by the direct sum of $24$ pairwise perpendicular sublattices of type $A_1$,
and we denote the roots by $\pm\!\ch E_1,\ldots, \pm\!\ch E_{24}$ with 
$\langle \ch E_j,\ch E_k\rangle = 2\delta_{jk}$ for $j,k\in\{1,\ldots,24\}$. We have
$(A_1^{24})^\ast/ A_1^{24}\cong \F_2^{24}$, and
the lattice $\ch N$ may by generated from its root lattice along with the twelve 
vectors $\frac{1}{2}\smash{\sum\limits_{k=1}^{24}} a_{jk} \ch E_k$, $j\in\{1,\ldots,12\}$,
where the $a_{jk}\in\{0,1\}$
are the entries of the matrix
$$
\left(
\begin{array}{cccccc|cccccc|cccccc|cccccc}
1 & 1 & 1 & 0 & 1 & 0 & 0 & 0 & 0 & 0 & 0 & 1 & 0 & 0 & 0 & 1 & 0 & 1 & 0 & 0 & 0 & 1 & 0 & 0\\
0 & 0 & 1 & 1 & 0 & 1 & 0 & 0 & 1 & 0 & 0 & 0 & 0 & 0 & 1 & 0 & 1 & 0 & 0 & 0 & 1 & 1 & 0 & 0\\
0 & 0 & 0 & 1 & 0 & 0 & 1 & 1 & 0 & 1 & 0 & 1 & 0 & 0 & 1 & 1 & 0 & 0 & 1 & 0 & 0 & 0 & 0 & 0\\
0 & 0 & 0 & 0 & 0 & 0 & 1 & 1 & 0 & 0 & 1 & 0 & 0 & 0 & 1 & 1 & 1 & 1 & 0 & 0 & 0 & 0 & 1 & 0\\
1 & 1 & 0 & 1 & 0 & 0 & 1 & 1 & 0 & 0 & 0 & 0 & 0 & 0 & 0 & 1 & 0 & 1 & 0 & 0 & 1 & 0 & 0 & 0\\
1 & 0 & 0 & 0 & 1 & 0 & 0 & 1 & 0 & 0 & 0 & 0 & 1 & 1 & 0 & 1 & 0 & 1 & 0 & 0 & 0 & 0 & 1 & 0\\
\hline
0 & 0 & 0 & 0 & 0 & 0 & 0 & 0 & 0 & 1 & 1 & 1 & 1 & 1 & 0 & 0 & 1 & 0 & 0 & 1 & 0 & 1 & 0 & 0\\
0 & 0 & 1 & 0 & 0 & 0 & 0 & 1 & 1 & 0 & 1 & 0 & 1 & 0 & 0 & 0 & 0 & 1 & 1 & 0 & 0 & 0 & 1 & 0\\
0 & 1 & 0 & 0 & 1 & 0 & 0 & 0 & 1 & 1 & 0 & 0 & 0 & 1 & 0 & 0 & 0 & 0 & 1 & 0 & 1 & 0 & 0 & 1\\
1 & 0 & 0 & 1 & 0 & 0 & 1 & 1 & 1 & 0 & 0 & 0 & 1 & 0 & 1 & 1 & 0 & 1 & 1 & 0 & 0 & 1 & 1 & 0\\ 
1 & 1 & 1 & 0 & 0 & 0 & 1 & 1 & 0 & 1 & 1 & 0 & 1 & 0 & 0 & 1 & 0 & 1 & 0 & 0 & 0 & 1 & 1 & 0\\
1 & 1 & 0 & 0 & 1 & 0 & 1 & 1 & 1 & 1 & 0 & 1 & 0 & 0 & 1 & 1 & 0 & 1 & 0 & 0 & 0 & 0 & 1 & 0\\
\end{array}
\right) \, .
$$
Interpreting the entries as elements of the field $\F_2$ with two elements,
this matrix gives the generating matrix of the \textsc{extended binary Golay code} 
$\m C_{24}\subset\F_2^{24}$.

The automorphism group of the Niemeier lattice $\ch N$ is $(\Z_2)^{24}\rtimes M_{24}$, where
$(\Z_2)^{24}$ is the group generated by the reflections in the hyperplanes whose
normal vectors are the
roots, while $M_{24}$ is the 
largest Mathieu group, which acts by permutations on the $24$ components of type $A_1$ of the
root lattice and thereby on the entries
of the codewords in $\m C_{24}$. 
Thus $M_{24}$ is also the automorphism group of the extended binary 
Golay code $\m C_{24}$. We therefore 
find it convenient to use the extended binary 
Golay code $\m C_{24}$ to illustrate the group $M_{24}$, along with
the very playful visualization of codewords in $\m C_{24}$ that is offered by the 
Miracle Octad Generator (MOG). The MOG was devised by Robert Curtis \cite{cu74},
and we use a variant which was developed by Conway shortly after. 
It is given by a $4 \times 6$ array whose entries are elements of 
$\mathbb{F}_2=\{0, 1\}$, and therefore provides binary words of length $24$. 
Following the convention introduced by Conway, we populate the entries of the 
$4 \times 6$ array with the $24$ coefficients of a vector in $\F_2^{24}$, labelled
by $\{1,\ldots,24\}$ according to

\be\label{ConwayMOG}
\begin{array}{|cc|cc|cc|}
\hline
23&24&1&11&2&22\\
19&3&20&4&10&18\\
\hline
15&6&14&16&17&8\\
5&9&21&13&7&12\\
\hline \end{array}
\ee
(cf.\ \cite[\S11.8, Figure 11.7]{conway1998sphere}).
For better readability, we omit the entries $0$ and replace the entries $1$ by 
the symbol $\bullet$. Thus the codeword given by the first row of the above
generating matrix is depicted by the MOG
\be\label{sampleword}%
\begin{array}{|cc|cc|cc|}
\hline
\vphantom{\bullet}&\hphantom{\bullet}&\bullet&\hphantom{\bullet}&\bullet&\bullet\\
\hphantom{\bullet}&\bullet&\hphantom{\bullet}&\hphantom{\bullet}&\hphantom{\bullet}&\bullet\\
\hline
\hphantom{\bullet}&\hphantom{\bullet}&\hphantom{\bullet}&\bullet&\hphantom{\bullet}&\hphantom{\bullet}\\
\bullet&\hphantom{\bullet}&\hphantom{\bullet}&\hphantom{\bullet}&\hphantom{\bullet}&\bullet\\
\hline \end{array}
\ee
The \textsc{hexacode} ${\m H}_6$ is used to determine whether or not a vector in 
$\F_2^{24}$, depicted in terms of the above array, belongs to $\m C_{24}$.
The hexacode is a $3$-dimensional code of length $6$ over the field of four 
element $\mathbb{F}_4=\{0, 1, \omega, \omega^2\}$, with $\omega^3=1$, $1+\omega=\omega^2$,
$1+1=\omega+\omega=\omega^2+\omega^2=0$.
It may be defined as
 $$
 {\m H}_6=\{ (a, b, \phi(0), \phi(1), \phi(\omega), \phi ({\omega^2})) | 
 a, b, \phi(0) \in \mathbb{F}_4, \, \phi(x):=ax^2+bx+\phi(0)\}.
 $$
Golay codewords are those configurations which obey the following two 
conditions:
\begin{enumerate}
\item 
The parity of each $4$-column of the 
MOG and the parity of the top row (the parity of a column or a row being the parity of the 
sum of its entries) are all equal.
\item  
To each 4-column with entries $\alpha, \beta, \gamma, \delta \in \mathbb{F}_2$, 
associate the $\mathbb{F}_4$ element $\beta +\gamma \omega +\delta {\omega^2}$, which 
is called its \textsc{score}.
The six scores calculated from a given MOG form a 
hexacode word. 
\end{enumerate}
One readily checks that the above example \eqref{sampleword}  depicts a
codeword in $\m C_{24}$: indeed, all columns have odd parity, as does the first row.
Moreover, the scores obtained from the 4-columns are 
$(\omega^2,1,0,\omega,0,\omega) = (\omega^2,1,0, \phi(1), \phi(\omega), \phi ({\omega^2}))$
with $\phi(x):=(\omega x)^2+x$, confirming that the scores yield a hexacode word.

The binary extended Golay code $\m C_{24}$ contains exactly one codeword 
zero and one codeword where all digits are one, together 
with $759$ octads (i.e.\ codewords with precisely eight entries ``1''), 
$2576$ dodecads (i.e.\ codewords with precisely twelve entries ``1''), 
and $759$ complement octads (i.e.\ codewords with precisely sixteen entries ``1'').
The following is a non-exhaustive list of octads in terms of their MOGs; 
we use the number that is attached to the South-East
corner of each MOG for reference in the main text:
\begin{eqnarray}\label{domino}
&&\hspace*{-1.5em}
\resizebox{0.95\hsize}{!}{$
\begin{array}{|cc|cc|cc|}
\hline
\hphantom{\bullet}&\hphantom{\bullet}&\hphantom{\bullet}&\hphantom{\bullet}&{\bullet}&{\bullet}\\
\hphantom{\bullet}&\hphantom{\bullet}&\hphantom{\bullet}&\hphantom{\bullet}&{\bullet}&{\bullet}\\
\hline
\hphantom{\bullet}&\hphantom{\bullet}&\hphantom{\bullet}&\hphantom{\bullet}&{\bullet}&{\bullet}\\
\hphantom{\bullet}&\hphantom{\bullet}&\hphantom{\bullet}&\hphantom{\bullet}&{\bullet}&{\bullet}\\
\hline \end{array}_{\domino}
\quad
\begin{array}{|cc|cc|cc|}
\hline
\hphantom{\bullet}&\hphantom{\bullet}&\hphantom{\bullet}&{\bullet}&\hphantom{\bullet}&{\bullet}\\
\hphantom{\bullet}&\hphantom{\bullet}&\hphantom{\bullet}&{\bullet}&\hphantom{\bullet}&{\bullet}\\
\hline
\hphantom{\bullet}&\hphantom{\bullet}&\hphantom{\bullet}&{\bullet}&\hphantom{\bullet}&{\bullet}\\
\hphantom{\bullet}&\hphantom{\bullet}&\hphantom{\bullet}&{\bullet}&\hphantom{\bullet}&{\bullet}\\
\hline \end{array}_{\domino}
\quad
\begin{array}{|cc|cc|cc|}
\hline
\hphantom{\bullet}&\hphantom{\bullet}&{\bullet}&\hphantom{\bullet}&{\bullet}&\hphantom{\bullet}\\
\hphantom{\bullet}&\hphantom{\bullet}&\hphantom{\bullet}&{\bullet}&\hphantom{\bullet}&{\bullet}\\
\hline
\hphantom{\bullet}&\hphantom{\bullet}&\hphantom{\bullet}&{\bullet}&\hphantom{\bullet}&{\bullet}\\
\hphantom{\bullet}&\hphantom{\bullet}&{\bullet}&\hphantom{\bullet}&{\bullet}&\hphantom{\bullet}\\
\hline \end{array}_{\domino}
\quad
\begin{array}{|cc|cc|cc|}
\hline
\hphantom{\bullet}&\hphantom{\bullet}&{\bullet}&\hphantom{\bullet}&{\bullet}&\hphantom{\bullet}\\
\hphantom{\bullet}&\hphantom{\bullet}&\hphantom{\bullet}&{\bullet}&\hphantom{\bullet}&{\bullet}\\
\hline
\hphantom{\bullet}&\hphantom{\bullet}&{\bullet}&\hphantom{\bullet}&{\bullet}&\hphantom{\bullet}\\
\hphantom{\bullet}&\hphantom{\bullet}&\hphantom{\bullet}&{\bullet}&\hphantom{\bullet}&{\bullet}\\
\hline \end{array}_{\domino}
\quad
\begin{array}{|cc|cc|cc|}
\hline
\hphantom{\bullet}&\hphantom{\bullet}&\hphantom{\bullet}&{\bullet}&{\bullet}&\hphantom{\bullet}\\
\hphantom{\bullet}&\hphantom{\bullet}&\hphantom{\bullet}&{\bullet}&\hphantom{\bullet}&{\bullet}\\
\hline
\hphantom{\bullet}&{\bullet}&\hphantom{\bullet}&\hphantom{\bullet}&{\bullet}&{\bullet}\\
\hphantom{\bullet}&{\bullet}&\hphantom{\bullet}&\hphantom{\bullet}&\hphantom{\bullet}&\hphantom{\bullet}\\
\hline \end{array}_{\domino}$}
\nonumber\\[10pt]
&&\hspace*{-1.5em}
\resizebox{0.95\hsize}{!}{$
\begin{array}{|cc|cc|cc|}
\hline
\hphantom{\bullet}&\hphantom{\bullet}&{\bullet}&\hphantom{\bullet}&{\bullet}&\hphantom{\bullet}\\
\hphantom{\bullet}&{\bullet}&\hphantom{\bullet}&\hphantom{\bullet}&\hphantom{\bullet}&\hphantom{\bullet}\\
\hline
\hphantom{\bullet}&{\bullet}&\hphantom{\bullet}&\hphantom{\bullet}&{\bullet}&{\bullet}\\
\hphantom{\bullet}&\hphantom{\bullet}&{\bullet}&\hphantom{\bullet}&\hphantom{\bullet}&{\bullet}\\
\hline \end{array}_{\domino}
\quad
\begin{array}{|cc|cc|cc|}
\hline
\hphantom{\bullet}&\hphantom{\bullet}&{\bullet}&{\bullet}&\hphantom{\bullet}&\hphantom{\bullet}\\
\hphantom{\bullet}&{\bullet}&{\bullet}&\hphantom{\bullet}&\hphantom{\bullet}&{\bullet}\\
\hline
\hphantom{\bullet}&{\bullet}&\hphantom{\bullet}&{\bullet}&\hphantom{\bullet}&{\bullet}\\
\hphantom{\bullet}&\hphantom{\bullet}&\hphantom{\bullet}&\hphantom{\bullet}&\hphantom{\bullet}&\hphantom{\bullet}\\
\hline \end{array}_{\domino}
\quad
\begin{array}{|cc|cc|cc|}
\hline
\hphantom{\bullet}&{\bullet}&{\bullet}&\hphantom{\bullet}&\hphantom{\bullet}&\hphantom{\bullet}\\
\hphantom{\bullet}&{\bullet}&\hphantom{\bullet}&{\bullet}&\hphantom{\bullet}&\hphantom{\bullet}\\
\hline
\hphantom{\bullet}&\hphantom{\bullet}&{\bullet}&{\bullet}&\hphantom{\bullet}&{\bullet}\\
\hphantom{\bullet}&\hphantom{\bullet}&\hphantom{\bullet}&\hphantom{\bullet}&\hphantom{\bullet}&{\bullet}\\
\hline \end{array}_{\domino}
\quad
\begin{array}{|cc|cc|cc|}
\hline
\hphantom{\bullet}&{\bullet}&\hphantom{\bullet}&\hphantom{\bullet}&{\bullet}&\hphantom{\bullet}\\
\hphantom{\bullet}&{\bullet}&\hphantom{\bullet}&\hphantom{\bullet}&{\bullet}&\hphantom{\bullet}\\
\hline
\hphantom{\bullet}&{\bullet}&\hphantom{\bullet}&\hphantom{\bullet}&{\bullet}&\hphantom{\bullet}\\
\hphantom{\bullet}&{\bullet}&\hphantom{\bullet}&\hphantom{\bullet}&{\bullet}&\hphantom{\bullet}\\
\hline \end{array}_{\domino}
\quad
\begin{array}{|cc|cc|cc|}
\hline
\hphantom{\bullet}&{\bullet}&\hphantom{\bullet}&{\bullet}&\hphantom{\bullet}&\hphantom{\bullet}\\
\hphantom{\bullet}&{\bullet}&\hphantom{\bullet}&{\bullet}&\hphantom{\bullet}&\hphantom{\bullet}\\
\hline
\hphantom{\bullet}&{\bullet}&\hphantom{\bullet}&{\bullet}&\hphantom{\bullet}&\hphantom{\bullet}\\
\hphantom{\bullet}&{\bullet}&\hphantom{\bullet}&{\bullet}&\hphantom{\bullet}&\hphantom{\bullet}\\
\hline \end{array}_{\domino}
$}\nonumber\\[10pt]
&&\hspace*{-1.5em}
\resizebox{0.95\hsize}{!}{$
\begin{array}{|cc|cc|cc|}
\hline
\hphantom{\bullet}&\hphantom{\bullet}&{\bullet}&{\bullet}&{\bullet}&\hphantom{\bullet}\\
\hphantom{\bullet}&{\bullet}&\hphantom{\bullet}&{\bullet}&\hphantom{\bullet}&\hphantom{\bullet}\\
\hline
\hphantom{\bullet}&\hphantom{\bullet}&\hphantom{\bullet}&\hphantom{\bullet}&\hphantom{\bullet}&{\bullet}\\
{\bullet}&\hphantom{\bullet}&\hphantom{\bullet}&{\bullet}&\hphantom{\bullet}&\hphantom{\bullet}\\
\hline \end{array}_{\domino}
\quad
\begin{array}{|cc|cc|cc|}
\hline
\hphantom{\bullet}&\hphantom{\bullet}&\hphantom{\bullet}&\hphantom{\bullet}&{\bullet}&\hphantom{\bullet}\\
\hphantom{\bullet}&{\bullet}&{\bullet}&\hphantom{\bullet}&\hphantom{\bullet}&\hphantom{\bullet}\\
\hline
{\bullet}&{\bullet}&\hphantom{\bullet}&\hphantom{\bullet}&\hphantom{\bullet}&{\bullet}\\
\hphantom{\bullet}&{\bullet}&\hphantom{\bullet}&{\bullet}&\hphantom{\bullet}&\hphantom{\bullet}\\
\hline \end{array}_{\domino}
\quad
\begin{array}{|cc|cc|cc|}
\hline
\hphantom{\bullet}&\hphantom{\bullet}&{\bullet}&\hphantom{\bullet}&\hphantom{\bullet}&\hphantom{\bullet}\\
\hphantom{\bullet}&{\bullet}&\hphantom{\bullet}&\hphantom{\bullet}&{\bullet}&\hphantom{\bullet}\\
\hline
{\bullet}&{\bullet}&\hphantom{\bullet}&{\bullet}&\hphantom{\bullet}&\hphantom{\bullet}\\
\hphantom{\bullet}&{\bullet}&\hphantom{\bullet}&\hphantom{\bullet}&\hphantom{\bullet}&{\bullet}\\
\hline \end{array}_{\domino}
\quad
\begin{array}{|cc|cc|cc|}
\hline
\hphantom{\bullet}&{\bullet}&{\bullet}&\hphantom{\bullet}&\hphantom{\bullet}&{\bullet}\\
\hphantom{\bullet}&\hphantom{\bullet}&\hphantom{\bullet}&\hphantom{\bullet}&\hphantom{\bullet}&\hphantom{\bullet}\\
\hline
{\bullet}&{\bullet}&\hphantom{\bullet}&\hphantom{\bullet}&\hphantom{\bullet}&\hphantom{\bullet}\\
\hphantom{\bullet}&{\bullet}&\hphantom{\bullet}&{\bullet}&{\bullet}&\hphantom{\bullet}\\
\hline \end{array}_{\domino}
\quad
\begin{array}{|cc|cc|cc|}
\hline
\hphantom{\bullet}&\hphantom{\bullet}&{\bullet}&\hphantom{\bullet}&\hphantom{\bullet}&\hphantom{\bullet}\\
{\bullet}&{\bullet}&\hphantom{\bullet}&{\bullet}&\hphantom{\bullet}&\hphantom{\bullet}\\
\hline
\hphantom{\bullet}&{\bullet}&\hphantom{\bullet}&\hphantom{\bullet}&\hphantom{\bullet}&{\bullet}\\
\hphantom{\bullet}&{\bullet}&\hphantom{\bullet}&\hphantom{\bullet}&{\bullet}&\hphantom{\bullet}\\
\hline \end{array}_{\domino}
$}\nonumber\\[10pt]
&&\hspace*{-1.5em}
\resizebox{0.95\hsize}{!}{$
\begin{array}{|cc|cc|cc|}
\hline
\hphantom{\bullet}&{\bullet}&\hphantom{\bullet}&\hphantom{\bullet}&\hphantom{\bullet}&\hphantom{\bullet}\\
{\bullet}&\hphantom{\bullet}&\hphantom{\bullet}&\hphantom{\bullet}&\hphantom{\bullet}&\hphantom{\bullet}\\
\hline
\hphantom{\bullet}&{\bullet}&{\bullet}&{\bullet}&\hphantom{\bullet}&\hphantom{\bullet}\\
\hphantom{\bullet}&{\bullet}&\hphantom{\bullet}&\hphantom{\bullet}&{\bullet}&{\bullet}\\
\hline \end{array}_{\domino}
\quad
\begin{array}{|cc|cc|cc|}
\hline
{\bullet}&\hphantom{\bullet}&\hphantom{\bullet}&{\bullet}&{\bullet}&\hphantom{\bullet}\\
\hphantom{\bullet}&{\bullet}&\hphantom{\bullet}&{\bullet}&\hphantom{\bullet}&{\bullet}\\
\hline
\hphantom{\bullet}&\hphantom{\bullet}&{\bullet}&{\bullet}&\hphantom{\bullet}&\hphantom{\bullet}\\
\hphantom{\bullet}&\hphantom{\bullet}&\hphantom{\bullet}&\hphantom{\bullet}&\hphantom{\bullet}&\hphantom{\bullet}\\
\hline \end{array}_{\domino}
\quad
\begin{array}{|cc|cc|cc|}
\hline
{\bullet}&{\bullet}&{\bullet}&\hphantom{\bullet}&\hphantom{\bullet}&\hphantom{\bullet}\\
\hphantom{\bullet}&\hphantom{\bullet}&{\bullet}&\hphantom{\bullet}&\hphantom{\bullet}&\hphantom{\bullet}\\
\hline
\hphantom{\bullet}&\hphantom{\bullet}&{\bullet}&\hphantom{\bullet}&\hphantom{\bullet}&\hphantom{\bullet}\\
\hphantom{\bullet}&\hphantom{\bullet}&\hphantom{\bullet}&{\bullet}&{\bullet}&{\bullet}\\
\hline \end{array}_{\domino}
\quad
\begin{array}{|cc|cc|cc|}
\hline
{\bullet}&{\bullet}&{\bullet}&\hphantom{\bullet}&\hphantom{\bullet}&\hphantom{\bullet}\\
\hphantom{\bullet}&\hphantom{\bullet}&\hphantom{\bullet}&{\bullet}&\hphantom{\bullet}&\hphantom{\bullet}\\
\hline
\hphantom{\bullet}&{\bullet}&\hphantom{\bullet}&\hphantom{\bullet}&{\bullet}&\hphantom{\bullet}\\
\hphantom{\bullet}&{\bullet}&\hphantom{\bullet}&\hphantom{\bullet}&\hphantom{\bullet}&{\bullet}\\
\hline \end{array}_{\domino}
\quad
\begin{array}{|cc|cc|cc|}
\hline
{\bullet}&{\bullet}&{\bullet}&\hphantom{\bullet}&\hphantom{\bullet}&\hphantom{\bullet}\\
\hphantom{\bullet}&{\bullet}&\hphantom{\bullet}&\hphantom{\bullet}&{\bullet}&\hphantom{\bullet}\\
\hline
\hphantom{\bullet}&{\bullet}&\hphantom{\bullet}&\hphantom{\bullet}&\hphantom{\bullet}&{\bullet}\\
\hphantom{\bullet}&\hphantom{\bullet}&\hphantom{\bullet}&{\bullet}&\hphantom{\bullet}&\hphantom{\bullet}\\
\hline \end{array}_{\domino}
$}\nonumber\\[10pt]
&&\hspace*{-1.5em}
\resizebox{0.95\hsize}{!}{$
\begin{array}{|cc|cc|cc|}
\hline
\hphantom{\bullet}&\hphantom{\bullet}&\hphantom{\bullet}&{\bullet}&{\bullet}&\hphantom{\bullet}\\
{\bullet}&\hphantom{\bullet}&\hphantom{\bullet}&\hphantom{\bullet}&\hphantom{\bullet}&\hphantom{\bullet}\\
\hline
\hphantom{\bullet}&{\bullet}&\hphantom{\bullet}&\hphantom{\bullet}&\hphantom{\bullet}&\hphantom{\bullet}\\
{\bullet}&{\bullet}&\hphantom{\bullet}&{\bullet}&{\bullet}&\hphantom{\bullet}\\
\hline \end{array}_{\domino}
\quad
\begin{array}{|cc|cc|cc|}
\hline
\hphantom{\bullet}&\hphantom{\bullet}&\hphantom{\bullet}&\hphantom{\bullet}&\hphantom{\bullet}&\hphantom{\bullet}\\
{\bullet}&\hphantom{\bullet}&\hphantom{\bullet}&{\bullet}&\hphantom{\bullet}&{\bullet}\\
\hline
{\bullet}&{\bullet}&\hphantom{\bullet}&\hphantom{\bullet}&\hphantom{\bullet}&\hphantom{\bullet}\\
\hphantom{\bullet}&{\bullet}&\hphantom{\bullet}&{\bullet}&\hphantom{\bullet}&{\bullet}\\
\hline \end{array}_{\domino}
\quad
\begin{array}{|cc|cc|cc|}
\hline
{\bullet}&{\bullet}&\hphantom{\bullet}&\hphantom{\bullet}&\hphantom{\bullet}&\hphantom{\bullet}\\
\hphantom{\bullet}&{\bullet}&\hphantom{\bullet}&{\bullet}&\hphantom{\bullet}&{\bullet}\\
\hline
\hphantom{\bullet}&\hphantom{\bullet}&\hphantom{\bullet}&\hphantom{\bullet}&\hphantom{\bullet}&\hphantom{\bullet}\\
{\bullet}&\hphantom{\bullet}&\hphantom{\bullet}&{\bullet}&\hphantom{\bullet}&{\bullet}\\
\hline \end{array}_{\domino}
\quad
\begin{array}{|cc|cc|cc|}
\hline
{\bullet}&\hphantom{\bullet}&{\bullet}&\hphantom{\bullet}&\hphantom{\bullet}&\hphantom{\bullet}\\
\hphantom{\bullet}&\hphantom{\bullet}&{\bullet}&{\bullet}&\hphantom{\bullet}&{\bullet}\\
\hline
{\bullet}&\hphantom{\bullet}&\hphantom{\bullet}&{\bullet}&\hphantom{\bullet}&\hphantom{\bullet}\\
\hphantom{\bullet}&\hphantom{\bullet}&\hphantom{\bullet}&\hphantom{\bullet}&\hphantom{\bullet}&{\bullet}\\
\hline \end{array}_{\domino}
\quad
\begin{array}{|cc|cc|cc|}
\hline
{\bullet}&\hphantom{\bullet}&\hphantom{\bullet}&\hphantom{\bullet}&{\bullet}&\hphantom{\bullet}\\
\hphantom{\bullet}&\hphantom{\bullet}&\hphantom{\bullet}&{\bullet}&\hphantom{\bullet}&\hphantom{\bullet}\\
\hline
{\bullet}&{\bullet}&\hphantom{\bullet}&{\bullet}&\hphantom{\bullet}&\hphantom{\bullet}\\
\hphantom{\bullet}&{\bullet}&\hphantom{\bullet}&\hphantom{\bullet}&{\bullet}&\hphantom{\bullet}\\
\hline \end{array}_{\domino}
$}\nonumber\\[10pt]
&&\hspace*{-1.5em}
\resizebox{0.95\hsize}{!}{$
\begin{array}{|cc|cc|cc|}
\hline
{\bullet}&\hphantom{\bullet}&\hphantom{\bullet}&{\bullet}&\hphantom{\bullet}&\hphantom{\bullet}\\
\hphantom{\bullet}&\hphantom{\bullet}&\hphantom{\bullet}&\hphantom{\bullet}&{\bullet}&\hphantom{\bullet}\\
\hline
{\bullet}&{\bullet}&\hphantom{\bullet}&\hphantom{\bullet}&{\bullet}&\hphantom{\bullet}\\
\hphantom{\bullet}&{\bullet}&\hphantom{\bullet}&{\bullet}&\hphantom{\bullet}&\hphantom{\bullet}\\
\hline \end{array}_{\domino}
\quad
\begin{array}{|cc|cc|cc|}
\hline
{\bullet}&\hphantom{{\bullet}}&\hphantom{\bullet}&\hphantom{\bullet}&{\bullet}&\hphantom{\bullet}\\
{\bullet}&\hphantom{\bullet}&\hphantom{\bullet}&\hphantom{\bullet}&\hphantom{\bullet}&{\bullet}\\
\hline
\hphantom{\bullet}&\hphantom{\bullet}&{\bullet}&\hphantom{\bullet}&{\bullet}&{\bullet}\\
\hphantom{\bullet}&\hphantom{\bullet}&{\bullet}&\hphantom{\bullet}&\hphantom{\bullet}&\hphantom{\bullet}\\
\hline \end{array}_{\domino}
\quad
\begin{array}{|cc|cc|cc|}
\hline
{\bullet}&{\bullet}&\hphantom{\bullet}&\hphantom{\bullet}&\hphantom{\bullet}&\hphantom{\bullet}\\
{\bullet}&{\bullet}&\hphantom{\bullet}&\hphantom{\bullet}&\hphantom{\bullet}&\hphantom{\bullet}\\
\hline
{\bullet}&{\bullet}&\hphantom{\bullet}&\hphantom{\bullet}&\hphantom{\bullet}&\hphantom{\bullet}\\
{\bullet}&{\bullet}&\hphantom{\bullet}&\hphantom{\bullet}&\hphantom{\bullet}&\hphantom{\bullet}\\
\hline \end{array}_{\domino}
\quad
\begin{array}{|cc|cc|cc|}
\hline
\hphantom{\bullet}&\hphantom{\bullet}&{\bullet}&\hphantom{\bullet}&{\bullet}&\hphantom{\bullet}\\
\hphantom{\bullet}&\hphantom{\bullet}&{\bullet}&\hphantom{\bullet}&\hphantom{\bullet}&{\bullet}\\
\hline
{\bullet}&\hphantom{\bullet}&\hphantom{\bullet}&\hphantom{\bullet}&{\bullet}&{\bullet}\\
{\bullet}&\hphantom{\bullet}&\hphantom{\bullet}&\hphantom{\bullet}&\hphantom{\bullet}&\hphantom{\bullet}\\
\hline \end{array}_{\domino}
\quad
\hphantom{\begin{array}{|cc|cc|cc|}
\hline
\hphantom{\bullet}&\hphantom{\bullet}&{\bullet}&\hphantom{\bullet}&{\bullet}&\hphantom{\bullet}\\
\hphantom{\bullet}&\hphantom{\bullet}&{\bullet}&\hphantom{\bullet}&\hphantom{\bullet}&{\bullet}\\
\hline
{\bullet}&\hphantom{\bullet}&\hphantom{\bullet}&\hphantom{\bullet}&{\bullet}&{\bullet}\\
{\bullet}&\hphantom{\bullet}&\hphantom{\bullet}&\hphantom{\bullet}&\hphantom{\bullet}&\hphantom{\bullet}\\
\hline \end{array}_{\domino}}
$}
\nonumber\\[10pt]
\end{eqnarray}

\section{Non-primitive embeddings of \texorpdfstring{$P(-1)$}{TEXT} into Niemeier lattices}\label{app:npembedding}

Let  $\widehat N$ denote a Niemeier lattice.  
In this appendix we show that there exists a 
non-primitive embedding of $P(-1)$ into $\widehat N$ if and only if $\widehat N$
is of type $E_6^4$.

\begin{figure}
\begin{center}  
\includegraphics[width=7cm,keepaspectratio]{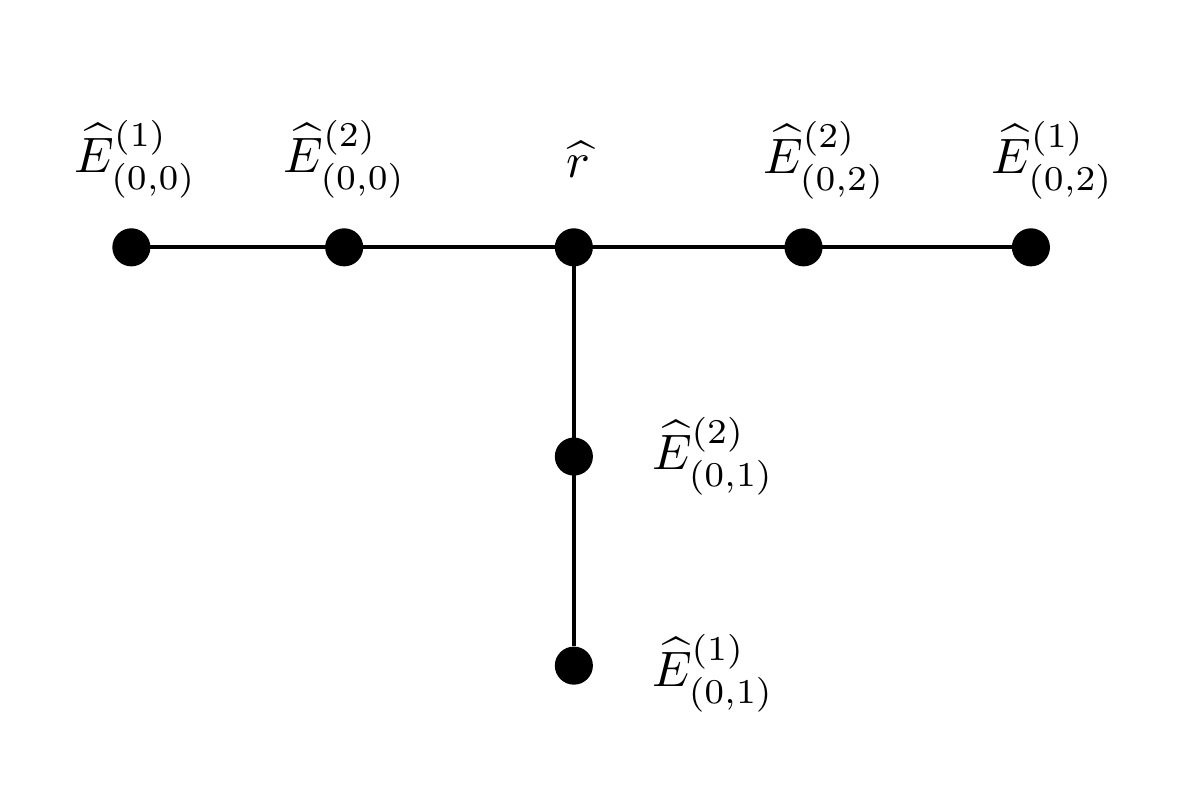}
\end{center}
\caption{{The extended Dynkin diagram of $E_6$ with nodes associated with  the configuration 
of the seven roots $\{\wh{r}, \wh{E}^{(j)}_t, t\in L_{12}, j\in \{1,2\}\}$.}}
\label{fig:extendedE6}
\end{figure}

To see this, assume first that $P(-1)$ is non-primitively embedded into $\widehat N$
by means of $\wh\iota\colon P(-1)\hookrightarrow\wh N$.
Using the same notations as in the proof of proposition \ref{primitiveNiemeierisunique} and
following the arguments given there, for every root $\widehat r\in\widehat N$ 
we either have $\widehat r\in \widehat P^\ast$ or  $\widehat r\perp\widehat P$. 
If  $\wh r\in \wh P$ for all roots $\widehat r\in \widehat P^\ast\cap \wh N$, then using the arguments 
employed in the proof of proposition \ref{primitiveNiemeierisunique}, it follows that $\wh N$ is of type $A_2^{12}$ and that $\wh P$ 
is a primitive sublattice. 
Assume therefore that $\widehat r\in \widehat P^\ast\setminus \widehat P$ is a root
of $\wh N$.
To determine the smallest primitive sublattice of $\wh N$ which contains
$\wh P$, that is, to determine the lattice
$(\widehat P\otimes_\Z\Q)\cap\widehat N=\wh P^\ast\cap\wh N$,
by proposition \ref{Pconstruct}, without loss of generality we may assume that 
$\widehat r=-\frac{1}{3}\sum\limits_{t\in L_{12}}\widehat E_t$. One now checks by a direct calculation that the 
Gram matrix for the roots 
$\widehat E_t^{(j)}$ with $t\in L_{12}$ and $j\in\{1,2\}$ together with $\widehat r$ agrees
with the adjacency matrix of the extended Dynkin diagram for $E_6$ (see figure~\ref{fig:extendedE6}).
In other words, these seven roots, which are 
linearly dependent by construction, generate a root lattice of type $E_6$.

Note that by proposition \ref{Pconstruct}, with $\wh P\subset\wh N$ and
$\widehat r\in\widehat N$ we also have 
$\widehat r^\prime:=-\frac{1}{3}\sum\limits_{t\in L^\prime}\widehat E_t\in\widehat N$ for 
any affine parallel
$L^\prime$  to $L_{12}$ in $\F_3^2$. By the same arguments as above, the roots 
$\widehat E_t^{(j)}$ with $t\in L^\prime$ and $j\in\{1,2\}$ together with $\widehat r^\prime$ 
generate a root lattice of type $E_6$.

We claim that the above implies that 
$(\widehat P\otimes_\Z\Q)\cap\widehat N=\widehat P^\ast\cap\widehat N$ is generated
by $\wh P$ and $\wh r$. Indeed, 
according to  proposition \ref{Pconstruct} the lattice $\wh P^\ast$
is generated by $\wh P$ together with the images $\wh p_1,\, \wh p_2,\, \wh p_3$
of $p_1,\,p_2,\, p_3\in P^\ast$ under $\wh\iota$, or equivalently by 
$\wh P$ together with $\wh r,\, \wh p_1,\, \wh p_2-\wh p_3$.
Hence we need to show that 
$a_1\widehat p_1+a_2(\wh p_2-\widehat p_3)\in\widehat N$ with $a_1,a_2\in\{-1,0,1\}$ 
implies $a_1=a_2=0$.
The latter follows from the fact that as an element of $\wh N$, the vector
$a_1\widehat p_1+a_2(\wh p_2-\widehat p_3)$ has integral
scalar product with $\widehat r$ and even squared length. This proves our claim 
about the lattice $\widehat P^\ast\cap\widehat N$, and it implies that the root 
sublattice of this lattice is of type $E_6^3$.

From what we have  shown so far,
we may  confirm that the root lattice of $\widehat N$ is of type $E_6^3\oplus \widehat K$, 
where $\widehat K$ is a root lattice of rank $6$. The list of root lattices for Niemeier lattices,
stated e.g.~in \cite[Table 16.1]{conway1998sphere}, reveals that $\widehat N$ is of type $E_6^4$.

It remains to show that  there does exist a non-primitive embedding of $P(-1)$
in the Niemeier lattice $\widehat N$ of type $E_6^4$. To see this, note that the 
above already shows that there exists a root sublattice  of type $E_6^3$ in $P^\ast(-1)$, 
which we may embed in $\widehat N$. Denote its image in 
$\widehat N$ by $\widehat R$ and extend the map $\Q$-linearly to obtain an embedding
of $P(-1)$ in $\wh N\otimes_\Z\Q$. We need to show that the image $\wh P$ of $P(-1)$
under this embedding is contained in $\wh N$.
Since $P(-1)$ does not contain a root lattice of 
type $E_6^3$, this will also imply that the embedding of $P(-1)$ is not primitive.

We use the same notations as above for the 
images $\widehat E_t$, $\widehat r$ and $\widehat r^\prime$
of the respective lattice vectors
in $\widehat N$. We know that $\wh N$ is obtained from its root sublattice
of type $E_6^4$ by including three additional generators, and 
using the explicit form of these generators according to \cite[Table 16.1]{conway1998sphere}
one shows that 
$\widehat R^\ast\cap\widehat N$ is obtained by adding one further 
generator to $\widehat R$. Concretely, for the lattice of type $E_6$ 
which is generated
by $\widehat E_t^{(j)}$ with $t\in L_{12}$ and $j\in\{1,2\}$ together with 
$\widehat r=-\frac{1}{3}\sum\limits_{t\in L_{12}}\widehat E_t$,
the discriminant group is isomorphic to $\Z_3$, and a generator of this group
is represented by $\check E_0:=\frac{1}{3}\left( \widehat E_{00}-\widehat E_{01}\right)$.
Similarly for the two other $E_6$-type summands of $\widehat R$, the discriminant
groups are generated by the classes represented by 
$\check E_1:=\frac{1}{3}\left( \widehat E_{10}-\widehat E_{11}\right)$
and $\check E_2:=\frac{1}{3}\left( \widehat E_{20}-\widehat E_{21}\right)$.
Now $\widehat R^\ast\cap\widehat N$ is obtained from the lattice $\widehat R$
by including the extra generator $\check E_0+\check E_1+\check E_2$.
To confirm that $\widehat P\subset \widehat R^\ast\cap\widehat N$, recall 
from \eqref{eq:v1v2v3} that $\widehat P$
is generated by the roots $\widehat E_t^{(j)}$ with $t\in\F_3^2$, all of which
are contained in $\widehat R\subset \widehat R^\ast\cap\widehat N$, and 
the images $\widehat v_j$ of the vectors
$v_j$, $j\in\{1,2,3\}$. Hence with
$\widehat r^\prime = -\frac{1}{3}\sum\limits_{t\in (1,0)+L_{12}} \widehat E_t$ 
and $\widehat r^{\prime\prime}= -\frac{1}{3}\sum\limits_{t\in (2,0)+L_{12}} \widehat E_t$,
the claim
$\widehat P\subset \widehat R^\ast\cap\widehat N$ follows from 
$$
\begin{array}{rccclcl}
\widehat v_1&=&\textstyle
\frac{1}{3}\left(\smash{\sum\limits_{t\in L_{12}} \widehat E_t-
\sum\limits_{t\in (1,0)+L_{12}} }\widehat E_t\right)
&=&\widehat r^\prime-\widehat r &\in& \widehat R\,,\\[1em]
\widehat v_2&=&\textstyle
\frac{1}{3}\left(\smash{\sum\limits_{t\in L_{34}} \widehat E_t-
\sum\limits_{t\in (0,1)+L_{34}} }\widehat E_t\right)
&=&
\check E_0+\check E_1+\check E_2 &\in&\widehat R^\ast\cap\widehat N\,,\\[1em]
\widehat v_3
&=&
\textstyle\frac{1}{3}\sum\limits_{t\in \F_3^2} \widehat E_t
&=&-\left(\widehat r+\widehat r^\prime+\widehat r^{\prime\prime}\right)&\in&\widehat R\,.
\end{array}
$$ 

\clearpage

\addcontentsline{toc}{section}{References}
\bibliographystyle{JHEP}
\bibliography{mono}

\end{document}